\documentclass[12pt]{article}
\usepackage{latexsym, amsbsy, amssymb,multirow, epsfig, amsmath}
\usepackage{hyperref}

\usepackage{graphicx,setspace,lscape,longtable}
\usepackage{epsfig,graphicx}

\usepackage{amsmath,amsthm,amssymb,color}

\RequirePackage[mathlines, displaymath]{lineno}

\def\beqr{\begin{eqnarray}}
\def\eeqr{\end{eqnarray}}
\def\beqrs{\begin{eqnarray*}}
\def\eeqrs{\end{eqnarray*}}
\def\bep{\begin{prop}}
\def\eep{\end{prop}}

\numberwithin{equation}{section}

\def\ba{\mbox{\boldmath$\alpha$}}
\def\bb{\mbox{\boldmath$\beta$}}
\def\bg{\mbox{\boldmath$\gamma$}}

\newcommand{\bSig}{\mbox{\boldmath $\Sigma$}}

\newcommand{\bLam}{\mbox{\boldmath $\Lambda$}}

\newtheorem{theo}{\bf Theorem}[section]

\newtheorem{prop}{\bf Proposition}[section]

\newtheorem{lemma}{Lemma}

\newtheorem{remark}{Remark}

\usepackage{latexsym}
\usepackage{epsfig}
\usepackage{bm}

\usepackage{threeparttable}
\usepackage{graphicx}
\usepackage[small]{caption2}
\usepackage{threeparttable}
\usepackage{dcolumn}
\usepackage{multirow}
\usepackage{booktabs}
\usepackage{footnote}
\usepackage[linesnumbered,boxed]{algorithm2e}
\usepackage[final]{pdfpages}

\DeclareSymbolFont{largesymbol}{OMX}{yhex}{m}{n}
\DeclareMathAccent{\Widehat}{\mathord}{largesymbol}{"62}

\begin{document}

\begin{center}
{\bf\large Model Checking for Parametric Ordinary Differential Equations System}
\\ \vskip0.6cm
Ran Liu$^1$, Yun Fang$^2$ and Lixing Zhu$^{1}$\footnote{The research was supported by a grant from the University Grants Council of Hong Kong and a SNFC grant. } 
\\
\textit{ $^1$ Hong Kong Baptist University\\
$^2$ Shanghai Normal University}

\today
\end{center}

\begin{singlespace}
\begin{abstract}
	Ordinary differential equations have been used to model dynamical systems in a broad range. Model checking for parametric ordinary differential equations is  a necessary step to check whether the assumed models are plausible. In this paper  we introduce three test statistics for their different purposes. We first give a trajectory matching-based test for the whole system. To further identify which component function(s) would be wrongly modelled,  we introduce two test statistics that are  based on integral matching and gradient matching respectively. We investigate the asymptotic properties of the three test statistics under the null, global and local alternative hypothesis. To achieve these purposes, we also investigate the asymptotic properties of nonlinear least squares estimation and two-step collocation estimation under both the null and alternatives. The results about the estimations are also new in the literature. To examine the performances of the tests, we conduct several numerical simulations. A real data example about immune cell kinetics and trafficking for influenza infection is analyzed for illustration.
\end{abstract}


\end{singlespace}

\newpage
\section{Introduction}
As they can model how the systems evolve with time, ordinary differential equations (ODEs) have been widely applied in many scientific fields such as physics, ecology (\cite{Lotka.1910}; \cite{Volterra.1928}; \cite{GOEL.1971}) and neuroscience (\cite{FitzHugh.1961}; \cite{Nagumo.1962}). A system of ODEs can be written as
\begin{equation}
X^{\prime}(t) \equiv \left[ \begin{array}{c}{\frac{d X_{1}(t)}{d t}} \\ {\vdots} \\ {\frac{d X_{p}(t)}{d t}}\end{array}\right]=\left[ \begin{array}{c}{g_{1}(t)} \\ {\vdots} \\ {g_{p}(t)}\end{array}\right] \\ \equiv g(t),
\end{equation}
where $X(t)=\left(X_{1}(t), \ldots, X_{p}(t)\right)^{\mathrm{T}}$ is a $p$-dimensional state vector. Typically, this system is measured on discrete time points with noises, say
\begin{equation}
Y_{i}=X\left(t_{i}\right)+\epsilon_{i}, \quad i=1, \ldots, n
\end{equation}
where the measurement errors $\epsilon_{i}$ satisfying $E(\epsilon_{i}|t_i)=0$ have nonsingular variance-covariance matrix $\Sigma_{\epsilon_i}$, and are independent with $\epsilon_{j}$ for every $j \neq i$.

In a large number of scientific questions, the vector of functions $g=\left(g_{1}, \ldots, g_{p}\right)^{\mathrm{T}}$ is supposed to belong to a given parametric family of functions $\mathcal{F}=\{f(\cdot, \theta)=f(t, X(t ; \theta); \theta)=\left(f_{1}(t, X(t ; \theta); \theta), \ldots, f_{p}(t, X(t ; \theta); \theta)\right)^{\mathrm{T}}: \theta \in \Theta \subset R^q \}$. Since the vector of parameters $\theta=\left(\theta_{1}, \ldots, \theta_{q}\right)^{T}$ is unknown, many efforts have been devoted to  parameter estimation and further statistical analysis. When the assumption on the parametric form does not hold, i.e, $g$ does not belong to $\mathcal{F}$, any further statistical analysis would be unreliable. Thus, a model checking for the assumed ODEs should be accompanied.
There is no method available in the literature about such a hypothesis testing problem. In this paper, we construct tests to fill up this gap.

To discuss how to check ODE models, we first review some relevant  methodologies of model checking for regressions in the literature to see whether those methods can motivate  test constructions we need. There exist two broad classes of tests. Tests in a class use nonparametric estimations, thus are called the local smoothing tests. Those tests include 
\cite{Hardle.1993}, \cite{Zheng.1996}, \cite{Lixing.1998}, \cite{Dette.1999} and \cite{Koul.2004} as examples. In general, the tests in this class  are sensitive to alternative models that are oscillating/ highly frequent. Tests in another class are based on residual-marked empirical processes and take averages over an index set. As averaging itself is a global smoothing step and thus, they are called the global smoothing tests such as \cite{Stute.1997}, \cite{Stute.1998a}, \cite{Stute.1998b} and \cite{Lixing.2003}. The tests in this class can have better asymptotic properties, but less sensitive to oscillating alternative models. \cite{GONZALEZMANTEIGA.2013} is a comprehensive reference.  \cite{Tan.2019} constructed a global smoothing test when the number of predictors is divergent as the sample size goes to infinity. 
The testing problem investigated in this paper is however rather different from the problems for regressions as the parametric model structure is not directly  on the unknown function $X(\cdot)$, but its derivative $X'(\cdot)$. Further, any component of $X'(\cdot)$ is also  related to the original function $X(\cdot)$. This structure then causes the problem much more complicated than the problems for regressions. We will discuss these issues in the next four sections.

To construct test statistics, we indispensably  need the parameter estimation of ODE models under both the null and alternative hypothesis. There are two commonly used types of methods: nonlinear least squares method and two-step collocation method (see, e.g. \cite{Chen.2017}; \cite{Ramsay.2017}). The nonlinear least squares method estimates parameters by matching the trajectory of ODEs. If we do not consider the numerical error of numerical solution, the classical nonlinear least squares theory is applicable and the corresponding estimator $\hat{\theta}_{NLS}$ is $\sqrt{n}$-consistent under certain regularity conditions. When  both numerical error and measurement error are involved, \cite{Xue.2010} pointed out that if the maximum step size of the $l$-order numerical algorithm for integration computation goes to zero at a rate faster than $n^{-1} /(l \wedge 4)$, the numerical error is negligible. We will use this estimator in constructing a trajectory matching-based test ($TM_{n}$) to check the whole ODE system.

Another testing problem is more challenging. That is, we wish to identify which component in the ODE system may not be correctly modeled. The main challenge is that actually any single component involves the same original function $X(\cdot)$. Any departure from the hypothetical model of $X$ could affect all components, but it is unclear how and at what degree the impact from the departure of $X$ is for the each component. We try to construct two tests to handle this issue. To this end,
nonlinear least squares estimation is no longer feasible as it involves the integral for all components, which cannot directly focus on the component we are going to check for. We then use  two-step collocation method to smooth data in the first stage and estimate parameters, in the second stage, by matching the gradient or the indefinite integral of ODEs. The $\sqrt{n}$-consistency of $\hat{\theta}_{TS}$ also holds under certain regularity conditions, but with less estimation accuracy (\cite{Brunel.2008}; \cite{Liang.2008}; \cite{Gugushvili.2012}; \cite{Dattner.2015}; \cite{Ramsay.2017}). We will study the asymptotic property of this estimator under different hypotheses and use it in constructing an integral matching-based test ($IM_{n}$) and a gradient matching-based test ($GM_{n}$) for, particularly, checking every component function of ODEs.

As the by-products, we investigate the asymptotic properties of the estimations, particularly  two-step collocation estimation under both the null and alternative hypothesis because the asymptotics of nonlinear least squares estimation can be similarly derived from existing results for existing estimation for regressions. These results are also new in the literature.

For test constructions, we propose an idea to solve the ODEs analytically or numerically and to convert them to multi-response regression models. Then  test statistics can be constructed,  grounded on the classical methods for  model checking, by using the residuals between the observed data and responses.
However, the significant difference and difficulty for ODE models from the ordinary multi-response regression models come from that any response of ODE-based multi-response models is also a function of other responses and thus the residuals are very complicated in function form. Thus the trajectory matching-based test can detect the alternatives distinct from the whole system under the null hypothesis, but cannot check which component(s) is (are) would not be tenable. We will discuss this phenomenon in detail in Section~3. To test the null hypothesis for  each component of ODE models, we construct integral matching-based test and a gradient matching-based test. These tests use  the two-step collocation methods for parameter estimation. More specifically, in the first step we estimate $X(t)$ and $X^{\prime}(t)$ non-parametrically, which decouples the connections among different components. Then we compute two types of pseudo-residuals connected to integral matching and gradient matching. The test statistics are the functionals of these two pseudo-residuals respectively. We also need to point out why we prefer to use such a local smoothing idea to construct tests rather than global smoothing method (see, e.g. \cite{Stute.1998a} and \cite{Stute.1998b}). This is because some of very useful ODE models are highly oscillating.

The rest of the paper is organized as follows. Section 2 contains the construction of $TM_{n}$ and the asymptotic properties under the null and alternative hypothesis. In Section 3 we will talk about some particularities
of checking ODE models and give the ideas to overcome the difficulties induced by these particularities. The construction of $IM_{n}$ and relevant results are discussed in Section 4. The construction and properties of $GM_{n}$ are  presented in Section 5. In Section 6, we report simulation results of the proposed tests and the analysis for a real data example concerning a model of influenza-specific CD8+ T cells (\cite{Wu.2011}; \cite{Ding.2014}). Section 7 contains a summary of the study and a brief discussion for further research. As the technical proofs are very tedious, we then put them in the supplementary materials. 

\section{Trajectory matching-based test}
\subsection{The hypotheses and test statistic}
Recall that $\mathcal{F}=\{f(\cdot, \theta)=f(t, X(t ; \theta); \theta)=\left(f_{1}(t, X(t ; \theta); \theta), \ldots, f_{p}(t, X(t ; \theta); \theta)\right)^{\mathrm{T}}: \theta \in \Theta \subset R^q \}$ is a given parametric family of functions. The hypotheses are as follows:
\begin{eqnarray*}	
	&H_0: & X^{\prime}(t) \equiv \left[ \begin{array}{c}{\frac{d X_{1}(t)}{d t}} \\ {\vdots} \\ {\frac{d X_{p}(t)}{d t}}\end{array}\right]=\left[ \begin{array}{c}{g_{1}(t)} \\ {\vdots} \\ {g_{p}(t)}\end{array}\right] \equiv g(t) = f(t, X(t ; \theta_0); \theta_0)  \in \mathcal{F},\\	
	&H_1: & X^{\prime}(t) = g(t) \notin \mathcal{F},
\end{eqnarray*}
where $\theta_0$ is an unknown parameter vector. Here we use $X(t ; \theta_0)$ to present $X(t)$ under the null hypothesis.

According to Cauchy - Lipschitz theorem, the equation
\begin{equation}	
X^{\prime}(t)=f(t, X(t ; \theta); \theta),
\end{equation}
has a unique solution $X(t; \theta)=F(t; \theta)$ under mild regularity conditions. Therefore, if we solve the equation analytically or numerically, the problem of checking ODE models is converted to the problem of testing whether $X(t)=F(t; \theta_0)$ for some $\theta_0 \in \Theta \subset R^q$.

We then transfer the ODE models to a multi-response regression when we have the observations $Y$, $t$ and the function form of $X' (\cdot)$ up to some unknown parameters under the null hypothesis. Consider the $p=1$ case to motivate our construction. 
Let  $\varepsilon_{i} \equiv Y_i-F(t_i; \theta^*)$ with $\theta^*=\arg \min_{\theta} E[\left\|Y_i-F(t_i;\theta)\right\|^2]$ be the residual. Note that under $H_0$, $\varepsilon_{i}=\epsilon_{i}$ and $E(\varepsilon_{i}|t_i)=0$ leads to $E[\varepsilon_{i}E(\varepsilon_{i}|t_i)p(t_i)]=0$, while under $H_1$,  $E(\varepsilon_{i}|t_i)=X(t_i)-F(t_i; \theta)\neq 0$, and $E[\varepsilon_{i}E(\varepsilon_{i}|t_i)p(t_i)]=E\{[E(\varepsilon_{i}|t_i)]^2p(t_i)\}>0$. Thus, letting  $e_i \equiv Y_i-F(t_i;\hat{\theta})$ be an estimator of $\varepsilon_i$, we can use the sample analogue of $E[\varepsilon_{i}E(\varepsilon_{i}|t_i)p(t_i)]$ to construct a test statistic
\begin{equation}	
V_{n}^{Zh} \equiv \frac{1}{n(n-1)} \sum_{i=1}^{n} \sum_{j=1 \atop j \neq i}^{n} \frac{1}{h} K\left(\frac{t_{i}-t_{j}}{h}\right) e_{i} e_{j},
\end{equation}
where $K$ is a kernel function, $h$ is a bandwidth parameter. This is in spirit similar to the test suggested by \cite{Zheng.1996}. A standardized test statistic $T_n^{Zh}$ can be easily obtained by using $V_n^{Zh}$ and its variance.

In the multi-response case, for every component $i$, we can construct a test statistic $V_{ni}^F$ to check whether $X_i(t)=F_i(t; \theta_0)$ for some $\theta_0 \in \Theta \subset R^q$ using the above idea. Therefore, we obtain a vector version of $V_n^{Zh}$, which is expressed as $V_n^F=\left(V_{n1}^F, \ldots, V_{np}^F\right)^{\mathrm{T}}$. To summarize the information contained in $V_n$, we aggregate $V_n$ to make a test statistic and write it as $TM_{n}$ in short:
\begin{equation}
TM_{n} \equiv n^2 h V_{n}^{FT}\widehat{\Sigma}^{F-1}V_{n}^F.
\end{equation}
Here $\widehat{\Sigma}^F$ is a symmetric matrix to normalize the test statistic:
\begin{equation}
\widehat{\Sigma}^F=\frac{2}{n(n-1)} \sum_{i=1}^{n} \sum_{j=1 \atop j \neq i}^{n} \frac{1}{h} K^{2}\left(\frac{t_{i}-t_{j}}{h}\right) (e_{i}\odot e_{j})(e_{i}\odot e_{j})^{\mathrm{T}}
\end{equation}
where $\odot$ denotes component-wise multiplication.

Let $\| \cdot \|$ represent the Frobenius norm. The unknown parameter $\theta$ is estimated using nonlinear least squares method:
\begin{equation}
\begin{aligned} \hat{\theta}_{NLS} &=\underset{\theta}{\arg \min } \sum_{i=1}^{n}\left\|Y_{i}-X\left(t_{i} ; \theta\right)\right\|^{2} \\ & \text {subject to } X^{\prime}(t ; \theta)=f(t,X(t ; \theta); \theta),\quad t \in[t_0,T] \end{aligned}
\end{equation}
where the trajectory $X(t ; \theta)$ is obtained by numerical methods such as Euler backward method and 4-stage Runge-Kutta algorithm when there is no closed-form solution. The following theorem gives the asymptotic properties of nonlinear least squares estimator under the null, global alternative and local alternative hypothesis, which is needed for studying the asymptotic results of $TM_{n}$. The results are of independent interest in estimation theory.
\begin{theo}
	Given  sets A and B of assumptions in Supplement A,  supposing the numerical error of numerical solution is negligible, the nonlinear least squares estimator $\hat{\theta}_{NLS}$ for ODE is a consistent estimator of $\theta^*_{NLS}$ with $\theta^*_{NLS}=\arg \min_{\theta \in \Theta} \mathrm{E}\left[ \left\|Y(t)-F\left(t ; \theta\right)\right\|^{2} \right] $. Further,  we have the following decomposition.
	
	1. Under the null hypothesis, we have $\theta^*_{NLS}=\theta_0$ and
	\begin{equation}
	\begin{aligned}
	\sqrt{n}(\hat{\theta}_{NLS}-\theta_{0})&=H_{\dot{F}}^{-1}\frac{1}{\sqrt{n}} \sum_{i=1}^{n}\sum_{k=1}^{p}\left[ \epsilon_{ik}\frac{\partial F_k\left(t_i;\theta_0 \right)}{\partial \theta}\right]+\mathrm{o}_{P}(1)\\
	\end{aligned}
	\end{equation}
	where
	\begin{equation}
	H_{\dot{F}}=\mathrm{E}\left[ \sum_{k=1}^{p}\frac{\partial F_k\left(t;\theta_{0}\right)}{\partial \theta} \frac{\partial F_k\left(t; \theta_{0}\right)}{\partial \theta^{T}} \right].
	\end{equation}
	
	2. Under the global alternative hypothesis $H_1$, we have $\theta^*_{NLS}=\theta_1^*$ with
	\begin{equation}
	\theta_1^*=\arg \min_{\theta \in \Theta} \mathrm{E}\left[ \left\|X(t)-F\left(t ; \theta\right)\right\|^{2} \right],
	\end{equation}
	and
	\begin{equation}
	\begin{aligned}
	&\sqrt{n}(\hat{\theta}_{NLS}-\theta_{1}^*)\\
	=&G_{\dot{F}}^{-1}\frac{1}{\sqrt{n}} \sum_{i=1}^{n} \sum_{k=1}^{p} \left\{ \left[Y_{ik}-F_k\left(t_i;\theta_1^* \right)\right]\frac{\partial F_k\left(t_i;\theta_1^* \right)}{\partial \theta}\right\}+\mathrm{o}_{P}(1)\\
	\end{aligned}
	\end{equation}
	where
	\begin{equation}
	\begin{aligned}
	G_{\dot{F}}=&\mathrm{E}\left[ \sum_{k=1}^{p}\frac{\partial F_k\left(t; \theta_{1}^*\right)}{\partial \theta} \frac{\partial F_k\left(t; \theta_{1}^*\right)}{\partial \theta^{T}} \right]\\
	&-\mathrm{E}\left\{ \sum_{k=1}^{p} \left[ X_k(t)-F_k\left(t ; \theta_1^*\right) \right] \frac{\partial^{2} F_k\left(t; \theta_1^* \right)}{\partial \theta \partial \theta^{\mathrm{T}}}  \right\}.\\
	\end{aligned}
	\end{equation}
	
	3. Consider a sequence of local alternatives via adding local disturbance to the trajectory of the ODE model:
	\begin{equation}\label{local}
	H_{1n}^F:X(t) = F\left(t; \theta_{0}\right)+\delta_{n} L\left(t\right)
	\end{equation}
	where $L(t)=(L_{1}(t),\ldots,L_{p}(t))^{T}$ is a bounded multiple response function, and $\delta_n \rightarrow 0$ as $n \rightarrow \infty$.
	
	Under this local alternative hypothesis, we have $\theta^*_{NLS}=\theta_0$ and
	\begin{equation}
	\begin{aligned}
	\sqrt{n}(\hat{\theta}_{NLS}-\theta_{0})=&H_{\dot{F}}^{-1}\frac{1}{\sqrt{n}} \sum_{i=1}^{n}\sum_{k=1}^{p}\left[ \epsilon_{ik}\frac{\partial F_k\left(t_i;\theta_0 \right)}{\partial \theta}\right]\\
	&+\sqrt{n}\delta_{n}H_{\dot{F}}^{-1}\mathrm{E}\left[ \sum_{k=1}^{p}L_k(t)\frac{\partial F_k\left(t;\theta_0 \right)}{\partial \theta} \right]+\mathrm{o}_{P}(1).\\
	\end{aligned}
	\end{equation}\label{NLS}
\end{theo}

Note that in  test construction, the estimator $\hat \theta_{NLS}$ involves all components of $X'(\cdot)$. Further, as the test is based on the whole original function $X(\cdot)$ to involve all components. This will cause the test only for the whole system of ODE's. We will give more detailed discussion in Section~3.

\subsection{Asymptotic properties}

It is easy to see that $V_n^F$ can be asymptotically written as  a vector of  U-statistics. 
We first state the asymptotic properties of $V_n$ under the null hypothesis.
\begin{lemma}
	Given sets A and B of assumptions in Supplement A. If $h \rightarrow 0$ and $n h \rightarrow \infty$, then under the null hypothesis,
	\begin{equation}
	n h^{1 / 2} V_{n}^F \stackrel{\mathbf{d}}{\longrightarrow} \mathrm{N}(0, \Sigma^F)
	\end{equation}
	where $\Sigma^F$ is a symmetric matrix with the entries: for any pair $(k_1, k_2)$ with $1\le k_1, k_2 \le p$,
	\begin{equation}
	\Sigma^F_{k_1k_2}=2\int K^{2}(u) \mathrm{d} u \cdot \int\left(\sigma_{k1k2}(t)\right)^{2} p^{2}(t) \mathrm{d} t .
	\end{equation}\label{LT.1}
\end{lemma}
Since $\Sigma^F$ is unknown, we use an estimator $\widehat{\Sigma}^F$ to replace it. The following lemma gives the consistency of this estimator.

\begin{lemma}
	Given sets A and B of assumptions in Supplement A. If $h \rightarrow 0$ and $n h \rightarrow \infty$, then under the null hypothesis,
	\begin{equation}
	\widehat{\Sigma}^F \stackrel{\mathbf{P}}{\longrightarrow} \Sigma^F
	\end{equation}
	where $\widehat{\Sigma}^F$ is a symmetric matrix as
	\begin{equation}
	\widehat{\Sigma}^F=\frac{2}{n(n-1)} \sum_{i=1}^{n} \sum_{j=1 \atop j \neq i}^{n} \frac{1}{h} K^{2}\left(\frac{t_{i}-t_{j}}{h}\right) (e_{i}\odot e_{j})(e_{i}\odot e_{j})^{\mathrm{T}}
	\end{equation}\label{LT.2}
\end{lemma}

Having Lemma~\ref{LT.1} and Lemma~\ref{LT.2}, it is easy to show the asymptotic property of $TM_{n}$ under the null hypothesis.
\begin{theo}
	Given sets A and B of assumptions in Supplement A. If $h \rightarrow 0$ and $n h \rightarrow \infty$, then under the null hypothesis,
	\begin{equation}
	TM_{n} \stackrel{\mathbf{d}}{\longrightarrow} \chi_p^2
	\end{equation}
	where $TM_{n} = n^2 h V_{n}^{FT}\widehat{\Sigma}^{F-1}V_{n}^F$.\label{TM_0}
\end{theo}

Theorem~\ref{TM_0} shows that the test statistic is asymptotically chi-square distributed with $p$ degrees of freedom and thus  the critical values can be easily determined by the limiting null distribution.

Next, we study the asymptotic power of the test under the global alternative. Hereafter we use the notation $v^2 \equiv v \odot v$ for the vector $v$. The following two lemmas give the asymptotic properties of $V_{n}^F$ and $\widehat{\Sigma}^F$.

\begin{lemma}
	Given sets A and B of assumptions in Supplement A. If $h \rightarrow 0$ and $n h \rightarrow \infty$, then under the global alternative $H_1$,
	\begin{equation}
	V_{n}^F \stackrel{\mathbf{P}}{\longrightarrow} E\left\{\left[X\left(t_{i}\right)-F\left(t_{i}, \theta_{1}\right)\right]^{2} \odot p\left(t_{i}\right)\right\} > 0
	\end{equation}
	where $\left[X\left(t_{i}\right)-F\left(t_{i}, \theta_{1}\right)\right]^{2}=\left[X\left(t_{i}\right)-F\left(t_{i}, \theta_{1}\right)\right] \odot \left[X\left(t_{i}\right)-F\left(t_{i}, \theta_{1}\right)\right]$.\label{LT.3}
\end{lemma}

\begin{lemma}
	Given Assumptions A and B in Supplement A. If $h \rightarrow 0$ and $n h \rightarrow \infty$, then under $H_1$,
	\begin{equation}
	\widehat{\Sigma}^F \stackrel{\mathbf{P}}{\longrightarrow} \Sigma^{F\prime} > 0,
	\end{equation}
	here $\Sigma^{F\prime}$ is defined as follows: for any element $(k_1,k_2)$ with $1\le k_1, k_2\le p$,
	\begin{equation}
	\begin{aligned}
	\Sigma^{F\prime}_{k_1k_2}=&2\int K^{2}(u) \mathrm{d} u \\
	&\times \int \left\{ \sigma_{k1k2}(t)+\left[X_{k_1}\left(t\right)-F_{k_1}\left(t, \theta_{1}\right)\right]\left[X_{k_2}\left(t\right)-F_{k_2}\left(t, \theta_{1}\right)\right] \right\}^{2} p^{2}(t) \mathrm{d} t .
	\end{aligned}
	\end{equation}\label{LT.4}
\end{lemma}

Therefore we have the following theorem to state the asymptotic property of $TM_{n}$ under $H_1$.
\begin{theo}
	Given sets A and B of assumptions in Supplement A. If $h \rightarrow 0$ and $n h \rightarrow \infty$, then under $H_1$,
	\begin{equation}
	TM_{n}/(n^2h) \stackrel{\mathbf{P}}{\rightarrow} V^{\prime T}\Sigma^{F\prime -1}V^{\prime}
	\end{equation}
	where $V^{\prime}=E\left\{\left[X\left(t_{i}\right)-F\left(t_{i}, \theta_{1}\right)\right]^{2} \odot p\left(t_{i}\right)\right\}$.\label{TM_1}
\end{theo}

Theorem~\ref{TM_1} shows that this test is consistent and sensitive to the global alternative in the sense that it can diverge to infinity at a very fast rate of order $n^2h$.

We now consider the power performance of the test under local alternatives of (\ref{local}):
\begin{equation*}	
H_{1n}^F:X(t) = F\left(t, \theta_{0}\right)+\delta_{n} L\left(t\right)
\end{equation*}
with the bounded multiple response function $L(t)$ and the $\mathrm{o}(1)$ term $\delta_n$. The following theorem states the the asymptotic property of $TM_{n}$ under $H_{1n}^F$.
\begin{theo}
	Given sets A and B of assumptions in Supplement A. If $h \rightarrow 0$ and $n h \rightarrow \infty$, then under  $H_{1n}$ with $n^{1 / 2} h^{1 / 4}\delta_{n} \rightarrow \infty$,
	\begin{equation}
	TM_{n}/(n^2h\delta_{n}^{4}) \stackrel{\mathbf{P}}{\longrightarrow} \mu^{\mathrm{T}}\Sigma^{F -1}\mu
	\end{equation}
	where $\mu$ is a $p$-dimensional vector with the $i$-th element $$\mu_i=\mathrm{E}\left\{\left[ L_i\left(t\right)-\frac{\partial F_i\left(t, \theta_{0}\right)}{\partial \theta^{\prime}} H^{-1}_{\dot{F}}\mathrm{E}\left(\sum_{k=1}^{p}L_k\left(t\right)\frac{\partial F_k\left(t, \theta_{0}\right)}{\partial \theta}\right)  \right]^2 p(x) \right\}^{\mathrm{T}},$$
	and $H_{\dot{F}}=\mathrm{E}\left(\sum_{k=1}^{p}\frac{\partial F_k\left(t, \theta_{0}\right)}{\partial \theta} \frac{\partial F_k\left(t, \theta_{0}\right)}{\partial \theta^{T}} \right)$.
	
	Particularly, if $\delta_{n}=n^{-1 / 2} h^{-1 / 4}$,
	\begin{equation}
	TM_{n} \stackrel{\mathbf{d}}{\longrightarrow} \chi_p^2(\lambda)
	\end{equation}
	where $\chi_p^2(\lambda)$ is noncentral chi-squared distribution with the noncentrality parameter  $\lambda=\mu^{\mathrm{T}}\Sigma^{F -1}\mu$.\label{TM_1n}
\end{theo}

This result shows that the test can detect the local alternatives distinct from the null at the  rate of order $n^{-1 / 2} h^{-1 / 4}$. This is the typical rate  of local smoothing tests for classical regressions in the literature, see \cite{Zheng.1996}.

\begin{remark}
	In the case that $p=1$, $V_{n}^F$ is similar to Zheng's statistic $V_n$ and it seems it follows the results of Theorem~3 of \cite{Zheng.1996}. However, the proof of Theorem~3 of \cite{Zheng.1996} needs a further condition that $\hat{\theta}-\theta_0=\mathrm{o}_{P}(1/\sqrt{n})$, which is not true for the least squares estimator under the local alternatives. Thus we give the corrected limiting result of $T^{Zh}_n$ and generalize it to obtain the results of $TM_{n}$ in the proof of Theorem~\ref{TM_1n}.
\end{remark}

\section{Particularity of checking ODE models}

If we reject $H_0$ in the above test, we may further want to identify the component(s) that is (are) wrongly modelled. In this situation, the hypotheses are as follows, for any $k$ with $1\le k\le p$,
\begin{eqnarray*}	
	&H_{0k}:& X^{\prime}_k(t) =\frac{d X_{k}(t)}{d t} = g_k(t) = f_k(t, X(t ; \theta_0); \theta_0)  \in \mathcal{F}_k,\\
	&H_{1k}:& X^{\prime}_k(t) = g_k(t) \notin \mathcal{F}_k,
\end{eqnarray*}
where $\theta_0$ is an unknown parameter vector.

However, as we briefly commented in the above section, although the trajectory matching-based test can detect the alternatives distinct form the whole system, it is quite incompetent to do this work. Actually, since the parametric ODE model structure is on $X^{\prime}(\cdot)$ instead of $X(\cdot)$, there are some extra difficulties for checking ODE models. In this section we will discuss three aspects of the particularity of checking ODE models and consider how to deal with them. We will use the idea proposed to construct two available tests for checking ODE component functions in the next two sections.

\begin{remark}
	Sometimes we may just have a model for some certain components instead of a model for all components. We wish to check whether the component(s) is (are) rightly modelled. This is another important case that we need to consider $H_{0k}$. In this case, since the model is not complete, we can not solve the ODEs and construct $TM_n$. However, the idea and tests proposed below are still available.
\end{remark}

\subsection{Mixed components}

To realize why $TM_{n}$ cannot identify the wrongly modelled component(s), think about the following toy example. Suppose that the ODE system to be tested are:
\begin{equation}	
X^{\prime}(t) = \left[ \begin{array}{c}{\frac{d X_{1}}{d t}} \\ {\frac{d X_{2}}{d t}} \end{array}\right]=\left[ \begin{array}{c}{f_{1}} \\ {f_{2}}\end{array}\right]=\left[ \begin{array}{c}{X_{1}} \\ {X_{1} + X_{2}}\end{array}\right]. \\
\end{equation}
Yet the true ODE system is:
\begin{equation}	
X^{\prime}(t) = \left[ \begin{array}{c}{\frac{d X_{1}}{d t}} \\ {\frac{d X_{2}}{d t}} \end{array}\right]=\left[ \begin{array}{c}{2 X_{1}} \\ {X_{1} + X_{2}}\end{array}\right] \\
\end{equation}
Here the first component is wrongly modelled. Recall the function $F(\cdot, \cdot)$ is a function of all components of $X(\cdot)$. If we use the trajectory matching-based test to check the second component, as it also involves the first component, the decision on whether rejecting or not for the null hypothesis will then make little sense. This shows the significant difference of this testing problem from the case with classical multi-response regression models.

To construct an available test, we should decouple the relationship among different components. Inspired by two-step collocation methods, we can do this by applying nonparametric techniques. Specifically, the nonparametric estimator $\hat{X}(t)$ is used to replace $X(t)$ in the parametric ODE model. Since the nonparametric estimator is model free, it always captures the true shape of the corresponding components and  decouple the relationship among different components. Then we can build the test grounded on the parametric model to be tested with the plug-in nonparametric estimator.

\subsection{Two types of local alternatives}

Researchers usually consider local misspecifications which convergence to the null model at a rate $\delta_{n}$ in the classical model checking tests. However, the things are more complicated for the ODE models. On the one hand, we can add the local disturbance to the trajectory of ODEs, as we set $H_{1n}^F$ in the last section. This setting is similar to the traditional regression model checking. On the other hand, since the ODEs model gradients instead of primitive functions, we may ponder the case that local misspecifications are added to the derivative functions. Thus, unlike the local alternatives $H_{1n}^F$ of (\ref{local})  about the original function $X(\cdot)$, we also consider the following sequence of local alternatives about the derivative $X'(\cdot)$:
\begin{equation}\label{local0}	
H_{1n}^f:X'(t) = f\left(t, X(t);\theta_{0}\right)+\delta_{n} l\left(t\right).
\end{equation}
To the best of our knowledge, such alternatives are never considered before in the literature. The following theorem states the relationship between $H^f_{1n}$ and $H^F_{1n}$.
\begin{theo}
	Given sets A-C of assumptions in Supplement A, then  under $H^F_{1n}$, the derivative has the form
	\begin{equation}
	X'(t) = f\left(t, X(t);\theta_{0}\right)+\delta_{n} v_1\left(t\right) + o(\delta_{n}) v_2(t).
	\end{equation}
	Under $H^f_{1n}$, the original function can be expressed as
	\begin{equation}
	X(t) = F\left(t, \theta_{0}\right)+\delta_{n}v_3\left(t\right) + o(\delta_{n}) v_4(t).
	\end{equation}\label{two_local}
\end{theo}
This phenomenon shows that adding the $\delta_n$ rate disturbance to the trajectory of ODEs is equivalent to adding a no slower than $\delta_{n}$ rate disturbance directly to ODEs and vice versa. Since the higher order little terms will vanish with a faster speed. They usually do not influence the asymptotic property of the test when we consider the local alternatives. Thus $H^F_{1n}$ is equivalent to $H^f_{1n}$ in this sense.

In the remaining part, we consider the two corresponding sequences of local alternatives in (\ref{local}) and (\ref{local0}) for any component function:
\begin{equation}\label{local1}
H_{1kn}^F: X_k(t) = F_k\left(t, \theta_{0}\right)+\delta_{n} L_k\left(t\right)
\end{equation}
with the counterpart function $l_k(t)\equiv v_{1k}(t)$, and
\begin{equation}\label{local2}	
H_{1kn}^f: X_k'(t) = f_k\left(t, X(t);\theta_{0}\right)+\delta_{n} l_k\left(t\right)
\end{equation}
with the counterpart function $L_k(t) \equiv v_{1k}(t)$.

\subsection{Mixed parameters}

Besides the phenomenon of mixed components mentioned above, the phenomenon of mixed parameters may also invalidate the trajectory matching-based test. This is because if different components share some same parameters, the wrongly modelled component(s) will make the estimators deviate from the true values when all  data is used. The estimators are not consistent and  thus the test relying on these estimators is ineffective. This problem can be solved if we use two-step collocation methods  to estimate the parameters when  the component to be tested is considered. Specifically, we estimate the parameters $\theta$ as follows:
\begin{equation}
\hat{\theta}_{TS}=\underset{\theta}{\operatorname{argmin}} \sum_{j=1}^{m}\left[\hat{X}_k^{\prime}\left(t_{j}^{*}\right)-f_k\left(t,\hat{X}\left(t_{j}^{*}\right) ; \theta\right)\right]^{2} \omega_k\left(t_{j}^{*}\right)
\end{equation}
with $\omega_k\left(t\right)$ being a selected weight function and $t_{j}^{*}$ being the selected time grid whose number $m$ can be larger than $n$. $\hat{X}(t)$ is the local linear estimator for $X(t)$ and $\hat{X}^{\prime}(t)$ is the local quadratic estimator for $X'(t)$ in the vector version, whose $k$-th components $\hat{X}_k(t)$ and $\hat{X}_k^{\prime}(t)$ are the corresponding local polynomial estimators for $X_k(t)$ and $X'_k(t)$. Define $h_e$ as the bandwidth.

Preparing for constructing tests,  we give the asymptotic properties of two-step collocation estimator under different hypotheses in the following theorem. Assume there exists a unique minimizer $\theta^*_{TS}$, such that
$$\theta^*_{TS}= \arg \min_{\theta} E_{p^*}\left\{ \left[X_k^{\prime}(t)-f_k\left(t, X(t), \theta\right)\right]^{2}w_k(t)\right\}$$
and denote $\Lambda\left(t\right)=\hat{X}(t)-X(t)$, $\Delta(t)=\hat{X}^{\prime}(t)-X^{\prime}(t)$. We have the following theorem.
\begin{theo}
	Given sets A and C of assumptions in Supplement A,  $\ln n/(n h_e^3)=o(1)$, the two-step collocation estimator $\hat{\theta}_{TS}$ is a consistent estimator of $\theta^*_{TS}$. Further,
	
	1. Under the null hypothesis, we have $\theta^*_{TS}=\theta_0$ and
	\begin{equation}
	\begin{aligned}
	\hat{\theta}_{TS}-\theta_{0}&=H_{\dot{f}}^{-1}\frac{1}{m} \sum_{j=1}^{m}\left[ \Delta_k(t_j^*) \omega_k(t_j^*)\frac{\partial f_k\left(t,X\left(t_{j}^*\right), \theta_0 \right)}{\partial \theta} \right.\\
	&\left.- \omega_k(t_j^*)\frac{\partial f_k\left(t,X\left(t_{j}^*\right), \theta_0 \right)}{\partial \theta} \frac{\partial f_k\left(t,X\left(t_{j}^*\right), \theta_0 \right)}{\partial X^\mathrm{T}} \Lambda\left(t_{j}^{*}\right)\right]+\mathrm{o}_{P}(n^{-1/2})\\
	\end{aligned}
	\end{equation}
	which is a term of order $\mathrm{o}_{P}(n^{-1/2})$. Here
	\begin{equation}
	H_{\dot{f}}=E_{p^*} \left[\omega_k(t)\frac{\partial f_k\left(t,X\left(t\right), \theta_0 \right)}{\partial \theta}\frac{\partial f_k\left(t,X\left(t\right), \theta_0 \right)}{\partial \theta^{\mathrm{T}}}\right].
	\end{equation}
	
	2. Under the global alternative hypothesis $H_1$, we have $\theta^*_{TS}=\theta_1$ and
	\begin{equation}
	\begin{aligned}
	&\sqrt{n}(\hat{\theta}_{TS}-\theta_{1})\\
	=&G^{-1}\frac{\sqrt{n}}{m} \sum_{j=1}^{m}\left[ \Delta_k(t_j^*) \omega_k(t_j^*)\frac{\partial f_k\left(t,X\left(t_{j}^*\right), \theta_1 \right)}{\partial \theta}\right.\\
	&\left.- \omega_k(t_j^*)\frac{\partial f_k\left(t,X\left(t_{j}^*\right), \theta_1 \right)}{\partial \theta} \frac{\partial f_k\left(t,X\left(t_{j}^*\right), \theta_1 \right)}{\partial X^\mathrm{T}}
	\Lambda\left(t_{j}^{*}\right)\right]+\mathrm{o}_{P}(1)\\
	\end{aligned}
	\end{equation}
	where
	\begin{equation}
	\begin{aligned}
	G=&E_{p^*} \left[\omega_k(t)\frac{\partial f_k\left(t,X\left(t\right), \theta_1 \right)}{\partial \theta}\frac{\partial f_k\left(t,X\left(t\right), \theta_1 \right)}{\partial \theta^{\mathrm{T}}}\right]\\
	&-E_{p^*} \left\{ \left[g_k(t)-f_k\left(t,X\left(t\right), \theta_1 \right)\right]\omega_k(t) \frac{\partial^{2} f_k\left(t,X\left(t\right), \theta_1 \right)}{\partial \theta \partial \theta^{\mathrm{T}}} \right\}.
	\end{aligned}
	\end{equation}
	
	3. Under the local alternative hypothesis $H_{1kn}^F$ in (\ref{local1}) or $H_{1kn}^f$ in (\ref{local2}) with $\delta_n \rightarrow 0$, we have $\theta^*_{TS}=\theta_0$ and
	\begin{equation}
	\begin{aligned}
	\sqrt{n}(\hat{\theta}_{TS}-\theta_{0})=&H_{\dot{f}}^{-1}\frac{\sqrt{n}}{m} \sum_{j=1}^{m}\left[ \Delta_k(t_j^*) \omega_k(t_j^*)\frac{\partial f_k\left(t,X\left(t_{j}^*\right), \theta_0 \right)}{\partial \theta}\right.\\
	&\left.- \omega_k(t_j^*)\frac{\partial f_k\left(t,X\left(t_{j}^*\right), \theta_0 \right)}{\partial \theta} \frac{\partial f_k\left(t,X\left(t_{j}^*\right), \theta_0 \right)}{\partial X^\mathrm{T}}
	\Lambda\left(t_{j}^{*}\right)\right]\\
	&+\sqrt{n} \delta_n H_{\dot{f}}^{-1} E_{p^*} \left[ l(t_j^*) \omega(t_j^*)\frac{\partial     f_k\left(t,X\left(t_{j}^*\right), \theta_0 \right)}{\partial \theta} \right]
	+\mathrm{o}_{P}(1).\\
	\end{aligned}
	\end{equation}\label{TS}	
	
\end{theo}

\section{Integral matching-based test}

\subsection{The hypotheses and test statistic}

Similar to integral matching, we use the indefinite integral instead of trajectory of ODEs to construct pseudo-residuals:
\begin{equation}	
\hat{e}_{ik} = Y_{ik}-X_k(t_0)-\int_{t_0}^{t_i} f_k\left(t,\hat{X}\left(t\right) ; \hat{\theta}\right)\mathrm{d}t. \\
\end{equation}
Here we use the local linear estimator $\hat{X}(t)$ with the bandwidth $h_0$ and the two-step collocation estimator $\hat{\theta}$.  Since $\hat{F}_k(t_i;\hat{\theta})=X_k(t_0)+\int_{t_0}^{t_i} f_k\left(t,\hat{X}\left(t\right); \hat{\theta}\right)\mathrm{d}t$ is expected to converge to $F_k(t_i,\theta_{0})$, $\hat{e}_{ik}$ can be used as a surrogate to replace $e_{ik}=Y_{ik}-F_k(t_i;\hat{\theta})$ in the trajectory matching-based test. Consequently we obtain an integral matching-based test and write it as $IM_{n}(k)$ in short.

However, to simplify notation  without confusion, we simply write $IM_{n}(k)$ as $IM_{n}$ and other statistics as ones without the subscript $k$ to indicate the corresponding component unless we need to stress its role in analysis. Define
\begin{equation}
\begin{aligned}
IM_{n}=&\sqrt{\frac{n-1}{n}} \frac{nh^{1/2}V_n^{\hat{F}}}{\sqrt{\widehat{\Sigma}^{\hat{F}}}}\\
=&\frac{\sum_{i=1}^{n} \sum_{j=1 \atop j \neq i}^{n} K\left(\frac{t_{i}-t_{j}}{h}\right) \hat{e}_{ik} \hat{e}_{jk}}{\left\{\sum_{i=1}^{n} \sum_{j=1 \atop j \neq i}^{n} 2 K^{2}\left(\frac{t_{i}-t_{j}}{h}\right) \hat{e}_{ik}^2\hat{e}_{jk}^2\right\}^{1 / 2}}
\end{aligned}	
\end{equation}
where
\begin{equation}	
V_{n}^{\hat{F}} = \frac{1}{n(n-1)} \sum_{i=1}^{n} \sum_{j=1 \atop j \neq i}^{n} \frac{1}{h} K\left(\frac{t_{i}-t_{j}}{h}\right) \hat{e}_{ik} \hat{e}_{jk},
\end{equation}
\begin{equation}
\widehat{\Sigma}^{\hat{F}}=\frac{2}{n(n-1)} \sum_{i=1}^{n} \sum_{j=1 \atop j \neq i}^{n} \frac{1}{h} K^{2}\left(\frac{t_{i}-t_{j}}{h}\right)\hat{e}_{ik}^2\hat{e}_{jk}^2.
\end{equation}

The form of this integral matching-based test is analogous to the trajectory matching-based test except that $e_i$ is substituted by $\hat{e}_i$. As we mentioned in the last section, this replacement is critical since $\hat{X}(t)$ always captures the true form $X(t)$ which eliminates the influence of latent wrong modelled components.



\subsection{Asymptotic properties}

We hereafter have to deal with a high-order U-statistic. Thus we first give a lemma to establish the asymptotic equivalence of $U_n$ and $\hat{U}_n$ and present the limiting distribution of a non-degenerate U-statistic of order $m^*$ with a kernel varying with $n$.
\begin{lemma}
	Suppose $U_{n}$ is an U-statistic with the kernel $h_{n}\left(z_{1},\cdots,z_{m^*}\right)$ of order $m$. If $E\left[\left\|h_{n}\left(z_{1},\cdots,z_{m^*}\right)\right\|^{2}\right]=o(n)$, then
	\begin{equation}
	\sqrt{n}\left(U_{n}-\hat{U}_{n}\right)=o_{p}(1)
	\end{equation}
	where $\hat{U}_{n}=E\left[h_{n}\left(z_{1},\cdots,z_{m^*}\right)\right]+\frac{m^*}{n} \sum_{i=1}^{n}\{ E\left[h_{n}\left(z_{1},\cdots,z_{m^*}\right) | z_{i}\right]-E\left[h_{n}\left(z_{1},\cdots,z_{m^*}\right)\right] \}$ is the projection of $U_n$.\label{LI.1}
\end{lemma}

By denoting $\hat{e}_i=\hat{e}_i+\varepsilon_{i}-\varepsilon_{i}$, $V_{n}^{\hat{F}}$ can be decomposed as an U-statistics plus remaining  terms. Applying Lemma~\ref{LI.1}, the remaining terms can be showed as $\mathrm{o}_{P}(n^{-1}h^{-1/2})$ and thus the replacement of $\hat{e}_i$ to $e_i$ have no influence to the limiting null distribution of $IM_{n}$ under certain regularity conditions. Define $a_n(h)=h^2+n^{-1/2}h^{-1/2}\log{n}^{-1/2}$ which is the uniform convergence rate of local linear estimator (\cite{Hansen.2008}). We state this property in the following lemma.
\begin{lemma}
	Given sets A-C of assumptions in Supplement A, if $h \rightarrow 0$, $n h \rightarrow \infty$, $ h_0=o(n^{-1/4}h^{-1/4})$, $n^{-1/2}h^{-1/2}=o(h_0)$ and $a_n^2(h_0)=o(n^{-1}h^{-1/2})$, then
	\begin{equation}
	nh^{1/2}V_{n}^{\hat{F}}=nh^{1/2}V_{1n}+o_{p}(1),
	\end{equation}
	where
	\begin{equation}	
	V_{1n}= \frac{1}{n(n-1)} \sum_{i=1}^{n} \sum_{j=1 \atop j \neq i}^{n} \frac{1}{h} K\left(\frac{t_{i}-t_{j}}{h}\right) e_{ik} e_{jk}.
	\end{equation}\label{LI.2}
\end{lemma}

Having Lemma~\ref{LI.2}, we can easily derive the asymptotic properties of $IM_{n}$ in the following theorems.
\begin{theo}
	Given sets A-C of assumptions in Supplement A, if $h \rightarrow 0$, $n h \rightarrow \infty$, $ h_0=o(n^{-1/4}h^{-1/4})$, $n^{-1/2}h^{-1/2}=o(h_0)$ and $a_n^2(h_0)=o(n^{-1}h^{-1/2})$, then under the null hypothesis,
	\begin{equation}
	IM_{n} \stackrel{\mathbf{d}}{\longrightarrow} N(0,1).
	\end{equation}\label{IM_0}
\end{theo}
\begin{theo}
	Given sets A-C of assumptions in Supplement A, if $h \rightarrow 0$, $n h \rightarrow \infty$, $ h_0=o(n^{-1/4}h^{-1/4})$, $n^{-1/2}h^{-1/2}=o(h_0)$ and $a_n^2(h_0)=o(n^{-1}h^{-1/2})$, then under $H_{1k}$,
	\begin{equation}
	IM_{n}/(nh^{1/2}) \stackrel{\mathbf{P}}{\longrightarrow} \frac{\mathrm{E} \left\{\left[X_k\left(t\right)-F_k\left(t, \theta_{1}\right)\right]^{2} p\left(t_{i}\right) \right\}}{\left\{2 \int K^{2}(u) \mathrm{d} u \int\left\{\sigma_k^{2}(t)+\left[X_k(t)-F_k\left(t, \theta_{1}\right)\right]^{2}\right\}^{2} p^{2}(t) \mathrm{d} t\right\}^{1 / 2}}.
	\end{equation}\label{IM_1}
\end{theo}

Theorem~\ref{IM_1} shows that this test is consistent and sensitive to the global alternative in the sense that it can diverge to infinity at the rate of order $nh^{1/2}$. Recall that under the global alternative, the trajectory-matching-based test can diverges to infinity at the of order $n^2h$, and thus seems more powerful than  the integral matching-based test developed here. But note that the trajectory-matching-based test is a quadratic form and thus, its critical value is also larger than that for the integral matching-based test. Thus this is not comparable.

The following theorem states the asymptotic property of $IM_{n}$ under $H_{1kn}$.
\begin{theo}
	Given sets A-C of assumptions in Supplement A, if $h \rightarrow 0$, $n h \rightarrow \infty$, $ h_0=o(n^{-1/4}h^{-1/4})$, $n^{-1/2}h^{-1/2}=o(h_0)$ and $a_n^2(h_0)=o(n^{-1}h^{-1/2})$, we have the following asymptotic property of $IM_{n}$ under $H_{1kn}^F$ or $H_{1kn}^f$.
	
	With $n^{1 / 2} h^{1 / 4}\delta_{n} \rightarrow \infty$,
	\begin{equation}
	IM_{n}/(nh^{1/2}\delta_{n}^{2}) \stackrel{\mathbf{P}}{\longrightarrow} \mu_I/ \sigma_k
	\end{equation}
	where
	\begin{equation}
	\begin{aligned}
	\mu_I=&\mathrm{E}\left\{\left\{ L_k\left(t\right)-\frac{\partial F_k\left(t, \theta_{0}\right)}{\partial \theta^{T}} H^{-1}_{\dot{f}}\mathrm{E}_g\left[l_k\left(t\right)\omega_k(t)\frac{\partial f_k\left(t, X(t); \theta_{0}\right)}{\partial \theta}\right]  \right\}^2 p(t) \right\}.
	\end{aligned}
	\end{equation}
	
	Particularly, when $\delta_{n}=n^{-1 / 2} h^{-1 / 4}$,
	\begin{equation}
	IM_{n} \stackrel{\mathbf{d}}{\longrightarrow} N(\mu_I/ \sigma_k,1).
	\end{equation}\label{IM_1n}
	
\end{theo}

This result again shows the similar sensitivity to the local alternatives as  classical local smoothing tests do for regressions.

\section{Gradient matching-based test}

\subsection{The hypotheses and test statistic}

Reminded by gradient matching, we can also check the component(s) in ODE models by directly using the gradient of ODEs instead of the trajectory. We then define a gradient matching-based test ($GM_{n}$). To be more specific, we first use Nadaraya-Watson kernel estimation to estimate $X(t)$ and $X^{\prime}(t)$ as
\begin{eqnarray}
\hat{X}(t)&=& \hat{h}(t) / \hat{p}\left(t\right),\nonumber\\
\hat{X}^{\prime}(t)&=& \frac{\hat{h}^{\prime}(t) \hat{p}(t)-\hat{h}(t) \hat{p}^{\prime}(t)}{\hat{p}^2\left(t\right)},
\end{eqnarray}
where
\begin{eqnarray}
\hat{h}(t)&=&\frac{1}{n} \sum_{i=1}^{n} \frac{1}{h}K\left(\frac{t-t_{i}}{h} \right)Y_i \quad
\hat{h}^{\prime}(t)=\frac{1}{n} \sum_{i=1}^{n} \frac{1}{h^2}K^{\prime}\left(\frac{t-t_{i}}{h} \right)Y_i \nonumber\\
\hat{p}(t)&=&\frac{1}{n} \sum_{i=1}^{n} \frac{1}{h}K\left(\frac{t-t_{i}}{h}\right) \quad
\hat{p}^{\prime}(t)=\frac{1}{n} \sum_{i=1}^{n} \frac{1}{h^2}K^{\prime}\left(\frac{t-t_{i}}{h}\right).
\end{eqnarray}

Similar to Section~4, we again simplify notation  without confusion, we simply write statistics as ones without the subscript $k$ to indicate the corresponding component unless we need to stress its role in analysis.	

Under the null hypothesis, $e_f(t)\equiv X^{\prime}_k(t) -f_k(t, X(t ; \theta_0); \theta_0) =0$ while it is not the case under the alternatives. Thus, if we replace $X^{\prime}_k(t)$ by $\hat{X}^{\prime}_k(t)$ and  $f_k(t, X(t ; \theta_0); \theta_0)$ by $f_k(t, \hat{X}(t); \hat{\theta})$, the preudo-residual $\hat{e}_f(t)=\hat{X}^{\prime}_k(t)-f_k(t, \hat{X}(t); \hat{\theta})$ is expected to converge to zero in probability. Therefore, $E\left[\hat{e}_f^{2}(t_i) \hat{p}^4\left(t_{i}\right)\right]$ is expected to converge to zero under the null hypothesis while to a positive constant under the alternative hypothesis, where $\hat{p}^4(t_i)$ is used to eliminate the denominator in the nonparametric estimation. Then we construct a test statistic:	
\begin{equation}
\begin{aligned}
V_n^f=& \frac{1}{nh^2}\sum^{n}_{d=1}\left[\hat{X}_k^{\prime}(t_d)-f_k(t, \hat{X}(t_d); \hat{\theta}) \right]^2\hat{p}^4\left(t_{d}\right) \\
=& \frac{1}{nh^2}\sum^{n}_{d=1}\left[\hat{h}_k^{\prime}(t_d) \hat{p}(t_d)-\hat{h}_k(t_d) \hat{p}^{\prime}(t_d)-\hat{p}^2\left(t_{d}\right)f_k(t, \hat{X}(t_d); \hat{\theta})\right]^2 \\
=&\frac{1}{nh^2}\sum^{n}_{d=1} \left\{ \frac{1}{(n-1)^2}\sum^{n}_{i=1} \sum^{n}_{j=1} \left[\frac{1}{h^3} K^{\prime}\left(\frac{t_{d}-t_{i}}{h}\right)K\left(\frac{t_{d}-t_{j}}{h}\right) (Y_{ik}-Y_{jk}) \right.\right.\\
& \left. \left.- \frac{1}{h^2}K\left(\frac{t_{d}-t_{i}}{h}\right)K\left(\frac{t_{d}-t_{j}}{h}\right)f_k(t, \hat{X}(t_d); \hat{\theta}) \right] \right\}^2.
\end{aligned}
\end{equation}
Note that we add an extra multiplier $1/h^2$ compared with the previous tests.  This is because the preudo-residual $\hat{e}_f(t)=\hat{X}_k^{\prime}(t)-f_k(t, \hat{X}(t); \hat{\theta})$ is expected to converge to zero rather than a zero mean random variable in probability under the null hypothesis. This is a significant difference from the previous ones. This is required when we need a non-degenerate limit. Here, the test statistic $V_n^f$ can also be asymptotically a V-statistic that is further asymptotically equivalent to the corresponding U-statistic according to  Theorem 1 of
\cite{MartinsFilho.2006}. Therefore its asymptotic properties can be derived by the U-statistics theory. However, as the  nonparametric density estimation is biased, which causes a non-negligible  bias term of this test statistic even under the null hypothesis. We now use a bias correction. Check the mean of this test statistic:
\begin{equation}
\begin{aligned}
V&=E(V^f_n)\\
=&\frac{1}{h^2}\int \left[ f_k(t, X(t); \theta^*)-X_k^{\prime}(t) \right]^2p(t)^5\mathrm{d}t
\\&+h^2 \left[\int \frac{u^3}{6}K^{\prime}(u)\mathrm{d}u\right]^2\int X_k^{(3)}(t)^2p(t)^5\mathrm{d}t+o(h^2).
\end{aligned}
\end{equation}
It is easy to see the term $S := 1/h^2\int \left[ f_k(t, X(t); \theta_0)-X_k^{\prime}(t) \right]^2p(t)^5\mathrm{d}t $ is zero under the null hypothesis while nonzero under the alternative hypothesis. We will use this piece of information to construct a new test. We randomly partition the original sample into $2$ subsamples. Using these two subsamples, we construct two test statistics $V^{f}_{\tilde{n}1}$ and $V^{f}_{n-\tilde{n}2}$, where $\tilde{n}=\lfloor n/2 \rfloor$. As $n-2\tilde n\le 1$, the asymptotic properties of $V^{f}_{(n-\tilde{n})2}$ should be the same as those of $V^{f}_{\tilde{n}2}$. Thus, we assume that, without loss of generality, $n=2\tilde{n}$ is even. The difference $V^{f}_{\tilde{n}1}- V^{f}_{\tilde{n}2}$ is a new statistic.   However, it is clearly not a powerful test as even under the alternatives, its mean is also zero. Therefore, we use $V^{f}_{\tilde{n}1}- V^{f}_{\tilde{n}2}+cS$ instead. This quantity fully uses the information provided above: $S=0$ under the null and $>0$ under the alternatives.  The estimator of $S$ is obtained as  $\hat{S}=1/h^2\int  \left[f_k(t,\hat{X}(t);\hat{\theta})-\hat{X}_k^{\prime}(t)\right]^2\hat{p}(t)^5\mathrm{d}t $. The local  linear smoother and local quadratic  smoother are used to obtain $\hat{X}(t)$ and $\hat{X}^{\prime}(t)$ respectively with the corresponding bandwidths $h_0$ and $h_1$ while the kernel density estimation is used to get $\hat{p}(t)$ with the bandwidth $h_0$. The parameters $\theta$ is again estimated by using the two-step collocation method. But, in our theoretical development, we need to assume the boundedness of the support for the density function, we then can  use $\hat{S}^{\prime}=1/h^2\int  \left[f_k(t,\hat{X}(t);\hat{\theta})-\hat{X}_k^{\prime}(t)\right]^2\mathrm{d}t$ without the estimator $\hat p (\cdot)$ of $p(\cdot)$ in the numerical studies. It is an estimator of  ${S}^{\prime}=1/h^2\int  \left[f_k(t,{X}(t);{\theta})-{X}_k^{\prime}(t)\right]^2\mathrm{d}t$ that is equal to zero under the null as well. Without notational confusion, we still write $S'$ as $S$ throughout the paper.
Note that we use two different bandwidths for $X$ and $X'$ to  ensure suitable rates of convergence.

By correcting the bias and dividing by the estimator of its variance, we can modify $V_n^f$ to construct the final test statistic $GM_{n}$ as
\begin{equation}
\begin{aligned}
GM_{n}=\frac{\sqrt{\tilde{n}}(V^{f}_{\tilde{n}1}-V^{f}_{\tilde{n}2}+c\hat{S})}{\sqrt{2\widehat{\Sigma}^f}},
\end{aligned}
\end{equation}
where
\begin{eqnarray}
\widehat{\Sigma}^f &=&  \frac{1}{n-1}\sum_{k}^{n}[\hat{w}_n(z_s)-\frac{1}{n}\sum_{i=1}^{n}\hat{w}_n(z_i)]^2,\nonumber\\
\hat{w}_n(z_s) &=& \frac{1}{\lfloor \frac{n-1}{4} \rfloor} \sum_{i=1}^{\lfloor \frac{n-1}{4} \rfloor} W_n(z_{1i},z_{2i},z_{3i},z_{4i},z_{s}),\nonumber\\
W_n(z_a,z_b,z_c,z_d,z_s)&=&\frac{1}{5!}\sum_{P}W_n^{\prime}(z_a,z_b,z_c,z_d,z_s)\nonumber\\
\end{eqnarray}
and
\begin{equation}
\begin{aligned}
&W_n^{\prime}(z_a,z_b,z_c,z_d,z_s)\\
=&\frac{1}{h^2} K\left(\frac{t_{s}-t_{a}}{h}\right)K\left(\frac{t_{s}-t_{b}}{h}\right)\\
&\times \left[\frac{1}{h^3} K^{\prime}\left(\frac{t_{s}-t_{c}}{h}\right) (Y_{ck}-Y_{ak}) - \frac{1}{h^2}K\left(\frac{t_{s}-t_{c}}{h}\right)f_k(t_s, \hat{X}(t_s); \hat{\theta})\right]\\
&\times \left[\frac{1}{h^3} K^{\prime}\left(\frac{t_{s}-t_{d}}{h}\right) (Y_{dk}-Y_{bk})  - \frac{1}{h^2}K\left(\frac{t_{s}-t_{d}}{h}\right)f_k(t_s, \hat{X}(t_s); \hat{\theta})\right].\\
\end{aligned}
\end{equation}

%

\subsection{Asymptotic properties}

Unlike the previous tests, we here have to deal with an U-statistic of order $5$. 
Let $b_n(h)=h^2+n^{-1/2}h^{-3/2}\log{n}$ which is the uniform convergence rate of $X^{\prime}(t)$ (\cite{Liang.2008}). Applying Lemma~\ref{LI.1}, we can provide the asymptotic properties of $V_n^f$ under the null hypothesis in the following lemma.
\begin{lemma}
	Given sets A and C of assumptions in Supplement A, if $h^{-12}=o(n)$, $a_n^2(h_0)h^{-2}=o(n^{-1/2})$ and $b_n^2(h_1)h^{-2}=o(n^{-1/2})$, then under the null hypothesis,
	\begin{equation}
	\sqrt{n}\left(V_n^f-V\right) \stackrel{\mathbf{d}}{\longrightarrow} N(0,\Sigma^f),
	\end{equation}
	where $\Sigma^f=\frac{1}{9}(\int u^3K^{\prime}(u)\mathrm{d}u)^2\int (X^{(4)}(t_k))^2 \sigma^2(t_k) p^8(t_k)\mathrm{d}t_k$.\label{LG.1}
\end{lemma}

The following two lemmas give the asymptotic properties of the estimators of $S$ and $\Sigma^f$ under $H_{0k}$.
\begin{lemma}
	Given sets A and C of assumptions in Supplement A, then under the null hypothesis, if $h^{-12}=o(n)$, $a_n^2(h_0)h^{-2}=o(n^{-1/2})$ and $b_n^2(h_1)h^{-2}=o(n^{-1/2})$,
	\begin{eqnarray*}
		\sqrt{n}\hat{S}  \stackrel{\mathbf{P}}{\longrightarrow} 0. 
	\end{eqnarray*}\label{LG.2}
\end{lemma}
\begin{lemma}
	Given sets A and C of assumptions in Supplement A, then under the null hypothesis, if $h^{-12}=o(n)$, $a_n^2(h_0)h^{-2}=o(n^{-1/2})$ and $b_n^2(h_1)h^{-2}=o(n^{-1/2})$,
	\begin{equation}
	\widehat{\Sigma}^f \stackrel{\mathbf{P}}{\longrightarrow} \Sigma^f.
	\end{equation}\label{LG.3}
\end{lemma}

Now we state the asymptotic property of $GM_{n}$ under the null hypothesis.
\begin{theo}
	Given sets A and C of assumptions in Supplement A, then under the null hypothesis, if $h^{-12}=o(n)$, $a_n^2(h_0)h^{-2}=o(n^{-1/2})$ and $b_n^2(h_1)h^{-2}=o(n^{-1/2})$, recalling that $\tilde{n}=[n/2]$,
	\begin{equation}
	GM_{n} \stackrel{\mathbf{d}}{\longrightarrow} N(0,1).
	\end{equation}\label{GM_0}
\end{theo}

Next, we present the asymptotic distribution of $GM_{n}$ under the global alternative hypothesis.
\begin{theo}
	Given sets A and C of assumptions in Supplement A, then under the global alternative hypothesis, if $h^{-12}=o(n)$, $a_n^2(h_0)h^{-2}=o(n^{-1/2})$ and $b_n^2(h_1)h^{-2}=o(n^{-1/2})$, recalling that $\tilde{n}=[n/2]$,
	\begin{equation}
	GM_{n}/\sqrt{\tilde{n}} \stackrel{\mathbf{P}}{\longrightarrow} \frac{c\int \left[f_k(t, X(t); \theta_1)-X_k^{\prime}(t)\right]^2p(t)^5\mathrm{d}t}{\sqrt{2\Sigma^{f\prime}}} >0
	\end{equation}
	where
	\begin{equation}
	\begin{aligned}
	\Sigma^{f\prime}=&\int\left\{ 25\left[f_k(t, X(t); \theta_1)-X_k^{\prime}(t)\right]^4p^8(t) \right.\\
	&\left. +4\left[f_k^{\prime}(t, X(t); \theta_1)-X_k^{(2)}(t)\right]^2 \sigma_k^2(t) p^8(t)\right\}\mathrm{d}t\\
	&-25\left\{\int \left[f_k(t, X(t); \theta_1)-X_k^{\prime}(t)\right]^2p^4(t)\mathrm{d}t \right\}^2.
	\end{aligned}
	\end{equation}\label{GM_1}
\end{theo}

Theorem~\ref{GM_1} shows that the test is consistent and diverges to infinity at the rate of $\sqrt{n}$. We will make a comparison of their performance in the numerical studies. 

The following theorem states the asymptotic power of $GM_{n}$ under $H_{1kn}^F$ and $H_{1kn}^f$.
\begin{theo}
	Given sets A and C of assumptions in Supplement A, if $h^{-12}=o(n)$, $a_n^2(h_0)h^{-2}=o(n^{-1/2})$ and $b_n^2(h_1)h^{-2}=o(n^{-1/2})$, recalling that $\tilde{n}=[n/2]$, then under $H_{1kn}^F$ or $H_{1kn}^f$, with $\tilde{n}^{1/4}h^{-1}\delta_n \rightarrow \infty$ and $\delta_n h^{-1}=o(1)$,
	\begin{equation}
	GM_{n}/(\tilde{n}^{1/2}h^{-2}\delta_n^2) \stackrel{\mathbf{P}}{\longrightarrow} c\mu_4/ \sqrt{2\Sigma^f}
	\end{equation}
	where
	\begin{equation}
	\begin{aligned}
	\mu_{4}=&\left\{H_{\dot{f}}^{-1} E_{p^*}\left[l_k\left(t\right) \omega_k\left(t\right) \frac{\partial f_{k}\left(t,X\left(t\right), \theta_{0}\right)}{\partial \theta}\right]\right\}^{\mathrm{T}}\\
	&\times \left[ \int \frac{\partial f_k(t, X(t) ; \theta_0)}{\partial \theta} \frac{\partial f_k(t, X(t) ; \theta_0)}{\partial \theta^{\mathrm{T}}} \mathrm{d} t\right]\\
	&\times \left\{H_{\dot{f}}^{-1} E_{p^*}\left[l_k\left(t\right) \omega_k\left(t\right) \frac{\partial f_{k}\left(t,X\left(t\right), \theta_{0}\right)}{\partial \theta}\right]\right\}+\int l_k^2(t) \mathrm{d} t\\
	&-2\left[\int l_k(t) \frac{\partial f_k(t, X(t) ; \theta_0)}{\partial \theta^{\mathrm{T}}} \mathrm{d} t \right] H_{\dot{f}}^{-1} E_{p^*}\left[ l_k\left(t\right) \omega_k\left(t\right) \frac{\partial f_{k}\left(t,X\left(t\right), \theta_{0}\right)}{\partial \theta}\right].
	\end{aligned}
	\end{equation}
	
	Particularly, if $\delta_{n}=\tilde{n}^{-1 / 4} h$,
	\begin{equation}
	GM_{n} \stackrel{\mathbf{d}}{\longrightarrow} N(c\mu_4/ \sqrt{2\Sigma^f}, 1).
	\end{equation}\label{GM_1n}
\end{theo}


\section{Numerical studies}

\subsection{Simulations}
We now conduct several simulations to evidence the performance of the proposed tests in finite sample scenarios. Three simulation studies are considered. In each study, we use $TM_{n}$ to check the whole ODE system while use $IM_{n}$ and $GM_{n}$ to check each component in ODE models. The number in subscript is used to denote which component the tests check. For example, $GM_{n1}$ is the gradient matching-based test for the first component in ODE models. Therefore, we give five tests in each study and examine their power and size. In Study 1, the null models are set to be the linear ODE system. Study 2 and Study 3 uses two nonlinear ODE models often used in neuroscience and ecology as the null models. 

Given the ODE models, we obtain the trajectory $X(t;\theta_{0})$ under the null by using the 4-stage Runge-Kutta algorithm. Then the observation values $Y(t)=X(t)+\epsilon(t)$ can be constructed. In the following studies, the observation time points $t_i$, $i=1,\ldots, n$ are independently generated from the uniform distribution $U(0,1)$. The error terms $\epsilon_{i}$, $i=1,\ldots, n$ are independent and identically distributed following the normal distribution $N(0,\sigma^2_{\epsilon} I_2)$. The initial values of ODE models are supposed to be given. In nonparametric estimation, we use Epanechnikov kernel $K(u)=\frac{3}{4}(1-u^2)$. The replication time is $1000$ for each simulation case. The significance level is $0.05$.

For constructing $TM_{n}$, we apply nonlinear least squares method to estimate $\theta$ which is implemented  using OPTI Toolbox (\cite{Currie.2012}). Then we obtain the trajectory $X(t;\hat{\theta})$ using the 4-stage Runge-Kutta algorithm. The bandwidth is chosen to be $h=0.05 \times n^{-2/5}$ by the rule of thumb.

For $IM_{n}$ and $GM_{n}$, the two-step collocation method is used to estimate $\theta$. In the estimation procedures, the local linear and quadratic smoother is applied to obtain $\hat{X}(t)$ and $\hat{X}^{\prime}(t)$ respectively. The time grid $t^*_j$ is chosen to be equidistantly distributed in the time interval. Following the advice in \cite{Fang.October2010}, we set $m=2\times \lfloor n^{4/3} \rfloor$ to improve the performance of the two-step collocation method. The weight function $\omega(t)$ is selected to be piecewise linear with the decreasing weights for points near boundary. To satisfy the condition of $h_e$, we select the bandwidth $h_e=\hat{h}_{opt} \times n^{-2/15} \times \ln^{1/2} n$, where $\hat{h}_{opt}$ is an estimator of the optimal bandwidth of kernel regression smoothing. We calculate this value from the R package `lokern' (\cite{MaechlerMartin.2010}).

Rather than estimate $\theta$, we also need the local linear estimator $X(t)$ and the local quadratic estimator $X'(t)$ to replace the true form of $X(t)$ and $X'(t)$ respectively in $IM_{n}$ and $GM_{n}$. We choose $h_0=\hat{h}_{opt}$ and $h_1=\hat{h}^f_{opt}$. Again, these optimal bandwidths are calculated from the R package `lokern' (\cite{MaechlerMartin.2010}).

The bandwidth for $IM_{n}$ is $h=0.025 \times n^{-3/5} \times \ln^{1/2} n$ while the bandwidth for $GM_{n}$ is $h=1 \times n^{-1/29}$. 
We have these choices of the rates because we need to meet the requirements for the consistency of the test statistics. As for the choice of the tuning parameter $c$ in $GM_n$, it has no significant influence for the asymptotic properties. However, it affects the power and size in finite sample performance. We found that a small positive $c$ is good for maintaining the significance level while a larger one is in favor of  power performance. Thus, a trade-off is in need. We recommend to use $1$ for linear ODE model setting, while use a more conservative value for the complex nonlinear ODE setting. Here we choose $c=1$ in the Study~1 while choose $c=0.2$ in the Study~2 and Study~3.

{ In particular, simulation results show that the empirical size of $IM_n$ tends to be very large in the complex ODE model settings. We have tried several sizes of sample, from $300$ to $10,000$ and found, from the unreported results, that with increasing the sample size, the empirical size can  gradually become smaller, though still large. This seems to show that the test can still be consistent, but the cumulative error by integration and involved nonparametric estimation very much affect the performance of the test $IM_n$.  This is because   the trajectory $X(t)$ of complex ODE system usually has a complex nonlinear function form and thus  the nonparametric estimation $\hat X(t)$ for all time points $t$ can have more serious estimation error, and the integral over  the surrogate $\hat{e}_{ik}$  of $e_{ik}$ in finite sample cases can  cause very large cumulative error of $IM_n$. Empirically, to  control the empirical size of $IM_n$, we make an adjusted version such that it can be applied at least for some simple models. First, 
	to alleviate the boundary effect of the estimation and the cumulative error by integration, we consider the following modification. First, we restrict the integral in the shorter interval $(0.1, 0.9)$ rather than the whole interval $(0,1)$ to avoid the boundary effect. Second, to reduce the error caused by the integration, we  split the  interval  into $n_l=8$ equidistant parts $\mathcal{T}_l=(l/10,(l+1)/10)$, $l=1,2,\cdots,8$ and define the corresponding residuals and test statistics as
	\begin{equation*}
	\begin{aligned}
	\hat{e}_{ik}^l=&\left[Y_{ik}-\hat{X}_k(\frac{l}{10})-\int_{min(\frac{l}{10},t_i)}^{t_i} f_k\left(t,\hat{X}\left(t\right) ; \hat{\theta}\right)\mathrm{d}t\right]I(t_i \in \mathcal{T}_l),\\
	IM_n^l=&\frac{\sum_{i=1}^{n} \sum_{j=1 \atop j \neq i}^{n} K\left(\frac{t_{i}-t_{j}}{h}\right) \hat{e}_{ik}^l \hat{e}_{jk}^l}{\left\{\sum_{i=1}^{n} \sum_{j=1 \atop j \neq i}^{n} 2 K^{2}\left(\frac{t_{i}-t_{j}}{h}\right) \hat{e}_{ik}^{l2}\hat{e}_{jk}^{l2}\right\}^{1 / 2}},
	\end{aligned}
	\end{equation*}
	where $I(\cdot)$ is the characteristic function. Then we define a test statistic as
}
\begin{equation*}
IM^*_n=\frac{\sum_{l=1}^{8} IM_n^l}{2\sqrt{2}}.
\end{equation*}
Note that the statistics $IM_n^l$ can be independent and their asymptotic properties can be the same as those of the original $IM_n$ and the new test statistic  can also have the same asymptotic properties.

Finally, we also adjust the test by using a factor $\mu_n=1+3n^{-1/2}$ to reduce the magnitude of the test statistic:
\begin{equation*}
\tilde{IM}_n=\frac{IM^*_n}{1+3n^{-1/2}}.
\end{equation*}
Again, it is easy to see that this test statistic has the same asymptotic normality  as the original one under the null hypothesis by using Cram\'er-Wald device and continuous mapping theorem. Without notational confusion, we still write this adjusted version as $IM_n$  in the following.

{\bf Study 1.} Data sets are generated from the following ODE models:

\begin{eqnarray*}	
	&H_{11}:& X^{\prime}(t) = \left[ \begin{array}{c}{\frac{d X_{1}}{d t}} \\ {\frac{d X_{2}}{d t}} \end{array}\right]=\tau \left[ \begin{array}{c}{a X_{1} + 0.4\alpha cos(a X_{1})} \\ {a X_{1} + b X_{2} +0.4\beta cos(a X_{1} + b X_{2})} \end{array}\right],\\	
	&H_{12}: & X^{\prime}(t) = \left[ \begin{array}{c}{\frac{d X_{1}}{d t}} \\ {\frac{d X_{2}}{d t}} \end{array}\right]=\tau \left[ \begin{array}{c}{a X_{1} + 0.1\alpha (a X_{1})^3} \\ {a X_{1} + b X_{2} +0.1\beta (a X_{1} + b X_{2})^3} \end{array}\right],\\
	&H_{13}: & X^{\prime}(t) = \left[ \begin{array}{c}{\frac{d X_{1}}{d t}} \\ {\frac{d X_{2}}{d t}} \end{array}\right]=\tau \left[ \begin{array}{c}{a X_{1} + 2\alpha exp(a X_{1})} \\ {a X_{1} + b X_{2} +5\beta exp(a X_{1} + b X_{2})} \end{array}\right].\\
\end{eqnarray*}

In this study, we consider three different cases in which the linear null ODE models are added with different disturbance terms to form alternative ODE models. The alternatives are oscillating functions of $X$ in $H_{11}$ while they are low-frequent functions of $X$ in both $H_{12}$ and $H_{13}$ under the alternatives. In each case, $\alpha=0$ and $\beta=0$ correspond to the null hypothesis, otherwise  to the alternative hypothesis. When  only one of  $\alpha$ and $\beta$ is nonzero, then  only one element ODE function is different under the alternative hypothesis. When $\alpha$ and $\beta$ are both nonzero, both components are then  changed under the alternative hypothesis. $\tau$ is a timescale parameter which transforms the arbitrary length of sample time interval to $1$. We set the true parameter $(a,b)=(-0.06,-0.24)$, $\tau=10$, $\sigma_{\epsilon}=0.05$ and the sample size is 300. The empirical sizes and powers are presented in Table~1.

The results show that the trajectory matching-based test $TM_{n}$ maintains the significance level when both $\alpha=0$ and $\beta=0$. It also has very good powers under all of the alternative models, which are significantly larger than $IM_n$ and $GM_n$. This is not surprised because $TM$ summarizes the deviation of all the components from the trajectory of null model.


By and large, the integral matching-based tests $IM_{n1}$ and $IM_{n2}$ can maintain the significance level, although in some cases the empirical sizes of $IM_{n1}$ are somewhat lower than the significance level. These departures may originate from the influence of the factor $\mu_n$. $IM_{n1}$ and $IM_{n2}$ have nice power in most settings, showing the effect of these tests.

The proposed tests $GM_{n1}$ and $GM_{n2}$ for the first and second component respectively, when the corresponding component ODE function is under the null hypothesis, can basically maintain the significance level. $GM_{n1}$ has good powers in all three cases while $GM_{n2}$ has varying powers in different cases. This confirms the developed theory.
In the last two cases, we observe that when $(\alpha,\beta)=(1,1)$, $GM_{n2}$ has low powers ($0.600, 0.095$), while when $(\alpha,\beta)=(0,1)$, it has higher powers ($0.734, 0.120$). This phenomenon is worthwhile to pay attention and very different from the classical testing for regressions.  A possible explanation would be that an extra $\alpha$ suppresses the influence of $\beta$ term, making the disturbance term in $\hat{X}^{\prime}(t)$ less important. However, this explanation is not based on any theoretical justification. This anyhow shows the complexity of the testing problem  and is worth of a further investigation.

In the third case, for the first component, $GM_{n1}$ shows greater powers than $IM_{n1}$. However, for the second component, the situation is just on the contrary. This astonishing phenomenon reminds us the  significant difference between matching integral and matching gradient. As is well known in the ODE literature, a relatively little disturbance to the gradient may totally change the form of trajectory of ODE systems while a relatively large disturbance to the gradient may cause little change to the form of trajectory. Thus it is reasonable that $IM_n$ and $GM_n$ have different sensitivity superiority in different settings.

\begin{remark}
	As the first component $X_1'(t)$ in this ODE model only contains $X_1(t)$, the test for this component need not consider the problem of mixed components. Thus, we can also use an adjust version of $TM_n$ which replaces the nonlinear least squares estimator with two-step collocation estimator to check the first component.
\end{remark}

\begin{table*}
	\caption{\linespread{1.5}Empirical sizes and powers in Study 1. n=300, $\alpha=0.05$}
	\small
	\bigskip
	\noindent\makebox[\textwidth][c]{
	\linespread{1.5}\selectfont
	\begin{tabular}{cccccccc}
		\hline
		\textit{Hypothesis} &\textit{$\alpha$} & \textit{$\beta$} & \textbf{$TM_{n}$} & \textit{$IM_{n1}$} & \textbf{$IM_{n2}$} & \textit{$GM_{n1}$} & \textbf{$GM_{n2}$}\\
		\hline
		$H_{11}$ &0 & 0 & 0.045 & 0.028 & 0.047 & 0.038 & 0.050\\
		&   0.5 & 0 & 1.000 & 0.562 &0.048 &0.279 &0.034 \\
		&	0 & 0.5 & 1.000 & 0.025 & 1.000 & 0.042 & 0.191\\
		&	0.5 & 0.5 & 1.000 & 0.555 &1.000 &0.270 &0.193\\
		&   1 & 0 & 1.000 & 1.000 &0.048 &1.000 &0.044 \\
		&	0 & 1 & 1.000 & 0.019 & 1.000 & 0.035 & 0.994\\
		&	1 & 1 & 1.000 & 1.000 &1.000 &1.000 &0.994\\
		$H_{12}$ &0 & 0 & 0.043 & 0.030 & 0.039 & 0.035 & 0.048\\
		&   0.5 & 0& 1.000 & 1.000 &0.046 &0.915 &0.036 \\
		&	0 & 0.5 & 1.000 & 0.022 & 1.000 & 0.039 & 0.661\\
		&	0.5 & 0.5 & 1.000 & 1.000 &1.000 &0.919 &0.606\\
		&   1 & 0 & 1.000 & 1.000 &0.041 &0.998 &0.045 \\
		&	0 & 1 & 1.000 & 0.021 & 0.999 & 0.031 & 0.734\\
		&	1 & 1 & 1.000 & 1.000 &0.996 &0.997 &0.600\\	
		$H_{13}$ &0 & 0 & 0.046 & 0.033 & 0.045 & 0.037 & 0.050\\
		&	0.5 & 0 & 1.000  & 0.115 & 0.043 & 0.294 & 0.034\\
		&	0 & 0.5 & 1.000 & 0.019 & 0.905 & 0.040 & 0.078\\
		&	0.5 & 0.5 & 1.000 &0.131 &0.783 &0.296 &0.062\\
		&   1 & 0 & 1.000 & 0.172 &0.043 &0.906 &0.048 \\
		&	0 & 1 & 1.000 & 0.021 & 0.992 & 0.031 & 0.120\\
		&	1 & 1 & 1.000 & 0.161 &0.925 &0.893 &0.095\\
		\hline
	\end{tabular}}
\end{table*}

{\bf Study 2.} The data sets are generated from the following ODE system:

\begin{eqnarray*}		
	&H_{2}: & X^{\prime}(t) = \left[ \begin{array}{c}{\frac{d X_{1}}{d t}} \\ {\frac{d X_{2}}{d t}} \end{array}\right]=\tau\left[ \begin{array}{c}{a (X_{1}+X_{2}-\frac{X_1^3}{3}) + \alpha X_{1}X_{2}} \\ {-\frac{X_{1} + b X_{2} - c}{a} +0.4\beta X_{1}X_{2}} \end{array}\right].\\
\end{eqnarray*}

\begin{figure}
	\vspace{6pc}
	\centering
	\includegraphics[scale=0.5]{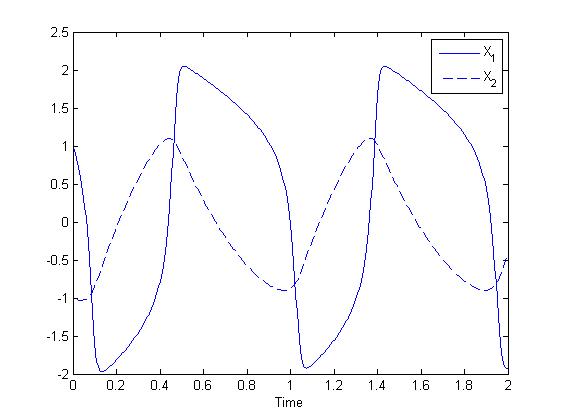}
	\caption[]{Time coarse of FitzHugh-Nagumo model.}
	\label{FN_system.jpg}
\end{figure}

This is the famous FitzHugh-Nagumo ODE system which describes the behavior of spike potentials in the giant axon of squid neurons (\cite{FitzHugh.1961}; \cite{Nagumo.1962}).  Following \cite{Ding.2014}, we set the true parameter $(a,b,c)=(3,0.2,0.34)$, $\tau=10$, $\sigma_{\epsilon}=0.05$, and the initial values $(X_1(0),X_2(0))=(1,-1)$. The sample size is 300. The time coarse of this ODE system is presented in Fig.1. The empirical sizes and powers are reported in Table~2.

The performances of $TM_n$ is still very well for checking this complex nonlinear ODE model. $GM_{n1}$ and $GM_{n2}$  also work well in most cases. Due to the complex interaction between the components of the ODE system, the performances of $GM_{n1}$ and $GM_{n2}$ when $(\alpha,\beta)=(0.5,0.5)$ seem totally different compared to the $(\alpha,\beta)=(1,1)$ setting. These simulation results again show the complexity of the ODE testing problem.

{ Obviously,  $IM_{n1}$ makes no sense at all for the testing. Some unreported results show that  when the sample size is even $10,000$, the empirical size can then be greatly reduced which suggests consistency, but is still too large to make sense.  Together with its performance for testing the linear model above, we must be careful to use $IM_n$ to check complex nonlinear ODE models. $IM_{n2}$ performs acceptably as the hypothetical model is now linear with $(\alpha,\beta)=(0,0)$. }

\begin{table*}
	\caption{\linespread{1.5}Empirical sizes and powers in Study 2. n=300, $\alpha=0.05$}
	\small
	\bigskip
	\noindent\makebox[\textwidth][c]{
	\linespread{1.5}\selectfont
	\begin{tabular}{cccccccc}
		\hline
		\textit{Hypothesis} &\textit{$\alpha$} & \textit{$\beta$} & \textbf{$TM_{n}$} & \textit{$IM_{n1}$} & \textbf{$IM_{n2}$} & \textit{$GM_{n1}$} & \textbf{$GM_{n2}$}\\
		\hline
		$H_{2}$&0 & 0 & 0.048 & 0.866 & 0.069 & 0.071 & 0.048\\
		&   0.5 & 0 & 1.000  & 1.000 &0.063 &0.131 &0.048 \\
		&	0 & 0.5 & 1.000  & 0.881 & 1.000 & 0.055 & 0.159\\
		&	0.5 & 0.5 & 1.000 & 1.000 &1.000 &0.129 &0.203\\
		&   1 & 0 & 1.000 & 1.000 &0.087 &0.514 &0.059 \\
		&	0 & 1 & 1.000 & 0.880 & 1.000 & 0.069 & 0.993\\
		&	1 & 1 & 1.000 & 1.000 &0.217 &0.990 &0.053\\
		\hline
	\end{tabular}}
\end{table*}

{\bf Study 3.} Data sets are generated from the following ODE models:

\begin{eqnarray*}		
	&H_{3}: & X^{\prime}(t) = \left[ \begin{array}{c}{\frac{d X_{1}}{d t}} \\ {\frac{d X_{2}}{d t}} \end{array}\right]=\tau\left[ \begin{array}{c}{a X_{1}+b X_{1}X_{2} + 0.8\alpha X_{2}} \\ {c X_{2}+d X_{1}X_{2} +4\beta X_{1}} \end{array}\right].\\
\end{eqnarray*}

\begin{figure}
	\vspace{6pc}
	\centering
	\includegraphics[scale=0.5]{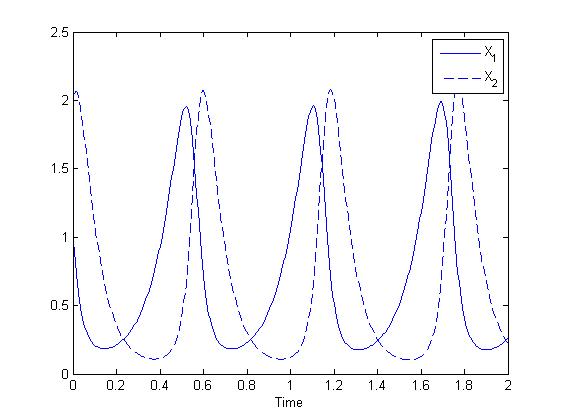}
	\caption[]{Time coarse of Lotka-Volterra model.}
	\label{LV_system.jpg}
\end{figure}

The null ODE system with $\alpha=0$ and $\beta=0$ is the standard Lotka-Volterra model which is well known for modeling the evolution of prey-predator populations (\cite{Lotka.1910}; \cite{Volterra.1928}; \cite{GOEL.1971}). Let the true parameters $(a,b,c,d)=(1,-1.5,-1.5,2)$ and the initial values $(X_1(0),X_2(0))=(1,2)$. The same setting was used in \cite{Brunel.2008} to check the performance of a two-step collocation estimator. We set $\tau=10$, $\sigma_{\epsilon}=0.05$ and the sample size $n =300$. The time coarse of this ODE system is summarized in Fig.~2. The empirical sizes and powers are reported in Table~3.

As can be seen, the performances of tests $TM_n$ and $GM_n$ are similar to those with the models in the last two studies. $IM_{n}$ still fails to maintain the significance level when $(\alpha,\beta)=(0,0)$. We omit the analysis details here.

\begin{table*}
	\caption{\linespread{1.15}Empirical sizes and powers in Study 3. n=300,
		$\alpha=0.05$}
	\label{sphericcase}
	\small
	\bigskip
	\noindent\makebox[\textwidth][c]{
		\linespread{1.5}\selectfont
	\begin{tabular}{cccccccc}
		\hline
		\textit{Hypothesis} &\textit{$\alpha$} & \textit{$\beta$} & \textbf{$TM_{n}$} & \textit{$IM_{n1}$} & \textbf{$IM_{n2}$} & \textit{$GM_{n1}$} & \textbf{$GM_{n2}$}\\
		\hline
		$H_{3}$&0 & 0 & 0.047 & 0.204 & 0.234 & 0.048 & 0.041\\
		&   0.5 & 0 & 1.000 & 0.597 &0.088 &0.113 &0.042 \\
		&	0 & 0.5 & 1.000 & 0.057 & 1.000 & 0.042 & 0.641\\
		&	0.5 & 0.5 & 1.000 & 0.280 &1.000 &0.120 &0.245\\
		&   1 & 0 & 1.000 & 0.255 &0.147 &1.000 &0.042 \\
		&	0 & 1 & 1.000 & 0.043 & 1.000 & 0.042 & 0.821\\
		&	1 & 1 & 1.000 & 0.253 &0.998 &0.095 &0.080\\
		\hline
	\end{tabular}}
\end{table*}

Summarizing the simulation results, we conclude the $TM_n$ and $GM_n$ tests have fine controlled sizes and good powers for extensive ODE systems. $IM_n$ is suitable to be used to check linear ODE system.

\subsection{A real data example}

Now apply our tests to a real data set downloadable from Hulin Wu Lab (https://sph.uth.edu/dotAsset/3ac61148-e59e-493c-bbda-0a38ffe111e5.zip). The data set has been analyzed to show the benefits of differential equation-constrained local polynomial regression for estimating parameters in an ODE model concerning influenza virus-specific effector CD$8$+ T cells (\cite{Ding.2014}). Here we employ the proposed tests to check the adequacy of this model.

The form of the mechanistic ODE system is as follows (\cite{Wu.2011}; \cite{Ding.2014}):
\begin{equation}\label{realdata}
\begin{aligned} \frac{d}{d t} X_{1} &=\tau\left[ \rho_{m} D^{m}(t-\delta_{t})-\delta_{m}-\gamma_{m s}-\gamma_{m l}\right] \\ \frac{d}{d t} X_{2} &=\tau\left[ \rho_{s} D^{s}(t-\delta_{t})-\delta_{s}-\gamma_{s l}+\gamma_{m s} \exp \left(X_{1}-X_{2}\right)\right] \\ \frac{d}{d t} X_{3} &= \tau\left[ \gamma_{m l} \exp \left(X_{1}-X_{3}\right)+\gamma_{s l} \exp \left(X_{2}-X_{3}\right)-\delta_{l}\right] \end{aligned}
\end{equation}
where $X=\left(X_{1}, X_{2}, X_{3}\right)^{\mathrm{T}}=\left(\log \left(T_{E}^{m}\right), \log \left(T_{E}^{s}\right), \log \left(T_{E}^{l}\right)\right)^{\mathrm{T}}$.

$T_{E}^{m}$, $T_{E}^{s}$ and $T_{E}^{l}$ are CD$8$+ T cells among lymph node, spleen and lung respectively. $D^m$ and $D^s$ represent the number of mature dendritic cells in the mediastinal lymph node and spleen separately. As in the simulation part, we add a timescale parameter $\tau=10$ to normalize the sample time interval. The data for $D^m$ is available in the data set and $D^s$ can be replaced by the smoothed estimates of $D^m$. $\theta=(\rho_{m}, \rho_{s}, \delta_{l}, \gamma_{ms}, \gamma_{sl})^{\mathrm{T}}$ is the parameter to be estimated. The other parameters are supposed to be known as $(\delta_{m}, \delta_{s}, \gamma_{ml}, \delta_{t})^{\mathrm{T}}=(0, 0, 0, 3.08)^{\mathrm{T}}$. There are $77$ observations at $9$ distinct time points for each component of $T=\left(T_{E}^{m}, T_{E}^{s}, T_{E}^{l}\right)$ in the data set.

The ODE model (\ref{realdata}) was used to fit data from day $5$ to day $14$, aiming to explain the mechanism of influenza virus-specific effector CD$8$+ T cells. However, to avoid model misspecification we should check the adequacy of this model. Thus, we apply the proposed three tests. We choose the same parameter values as in the simulation part. As the last two component functions contain $X(t)$, to help control the empirical size, we still use the adjusted version of $IM_n$ with the restricted interval $(0.25,0.75)$ and $n_l=2$.

Applying our trajectory-based test with the value $84.10$ of $TM_{n}$, the corresponding $p$-value is abut $0$. This result shows that the whole ODE model under the null is not plausible. Next we use our integral-based test and gradient-based test to check each component function. The values of $IM_{n}$ for the three component functions are $(3.17,2.96,4.04)$ and the p-values are $(0.00077,0.0016,0)$. Since this ODE model is somewhat complicated, the results of $IM_n$ need to be carefully treated. The values of $GM_{n}$ for the three component functions are $(13.44,2.68,25.96)$ and the p-values are $(0,0.0037,0)$. These results suggest all of the three component functions under the null are not tenable. Thus we may need to modify the models to fit the data if necessary.

\begin{figure}
	\vspace{6pc}
	\centering
	\includegraphics[scale=0.2]{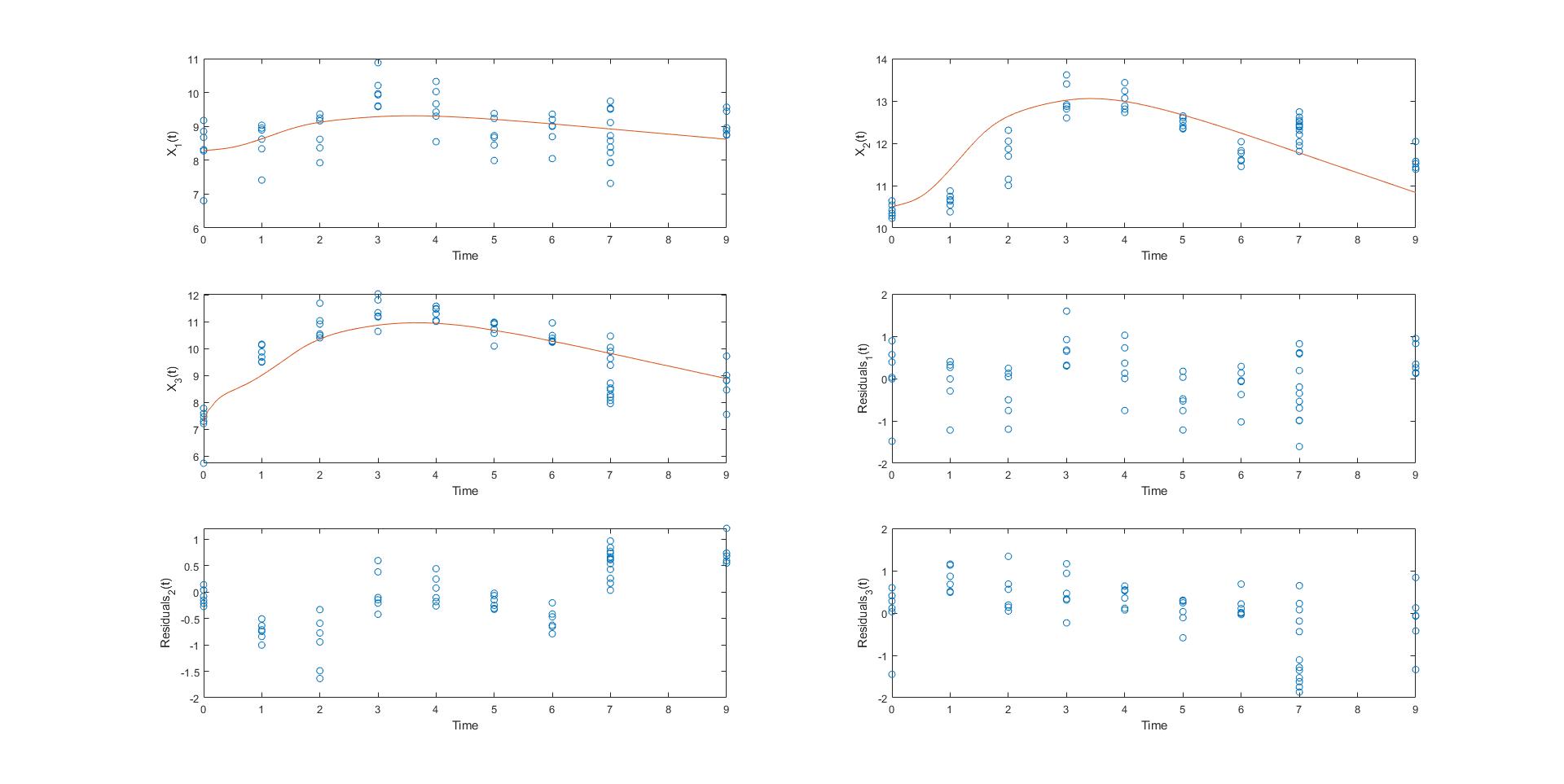}
	\caption[]{Time coarse of response and residuals.}
	\label{residual}
\end{figure}

\section{Conclusion}

In this paper, we investigate  model checking for parametric ordinary differential equations system and propose three tests to respectively check the whole system and their components. The trajectory matching-based test is for the whole ODE system and the other two integral matching-based and gradient matching-based tests for every single component function in ODEs. The tests can detect global as well as local alternatives.

There are four issues worthwhile to investigate in the future studies. First, due to the complicated structure, the tests involve delicately selected bandwidths that affect the performances of the tests. As briefly mentioned in Section~1,  we do not apply the idea of global smoothing test to construct a test for this problem. From Fig. 1 and Fig. 2, we can see that some famous ODE models are very highly oscillating and thus local smoothing test using nonparametric estimation may more sensitively capture local departures of ODE models. But, it deserves a study to see whether global smoothing test could be more powerful for low frequency ODE models. Second, we can see that due to any single component of the ODE system actually shares the same original function $X(\cdot)$, the corresponding tests are expectably not very powerful as we commented in Section~1. How to solve this problem is a big issue. Third, as seen in simulation, $IM_n$ is hard to control the significance level due to its sensitivity to the nonparametric estimator. How to modify it is a nontrivial task. Fourth,  for ODE models, it is also the case where the ODE system is large, that is, $p$ is large. This is a very challenging problem in effect. The research is ongoing.

\bibliographystyle{apalike}
\bibliography{document}

\begin{thebibliography}{}

\bibitem[Brunel, 2008]{Brunel.2008}
Brunel, N. J.-B. (2008).
\newblock Parameter estimation of ode's via nonparametric estimators.
\newblock {\em Electronic Journal of Statistics}, 2(0):1242--1267.

\bibitem[Chen et~al., 2017]{Chen.2017}
Chen, S., Shojaie, A., and Witten, D.~M. (2017).
\newblock Network reconstruction from high-dimensional ordinary differential
  equations.
\newblock {\em Journal of the American Statistical Association},
  112(520):1697--1707.

\bibitem[Currie and Wilson, 2012]{Currie.2012}
Currie, J. and Wilson, D.~I. (2012).
\newblock Opti: lowering the barrier between open source optimizers and the
  industrial matlab user.

\bibitem[Dattner and Klaassen, 2015]{Dattner.2015}
Dattner, I. and Klaassen, C. A.~J. (2015).
\newblock Optimal rate of direct estimators in systems of ordinary differential
  equations linear in functions of the parameters.
\newblock {\em Electronic Journal of Statistics}, 9(2):1939--1973.

\bibitem[Dette, 1999]{Dette.1999}
Dette, H. (1999).
\newblock A consistent test for the functional form of a regression based on a
  difference of variance estimators.
\newblock {\em The Annals of Statistics}, 27(3):1012--1040.

\bibitem[Ding and Wu, 2014]{Ding.2014}
Ding, A.~A. and Wu, H. (2014).
\newblock Estimation of ordinary differential equation parameters using
  constrained local polynomial regression.
\newblock {\em Statistica Sinica}, 24(4):1613--1631.

\bibitem[Fang, 2010]{Fang.October2010}
Fang, Y. (October, 2010).
\newblock {\em Several studies of models and methods in biostatistics
  (Chinese)}.
\newblock Disserttation for doctor degree, Shanghai, China.

\bibitem[FitzHugh, 1961]{FitzHugh.1961}
FitzHugh, R. (1961).
\newblock Impulses and physiological states in theoretical models of nerve
  membrane.
\newblock {\em Biophysical Journal}, 1(6):445--466.

\bibitem[Goel et~al., 1971]{GOEL.1971}
Goel, N.~S., Maitra, S.~C., and Montroll, E.~W. (1971).
\newblock On the volterra and other nonlinear models of interacting
  populations.
\newblock {\em Reviews of Modern Physics}, 43(2):231--276.

\bibitem[Gonz{\'a}lez-Manteiga and Crujeiras, 2013]{GONZALEZMANTEIGA.2013}
Gonz{\'a}lez-Manteiga, W. and Crujeiras, R.~M. (2013).
\newblock An updated review of goodness-of-fit tests for regression models.
\newblock {\em Test}, 22(3):361--411.

\bibitem[Gugushvili and Klaassen, 2012]{Gugushvili.2012}
Gugushvili, S. and Klaassen, C.~A. (2012).
\newblock $\sqrt{n}$ consistent parameter estimation for systems of ordinary
  differential equations: bypassing numerical integration via smoothing.
\newblock {\em Bernoulli}, 18(3):1061--1098.

\bibitem[Hansen, 2008]{Hansen.2008}
Hansen, B.~E. (2008).
\newblock Uniform convergence rates for kernel estimation with dependent data.
\newblock {\em Econometric Theory}, 24(3):726--748.

\bibitem[H{\"a}rdle and Mammen, 1993]{Hardle.1993}
H{\"a}rdle, W. and Mammen, E. (1993).
\newblock Comparing nonparametric versus parametric regression fits.
\newblock {\em The Annals of Statistics}, 21(4):1926--1947.

\bibitem[Koul and Ni, 2004]{Koul.2004}
Koul, H.~L. and Ni, P. (2004).
\newblock Minimum distance regression model checking.
\newblock {\em Journal of Statistical Planning and Inference}, 119(1):109--141.

\bibitem[Liang and Wu, 2008]{Liang.2008}
Liang, H. and Wu, H. (2008).
\newblock Parameter estimation for differential equation models using a
  framework of measurement error in regression models.
\newblock {\em Journal of the American Statistical Association},
  103(484):1570--1583.

\bibitem[Lotka, 1910]{Lotka.1910}
Lotka, A.~J. (1910).
\newblock Contribution to the theory of periodic reactions.
\newblock {\em The Journal of Physical Chemistry}, 14(3):271--274.

\bibitem[Maechler, 2010]{MaechlerMartin.2010}
Maechler, M. (2010).
\newblock Kernel regression smoothing with adaptive local or global plug-in
  bandwidth selection.

\bibitem[Martins-Filho and Yao, 2006]{MartinsFilho.2006}
Martins-Filho, C. and Yao, F. (2006).
\newblock A note on the use of v and u statistics in nonparametric models of
  regression.
\newblock {\em Annals of the Institute of Statistical Mathematics},
  58(2):389--406.

\bibitem[Nagumo et~al., 1962]{Nagumo.1962}
Nagumo, J., Arimoto, S., and Yoshizawa, S. (1962).
\newblock An active pulse transmission line simulating nerve axon.
\newblock {\em Proceedings of the IRE}, 50(10):2061--2070.

\bibitem[Ramsay and Hooker, 2017]{Ramsay.2017}
Ramsay, J. and Hooker, G. (2017).
\newblock {\em Dynamic data analysis}.
\newblock {Springer New York}, New York, NY.

\bibitem[Stute, 1997]{Stute.1997}
Stute, W. (1997).
\newblock Nonparametric model checks for regression.
\newblock {\em The Annals of Statistics}, 25(2):613--641.

\bibitem[Stute et~al., 1998a]{Stute.1998a}
Stute, W., Manteiga, W.~G., and Quindimil, M.~P. (1998a).
\newblock Bootstrap approximations in model checks for regression.
\newblock {\em Journal of the American Statistical Association},
  93(441):141--149.

\bibitem[Stute et~al., 1998b]{Stute.1998b}
Stute, W., Thies, S., and Zhu, L.-X. (1998b).
\newblock Model checks for regression: an innovation process approach.
\newblock {\em The Annals of Statistics}, 26(5):1916--1934.

\bibitem[Tan and Zhu, 2019]{Tan.2019}
Tan, F. and Zhu, L.-X. (2019).
\newblock Adaptive-to-model checking for regressions with diverging number of
  predictors.
\newblock {\em The Annals of Statistics}, 47(4):1960--1994.

\bibitem[Volterra, 1928]{Volterra.1928}
Volterra, V. (1928).
\newblock Variations and fluctuations of the number of individuals in animal
  species living together.
\newblock {\em ICES Journal of Marine Science}, 3(1):3--51.

\bibitem[Wu et~al., 2011]{Wu.2011}
Wu, H., Kumar, A., Miao, H., Holden-Wiltse, J., Mosmann, T.~R., Livingstone,
  A.~M., Belz, G.~T., Perelson, A.~S., Zand, M.~S., and Topham, D.~J. (2011).
\newblock Modeling of influenza-specific cd8+ t cells during the primary
  response indicates that the spleen is a major source of effectors.
\newblock {\em Journal of immunology}, 187(9):4474--4482.

\bibitem[Xue et~al., 2010]{Xue.2010}
Xue, H., Miao, H., and Wu, H. (2010).
\newblock Sieve estimation of constant and time-varying coefficients in
  nonlinear ordinary differential equation models by considering both numerical
  error and measurement error.
\newblock {\em The Annals of Statistics}, 38(4):2351--2387.

\bibitem[Zheng, 1996]{Zheng.1996}
Zheng, J.~X. (1996).
\newblock A consistent test of functional form via nonparametric estimation
  techniques.
\newblock {\em Journal of Econometrics}, 75(2):263--289.

\bibitem[Zhu, 2003]{Lixing.2003}
Zhu, L.-X. (2003).
\newblock Model checking of dimension-reduction type for regression.
\newblock {\em Statistica Sinica}, 13(2):283--296.

\bibitem[Zhu and Li, 1998]{Lixing.1998}
Zhu, L.-X. and Li, R. (1998).
\newblock Dimension-reduction type test for linearity of a stochastic
  regression model.
\newblock {\em Acta Mathematicae Applicatae Sinica}, 14(2):165--175.

\end{thebibliography}


\begin{thebibliography}{}

\bibitem[Hall, 1984]{Hall.1984}
Hall, P. (1984).
\newblock Central limit theorem for integrated square error of multivariate
  nonparametric density estimators.
\newblock {\em Journal of Multivariate Analysis}, 14(1):1--16.

\bibitem[Hansen, 2008]{Hansen.2008}
Hansen, B.~E. (2008).
\newblock Uniform convergence rates for kernel estimation with dependent data.
\newblock {\em Econometric Theory}, 24(3):726--748.

\bibitem[Li et~al., 2019]{Li.2019}
Li, L., Chiu, S.~N., and Zhu, L.-X. (2019).
\newblock Model checking for regressions: An approach bridging between local
  smoothing and global smoothing methods.
\newblock {\em Computational Statistics {\&} Data Analysis}, 138:64--82.

\bibitem[Liang and Wu, 2008]{Liang.2008}
Liang, H. and Wu, H. (2008).
\newblock Parameter estimation for differential equation models using a
  framework of measurement error in regression models.
\newblock {\em Journal of the American Statistical Association},
  103(484):1570--1583.

\bibitem[Martins-Filho and Yao, 2006]{MartinsFilho.2006}
Martins-Filho, C. and Yao, F. (2006).
\newblock A note on the use of v and u statistics in nonparametric models of
  regression.
\newblock {\em Annals of the Institute of Statistical Mathematics},
  58(2):389--406.

\bibitem[Powell et~al., 1989]{Powell.1989}
Powell, J.~L., Stock, J.~H., and Stoker, T.~M. (1989).
\newblock Semiparametric estimation of index coefficients.
\newblock {\em Econometrica}, 57(6):1403.

\bibitem[Zheng, 1996]{Zheng.1996}
Zheng, J.~X. (1996).
\newblock A consistent test of functional form via nonparametric estimation
  techniques.
\newblock {\em Journal of Econometrics}, 75(2):263--289.

\end{thebibliography}

\end{document}


\begin{center}
{\bf\large Supplementary Material to ``Model Checking for Parametric Ordinary Differential Equations System"}
\\ \vskip0.6cm
Ran Liu$^1$, Yun Fang$^2$ and Lixing Zhu$^{1}$
\\
\textit{ $^1$ Hong Kong Baptist University\\
$^2$ Shanghai Normal University}

\today
\end{center}

\section{Additional assumptions}

Other than the assumptions given in each lemma and theorem, we give some other assumptions here. These assumptions are divided to three sets. Assumptions in set A give the basic setting of $(t,Y)$ and the conditions on kernel function. These assumptions assume ensure the uniform convergence rate for kernel estimation. Assumptions in set B place restrictions on the primitive function $F(t;\theta)$ which include some conditions needed for the nonlinear least squares estimation. Assumptions in set C contain the conditions on $f(t,X(t);\theta)$ that are necessary for the two-step collocation method.\\

Set A.

1. $t_i$ are i.i.d. random samples and have a common compact support $[t_0,T]$. The density function $p(t)$ is bounded and bounded away from $0$. The first and second derivative of $p(t)$ are bounded and continuous.

2. For all $1\leq k \leq p$, $E(y_{ik}^4|t_i)$ is continuously differentiable and bounded by a measurable function $b(t)$ with $E(b^2(t_i))< \infty$. Furthermore, there exists $s>2$ such that
$$E(|y_{ik}|^s)<\infty,$$
$$\sup_{t}E(|y_{ik}|^s|t_i=t)p(t) \leq \Lambda_{1}.$$

3. The kernel function $K(u)$ is a nonnegative, bounded, continuous, symmetric function and is supported on $[-1,1]$ with $\int K(u) d u=1$. For all $u, u^{\prime} \in R$, $\left|K(u)-K\left(u^{\prime}\right)\right| \leq \Lambda_{2}\left\|u-u^{\prime}\right\|$ for some $\Lambda_{2}< \infty$.

4. The parameter space $\Theta$ is a closed, convex, bounded compact subset of $R^q$.

\

Set B.

1. $F_k(t;\theta)$ is a Borel measurable real function on $\mathrm{R}^p$ for each $\theta$ and is twice continuously differentiable with respect to $\theta$ for each $t$.

2. For all $1\leq k \leq p$,
\begin{eqnarray*}	
	&E&\left[\sup _{\theta \in \Theta} F_k^{2}(t; \theta)\right]<\infty,\\
	&E&\left[ \sup _{\theta \in \Theta}\left\|\frac{\partial F_k(t;\theta)}{\partial \theta} \frac{\partial F_k(t; \theta)}{\partial \theta^{\mathrm{T}}}\right\|\right]<\infty,\\
	&E&\left\{\sup _{\theta \in \Theta}\left\|[Y_k-F_k(t; \theta)]^{2} \frac{\partial F_k(t; \theta)}{\partial \theta} \frac{\partial F_k(t; \theta)}{\partial \theta^{\mathrm{T}}}\right\|\right\}<\infty,\\
	&E&\left\{\sup _{\theta \in \Theta}\left\|[Y_k-F_k(t; \theta)]^{2} \frac{\partial^{2} F_k(t; \theta)}{\partial \theta \partial \theta^{\mathrm{T}}}\right\|\right\}<\infty.
\end{eqnarray*}

3. $\mathrm{E}\left[\sum_{k=1}^{p}\left(Y_{ki}-F_k\left(t_{i}; \theta\right)\right)^{2}\right] \text { takes a unique minimum at } \theta^* \in \Theta$.

4. The matrix $E\left[\sum_{k=1}^{p}\frac{\partial F_k\left(t_{i}; \theta_{0}\right)}{\partial \theta} \frac{\partial F_k\left(t_{i};\theta_{0}\right)}{\partial \theta^{\mathrm{T}}}\right]$ is nonsingular.


\

Set C.

1. The function $X_k^{(3)}(t)$ is continuous on $[t_0,T]$.


2. $f_k(t, X(t); \theta)$ is a continuous function of $\theta$ for $\theta \in \Omega_{\theta}$.

3. $E_{p^*}\left\{ \left[X_k^{\prime}(t)-f_k\left(t, X(t), \theta\right)\right]^{2}w_k(t)\right\}$ takes a unique miminum at $\theta^* \in \Theta$.

4. The first and second partial derivatives, $\frac{\partial f_k(t, X(t); \theta)}{\partial \theta}$, $\frac{\partial^2 f_k(t, X(t); \theta)}{\partial X \partial \theta}$, and $\frac{\partial^2 f_k(t, X(t); \theta)}{\partial \theta \partial \theta^{\mathrm{T}}}$, exist and are continuous for all $\theta \in \Theta$, $X \in \mathcal{X}$, and
$$\left\|\frac{\partial f_k(t, X_1(t); \theta)}{\partial \theta}-\frac{\partial f_k(t, X_2(t); \theta)}{\partial \theta}\right\| \leq C_1 \left\| X_1-X_2 \right\|^{\zeta}$$
for some $0 \leq \zeta \leq 1$.

5. The first partial derivative $\frac{\partial f_k(t, X(t); \theta)}{\partial X}$ is continuous for $X \in \mathcal{X}$ and satisfies
$$\sup_{X \in \mathcal{X}} \left\| \frac{\partial f_k(t, X(t); \theta)}{\partial X} \right\| \leq M(t;\theta).$$

6. If $t^*_i$ is randomly designed,   its density function $p^*(t)$ is bounded away from $0$ and has  bounded and continuous first derivative on $[t_0,T]$. If $t^*_i$ is in fixed design,  there exists a distribution $P^*(t)$ with  the corresponding density function $p^*(t)$ satisfying the above conditions such that
$$\sup_{t}|P_m^*(t)-P^*(t)|=O(m^{-1})$$
where $P_m^*(t)$ is the empirical distribution of $(t_1^*,\cdots,t_m^*)$.

\section{Remark of notations}
In the following proofs, we omit the corresponding superscripts $F,\hat{F},f$ and subscripts $NLS$, $TS$ for simplicity. The notations $V_i$ and $S_i$ will  present statistics used in the proofs, which may have different meanings for each appearance.

\section{Preliminary Lemmas}

Before giving the proofs of Theorems and Lemmas, we provide some results about the uniform convergence rate of kernel estimation and local polynominal estimation as preliminary lemmas. These preliminary lemmas are useful for the proofs of the lemmas and theorems. The proof of Lemma~\ref{LP.1}  can be founded in \cite{Hansen.2008}, the proofs of Lemma~\ref{LP.2} and Lemma~\ref{LP.3} can be extended to the vector version under the Frobenius norm $\| \cdot \|$  by the proofs in \cite{Hansen.2008} and \cite{Liang.2008}.

\setcounter{lemma}{9}
\begin{lemma}
	(\cite{Hansen.2008}) Under sets A-C of assumptions in Appendix A, $\frac{\ln n}{n h}=o(1)$, for the kernel density estimator $\hat{p}(t)$, we have
	\begin{equation}
	\sup_{t} |\hat{p}(t)-p(t)|=\mathrm{o}_{P}(a_n)
	\end{equation}
	where $a_n=h^2+n^{-1/2}h^{-1/2}\log{n}^{-1/2}$.\label{LP.1}
\end{lemma}

\begin{lemma}
	Under sets A-C of assumptions in Appendix A, $\frac{\ln n}{n h}=o(1)$, for the local linear estimator,
	\begin{equation}
	\begin{aligned}
	\hat{X}(t)=&\frac{\frac{1}{n^2h^2}\sum_k^n \sum_l^n \left[(\frac{t-t_k}{h})^2K(\frac{t-t_k}{h})K(\frac{t-t_l}{h})Y_l - \frac{t-t_k}{h}K(\frac{t-t_k}{h}) \frac{t-t_l}{h}K(\frac{t-t_l}{h})Y_l \right]}{\frac{1}{n^2h^2}\sum_k^n \sum_l^n \left[(\frac{t-t_k}{h})^2K(\frac{t-t_k}{h})K(\frac{t-t_l}{h}) - \frac{t-t_k}{h}K(\frac{t-t_k}{h}) \frac{t-t_l}{h}K(\frac{t-t_l}{h})\right]}\\
	\equiv &\frac{N_n(t)}{M_n(t)},
	\end{aligned} \label{LPE.1}
	\end{equation}
	we have
	\begin{equation}
	\sup_{t} |M_n(t)-M(t)|=\mathrm{o}_{P}(a_n),
	\end{equation}
	with $M(t)=\int u^2 K(u) \mathrm{d}u \cdot p^2(t)=\mu_{2}(K)p^2(t)$, and
	\begin{equation}
	\sup_{t} \|\hat{X}(t)-X(t)\|=\mathrm{o}_{P}(a_n)
	\end{equation}
	where $a_n=h^2+n^{-1/2}h^{-1/2}\log{n}^{-1/2}$.\label{LP.2}
\end{lemma}

\begin{lemma}
	Under sets A-C of assumptions in Appendix A, $\frac{\ln n}{n h^3}=o(1)$, for the local quadratic estimator $\hat{X}^{\prime}(t)$, we have
	\begin{equation}
	\sup_{t} \|\hat{X}^{\prime}(t)-X^{\prime}(t)\|=\mathrm{o}_{P}(b_n)
	\end{equation}
	where $a_n=h^2+n^{-1/2}h^{-3/2}\log{n}$.\label{LP.3}
\end{lemma}


\section{The results in Section~2}
\subsection{Proof of Theorem~2.1}

\begin{proof}
	This theorem can be regarded as a straightforward extension of Lemma~3 of \cite{Li.2019} to the multi-response case. Since the proof is similar,  we omit it here.
\end{proof}

\subsection{Proof of Lemma~1}

\begin{proof}
	For every component $k$, we decompose $V_{nk}$ into three terms:
	\begin{equation}
	\begin{aligned}
	V_{nk}=&\frac{1}{n(n-1)} \sum_{i=1}^{n} \sum_{j=1 \atop j \neq i}^{n} \frac{1}{h} K\left(\frac{t_{i}-t_{j}}{h}\right) e_{ik} e_{jk}\\
	=&\left[ \frac{1}{n(n-1)} \sum_{i=1}^{n} \sum_{j=1 \atop j \neq i}^{n} \frac{1}{h} K\left(\frac{t_{i}-t_{j}}{h}\right) \varepsilon_{ik} \varepsilon_{jk}\right]\\
	&-2\left\{\frac{1}{n(n-1)} \sum_{i=1}^{n} \sum_{j=1 \atop j \neq i}^{n} \frac{1}{h} K\left(\frac{t_{i}-t_{j}}{h}\right) \varepsilon_{ik}\left[F_k\left(t_{j}, \hat{\theta}\right)-F_k\left(t_{j}, \theta_{0}\right)\right]\right\}\\
	&+\left\{\frac{1}{n(n-1)} \sum_{i=1}^{n} \sum_{j=1 \atop j \neq i}^{n} \frac{1}{h} K\left(\frac{t_{i}-t_{j}}{h}\right)\right.\\
	&\left.\times\left[F_k\left(t_{j}, \hat{\theta}\right)-F_k\left(t_{j}, \theta_{0}\right)\right]\left[F_k\left(t_{j}, \hat{\theta}\right)-F_k\left(t_{j}, \theta_{0}\right)\right]\right\}\\
	\equiv & V_{1 n k}-2 V_{2 n k}+V_{3 n k}.\label{A1}
	\end{aligned}
	\end{equation}
	
	Therefore  the vector $V_{n}$ can be written as $V_{1 n}-2 V_{2 n}+V_{3 n}$. We now show that $n h^{1 / 2} V_{1 n} \stackrel{\mathbf{d}}{\longrightarrow} \mathrm{N}(0, \Sigma)$ while $n h^{1 / 2} V_{2 n}$ and $n h^{1 / 2} V_{3 n}$ converge to zero in probability.
	
	To prove that $n h^{1 / 2} V_{1 n} \stackrel{\mathbf{d}}{\longrightarrow} \mathrm{N}(0, \Sigma)$, we only need to verify that for every $\lambda \in R^p$, $n h^{1 / 2} \lambda^{\mathrm{T}} V_{1 n} \stackrel{\mathbf{d}}{\longrightarrow} \mathrm{N}(0, \lambda^{\mathrm{T}} \Sigma \lambda)$ according to the Cram\'er-Wald device. To confirm this statement, write $\lambda^{\mathrm{T}} V_{1 n}$ in a U-statistic form with the kernel:
	\begin{equation*}
	\tilde{H}_{n}\left(z_{i}, z_{j}\right)=\sum_{k=1}^{p}\frac{\lambda_k}{h} K\left(\frac{t_{i}-t_{j}}{h}\right) \varepsilon_{ik} \varepsilon_{jk}
	\end{equation*}
	where $z_i=(t_i,\varepsilon_i)$. Since it is a one-dimensional degenerate U-statistic,  Theorem~1 of \cite{Hall.1984} can be applied to obtain its asymptotic distribution. We then verity the conditions of this theorem. To this end, we have the following equations:	
	\begin{equation*}
	\begin{aligned}
	&\mathrm{E}\left[L_{n}^{2}\left(z_{1}, z_{2}\right)\right]\\ =&\mathrm{E}\left\{\mathrm{E}\left[\tilde{H}_{n}\left(z_{3}, z_{1}\right) \tilde{H}_{n}\left(z_{3}, z_{2}\right) | z_{1}, z_{2}\right]\right\}^{2} \\
	=&\mathrm{E}\left\{\mathrm{E}\left[\sum_{k_1=1}^{p}\sum_{k_2=1}^{p}\frac{\lambda_{k_1}\lambda_{k_2}}{h^{2}} K\left(\frac{t_{3}-t_{1}}{h}\right) K\left(\frac{t_{3}-t_{2}}{h}\right) \varepsilon_{1k_1} \varepsilon_{2k_2} \varepsilon_{3k_1}\varepsilon_{3k_2} | z_{1}, z_{2}\right]\right\}^{2}\\
	=&\frac{1}{h^{4}} \mathrm{E}\left\{\sum_{k_1=1}^{p}\sum_{k_2=1}^{p}\lambda_{k_1}\lambda_{k_2}\varepsilon_{1k_1} \varepsilon_{2k_2} \mathrm{E}\left[K\left(\frac{t_{3}-t_{1}}{h}\right) K\left(\frac{t_{3}-t_{2}}{h}\right) \sigma_{k_1k_2}\left(t_{3}\right) | t_{1}, t_{2}\right]\right\}^{2} \\
	=&\frac{1}{h^{4}} \mathrm{E}\left[ \sum_{k_1=1}^{p}\sum_{k_2=1}^{p}\lambda_{k_1}\lambda_{k_2}\varepsilon_{1k_1} \varepsilon_{2k_1}  \int K\left(\frac{t_{3}-t_{1}}{h}\right) K\left(\frac{t_{3}-t_{2}}{h}\right) \sigma_{k_1k_2}\left(t_{3}\right) p\left(t_{3}\right) \mathrm{d} t_{3}\right]^{2}\\
	=&\frac{1}{h^{4}} \mathrm{E}\left[ \sum_{k_1=1}^{p}\sum_{k_2=1}^{p}\lambda_{k_1}\lambda_{k_2}\varepsilon_{1k_1} \varepsilon_{2k_1}  \int K\left(u\right) K\left(u+\frac{t_{1}-t_{2}}{h}\right) \right.\\
	&\left. \times \sigma_{k_1k_2}\left(t_{1}+h u\right) p\left(t_{1}+h u\right) h \mathrm{d} u \right]^{2}\\
	=&\frac{1}{h^{2}} \mathrm{E}\left\{\mathrm{E} \left\{\sum_{k_1=1}^{p}\sum_{k_2=1}^{p} \sum_{k_3=1}^{p} \sum_{k_4=1}^{p} \lambda_{k_1}\lambda_{k_2}\lambda_{k_3}\lambda_{k_4}\varepsilon_{1k_1} \varepsilon_{2k_2} \varepsilon_{1k_3} \varepsilon_{2k_4} \right.\right.\\
	& \times \left[\int K(u) K\left(u+\frac{t_{1}-t_{2}}{h}\right)\sigma_{k_1k_2}\left(t_{1}+h u\right) p\left(t_{1}+h u\right) d u \right]\\
	& \left. \left. \times \left[\int K(u) K\left(u+\frac{t_{1}-t_{2}}{h}\right)\sigma_{k_3k_4}\left(t_{1}+h u\right) p\left(t_{1}+h u\right) d u \right] | t_{1}, t_{2} \right\} \right\}\\
	=&\frac{1}{h^{2}} \mathrm{E} \left\{ \sum_{k_1=1}^{p}\sum_{k_2=1}^{p} \sum_{k_3=1}^{p} \sum_{k_4=1}^{p}\lambda_{k_1}\lambda_{k_2}\lambda_{k_3}\lambda_{k_4} \sigma_{k_1k_3}\left(t_{1}\right) \sigma_{k_2k_4}\left(t_{2}\right) \right.  \\
	&\times \left[\int K(u) K\left(u+\frac{t_{1}-t_{2}}{h}\right)\sigma_{k_1k_2}\left(t_{1}+h u\right) p\left(t_{1}+h u\right) d u \right]\\
	&\left. \times \left[\int K(u) K\left(u+\frac{t_{1}-t_{2}}{h}\right) \sigma_{k_3k_4}\left(t_{1}+h u\right) p\left(t_{1}+h u\right) d u \right] \right\} \\
	\leq & \frac{p^4 \lambda_{max}^4}{h^{2}} \int \sigma^{2}_{max}\left[\int K(u) K\left(u+\frac{t_{1}-t_{2}}{h}\right)\right. \\
	&\times \sigma_{max} p\left(t_{1}+h u\right) \mathrm{d} u ]^{2} p\left(t_{1}\right) p\left(t_{2}\right) \mathrm{d} t_{1} \mathrm{d} t_{2}\\
	=&\frac{p^4 \lambda_{max}^4}{h^{2}} \sigma_{max}^2 \int \left[\int K(u) K(u+v)\right.\\
	&\times \sigma_{max} p\left(t_{1}+h u\right) \mathrm{d} u ]^{2} h p\left(t_{1}\right) p\left(t_{1}-h v\right) \mathrm{d} t_{1} \mathrm{d} v \\
	=& \frac{p^4 \lambda_{max}^4}{h} \sigma_{max}^4\int\left[\int K(u) K(u+v) \mathrm{d} u\right]^{2} \mathrm{d} v \int p^{4}(t) \mathrm{d} t+o\left(1 / h\right) \\
	=& \mathrm{O}\left(1 / h\right).
	\end{aligned}
	\end{equation*}
	
	\begin{equation*}
	\begin{aligned}
	&\mathrm{E}\left[\tilde{H}_{n}^{2}\left(z_{1}, z_{2}\right)\right]\\ &=\mathrm{E}\left\{\mathrm{E}\left[\tilde{H}_{n}^{2}\left(z_{1}, z_{2}\right) | t_{1}, t_{2}\right]\right\} \\
	&=\int \frac{1}{h^{2}}\sum_{k_1=1}^{p}\sum_{k_2=1}^{p} \lambda_{k_1}\lambda_{k_2} K^{2}\left(\frac{t_{1}-t_{2}}{h}\right) \sigma_{k_1k_2}\left(t_{1}\right) \sigma_{k_1k_2}\left(t_{2}\right) p\left(t_{1}\right) p\left(t_{2}\right) \mathrm{d} t_{1} \mathrm{d} t_{2} \\
	&=\frac{1}{h^{2}}\sum_{k_1=1}^{p} \sum_{k_2=1}^{p} \lambda_{k_1}\lambda_{k_2}\int K^{2}(u) \sigma_{k_1k_2}(t_1) \sigma_{k_1k_2}(t_1-h u) p(t) p(t-h u)h \mathrm{d} t \mathrm{d} u \\
	&=\frac{1}{h}\sum_{k_1=1}^{p} \sum_{k_2=1}^{p} \lambda_{k_1}\lambda_{k_2}\int K^{2}(u) \mathrm{d} u \int\left[\sigma_{k1k2}(t)\right]^{2} p^{2}(t) \mathrm{d} t+\mathrm{o}\left(1 / h\right)\\
	&=\mathrm{O}\left(1 / h\right).
	\end{aligned}
	\end{equation*}
	
	Also	
	\begin{equation*}
	\begin{aligned}
	&\mathrm{E}\left[\tilde{H}_{n}^{4}\left(z_{1}, z_{2}\right)\right] \\
	=&\frac{1}{h^{4}}\int  K^{4}\left(\frac{t_{1}-t_{2}}{h}\right) \sum_{k_1=1}^{p}\sum_{k_2=1}^{p}\sum_{k_3=1}^{p}\sum_{k_4=1}^{p} \lambda_{k_1}\lambda_{k_2}\lambda_{k_3}\lambda_{k_4} \{ \mathrm{E}\left[ \varepsilon_{1k_1}\varepsilon_{1k_2}\varepsilon_{1k_3}\varepsilon_{1k_4} | t_1\right]\\
	&\times \mathrm{E}\left[ \varepsilon_{2k_1}\varepsilon_{2k_2}\varepsilon_{2k_3}\varepsilon_{2k_4} | t_2\right] \} p\left(t_{1}\right) p\left(t_{2} \right) \mathrm{d} t_{1} \mathrm{d} t_{2} \\
	=&\frac{1}{h^{4}} \{\sum_{k_1=1}^{p}\sum_{k_2=1}^{p}\sum_{k_3=1}^{p}\sum_{k_4=1}^{p}\lambda_{k_1}\lambda_{k_2}\lambda_{k_3}\lambda_{k_4} \int K^{4}(u) \{ \sigma_{k_1k_2k_3k_4}(t_1)\\
	&\times \sigma_{k_1k_2k_3k_4}(t_1-hu) \} p(t_1) p(t_1-hu) h\mathrm{d} t_1 \mathrm{d} u \} \\
	=&\mathrm{O}\left(1 / h^{3}\right).
	\end{aligned}
	\end{equation*}
	
	From these equations, we have
	\begin{equation*}
	\begin{aligned}
	\frac{\mathrm{E}\left[L_{n}^{2}\left(z_{1}, z_{2}\right)\right]+n^{-1} \mathrm{E}\left[\tilde{H}_{n}^{4}\left(z_{1}, z_{2}\right)\right]}{\left\{\mathrm{E}\left[\tilde{H}_{n}^{2}\left(z_{1}, z_{2}\right)\right]\right\}^{2}} &=\frac{\mathrm{O}\left(1 / h\right)+n^{-1} \mathrm{O}\left(1 / h^{3}\right)}{\mathrm{O}\left(1 / h^{2}\right)} \\ &=\mathrm{O}\left(h\right)+\mathrm{O}\left(1 / n h\right) \longrightarrow 0 .
	\end{aligned}
	\end{equation*}
	
	Since the conditions in Theorem~1 of \cite{Hall.1984}  are verified, we then have
	\begin{equation*}
	n \lambda^{\mathrm{T}} \cdot V_{1 n} /\left\{2 \mathrm{E}\left[\tilde{H}_{n}^{2}\left(z_{i}, z_{j}\right)\right]\right\}^{1 / 2} \stackrel{\mathbf{d}}{\longrightarrow} \mathrm{N}(0,1).
	\end{equation*}
	
	This implies that
	\begin{equation}
	n h^{1 / 2} \lambda^{\mathrm{T}} V_{1 n} \stackrel{\mathbf{d}}{\longrightarrow} \mathrm{N}\left(0, 2\sum_{k_1=1}^{p} \sum_{k_2=1}^{p} \lambda_{k_1}\lambda_{k_2}\int K^{2}(u) \mathrm{d} u \cdot \int\left[\sigma_{k1k2}(t)\right]^{2} p^{2}(t) \mathrm{d} t\right).\label{A2}
	\end{equation}
	Write the asymptotic variance as $\lambda^{\mathrm{T}} \Sigma \lambda$. The asymptotic normality is derived.
	
	For every component $k$, following the proof in \cite{Zheng.1996}, we can easily show that,
	\begin{equation}
	n h^{1 / 2} V_{2 nk} \stackrel{\mathbf{P}}{\longrightarrow} 0 \quad \text { and } \quad n h^{1 / 2} V_{3 nk} \stackrel{\mathbf{P}}{\longrightarrow} 0 .\label{A3}
	\end{equation}
	As $p$ is fixed, we then have that $n h^{1 / 2} V_{2 n}$ and $n h^{1 / 2} V_{3 n}$ converge to zero  in probability. The details are omitted here.
	
	Summarizing the results (\ref{A1}), (\ref{A2}) and (\ref{A3}), we conclude $n h^{1 / 2} V_{n} \stackrel{\mathbf{d}}{\longrightarrow} \mathrm{N}(0, \Sigma)$.
\end{proof}

\subsection{Proof of Lemma~2}

\begin{proof}
	Similarly as the  proof of Lemma~1, it is easy to decompose every component of $\hat{\Sigma}$ as
	\begin{equation*}
	\begin{aligned}
	\hat{\Sigma}_{k_1k_2}=& \frac{2}{n(n-1)} \sum_{i=1}^{n} \sum_{j=1 \atop j \neq i}^{n} \frac{1}{h} K^{2}\left(\frac{t_{i}-t_{j}}{h}\right) e_{ik_1}e_{ik_2}e_{jk_1}e_{jk_2} \\
	=&2\left[\frac{1}{n(n-1)} \sum_{i=1}^{n} \sum_{j=1 \atop j \neq i}^{n} \frac{1}{h} K^{2}\left(\frac{t_{i}-t_{j}}{h}\right) \varepsilon_{ik_1} \varepsilon_{ik_2}\varepsilon_{jk_1} \varepsilon_{jk_2}\right]+\mathrm{o}_{P}(1) \\
	\equiv & 2 S_{nk_1k_2}^1+\mathrm{o}_{P}(1).
	\end{aligned}
	\end{equation*}
	Here $S_{nk_1k_2}^1$ is a standard U-statistic with the kernel:
	\begin{equation*}
	H_{n}\left(z_{i}, z_{j}\right)=\frac{1}{h} K^{2}\left(\frac{t_{i}-t_{j}}{h}\right) \varepsilon_{ik_1} \varepsilon_{ik_2}\varepsilon_{jk_1} \varepsilon_{jk_2} .
	\end{equation*}
	
	As in the proof of Lemma~1, it is easy find that $\mathrm{E}\left[\left\|H_{n}\left(z_{i}, z_{j}\right)\right\|^{2}\right]=\mathrm{o}(n)$. Applying Lemma~3.1 of \cite{Zheng.1996}, we have
	\begin{equation*}
	\begin{aligned}
	S_{nk_1k_2}^1&=\overline{r}_{nk_1k_2}+\mathrm{o}_{P}(1)\\
	&=\mathrm{E}\left[\frac{1}{h} K^{2}\left(\frac{t_{i}-t_{j}}{h}\right) \varepsilon_{ik_1} \varepsilon_{ik_2}\varepsilon_{jk_1} \varepsilon_{jk_2}\right] +\mathrm{o}_{P}(1)\\
	&=\int K^{2}(u) \mathrm{d} u \int\left[\sigma_{k_1k_2}(t)\right]^{2} p^{2}(t) \mathrm{d} t+\mathrm{o}_{P}(1)\\
	&=\Sigma_{k_1k_2} / 2+o_{P}(1).
	\end{aligned}
	\end{equation*}
	
	Thus we conclude
	\begin{equation*}
	\hat{\Sigma}=2 S_{n}^1+\mathrm{o}_{P}(1)=\Sigma+\mathrm{o}_{P}(1).
	\end{equation*}
\end{proof}

\subsection{Proof of Theorem~2.2 (under the null hypothesis)}

\begin{proof}
	The result is an easy consequence of Lemma~1 and Lemma~2 by using Slutsky's theorem and continuous mapping theorem.
\end{proof}

\subsection{Proof of Lemma~3}

\begin{proof}
	Again similar to the proof of Lemma~1, $V_n$ can be decomposed as
	\begin{equation*}
	\begin{aligned}
	V_{n} &=\frac{1}{n(n-1)} \sum_{i=1}^{n} \sum_{j=1 \atop j \neq i}^{n} \frac{1}{h} K\left(\frac{t_{i}-t_{j}}{h}\right) e_{i} \odot e_{j} \\ &=\left\{\frac{1}{n(n-1)} \sum_{i=1}^{n} \sum_{j=1 \atop j \neq i}^{n} \frac{1}{h} K\left(\frac{t_{i}-t_{j}}{h}\right) \varepsilon_{i} \odot \varepsilon_{j}\right\}+\mathrm{o}_{P}(1) \\
	& \equiv S_{n}^2+\mathrm{o}_{P}(1) .
	\end{aligned}
	\end{equation*}
	Here $S_{n}^2$ is also a standard U-statistic with the kernel:
	\begin{equation*}
	H_{n}\left(z_{i}, z_{j}\right)=\frac{1}{h} K\left(\frac{t_{i}-t_{j}}{h}\right) \varepsilon_{i} \odot \varepsilon_{j}.
	\end{equation*}
	
	The conditions in Lemma~3.1 of \cite{Zheng.1996} can be easily verified using the same methods in the proof of Lemma~1. Since
	\begin{equation*}
	\mathrm{E} (\varepsilon_i | t_i)=X\left(t_{i}\right)-F\left(t_{i}, \theta_{0}\right)
	\end{equation*}
	under the global alternative. We then have
	\begin{equation*}
	\begin{aligned}
	\overline{r}_{n} =&\mathrm{E}\left\{\mathrm{E}\left[H_{n}\left(z_{i}, z_{j}\right) | t_{i}, t_{j}\right]\right\} \\
	=&\frac{1}{h} \mathrm{E}\left\{K\left(\frac{t_{i}-t_{j}}{h}\right)\left[X\left(t_{i}\right)-F\left(t_{i}, \theta_{0}\right)\right]\odot \left[X\left(t_{j}\right)-F\left(t_{j}, \theta_{0}\right)\right]\right\}\\
	=&\frac{1}{h} \int K\left(\frac{t_{i}-t_{j}}{h}\right)\left[X\left(t_{i}\right)-F\left(t_{i}, \theta_{0}\right)\right]\odot \left[X\left(t_{j}\right)-F\left(t_{j}, \theta_{0}\right)\right] p\left(t_{i}\right) p\left(t_{j}\right) \mathrm{d} t_{i} \mathrm{d} t_{j}\\
	=& \frac{1}{h} \int K(u)\left[X\left(t_{i}\right)-F\left(t_{i}, \theta_{0}\right)\right]\odot \left[X\left(t_{i}-h u\right)-F\left(t_{i}-h u, \theta_{0}\right)\right] \\ & \times p\left(t_{i}\right) p\left(t_{i}-h u\right) \mathrm{d} t_{i} h \mathrm{d} u\\
	=&\int\left[X(t)-F\left(t, \theta_{0}\right)\right]^{2} \odot p^{2}(t) \mathrm{d} t+\mathrm{o}(1)\\ =&\mathrm{E}\left\{\left[X\left(t_{i}\right)-F\left(t_{i}, \theta_{0}\right)\right]^{2} \odot p\left(t_{i}\right)\right\}+\mathrm{o}(1) .
	\end{aligned}
	\end{equation*}
	
	Thus
	\begin{equation*}
	V_{n} \stackrel{\mathbf{P}}{\longrightarrow} E\left\{\left[X\left(t_{i}\right)-F\left(t_{i}, \theta_{0}\right)\right]^{2} \odot p\left(t_{i}\right)\right\}.
	\end{equation*}
\end{proof}

\subsection{Proof of Lemma~4}

\begin{proof}
	By the similar proof of Lemma~3, it is easy to show that
	\begin{equation*}
	\begin{aligned}
	\hat{\Sigma}_{k_1k_2}=&2\left\{\frac{1}{n(n-1)} \sum_{i=1}^{n} \sum_{j=1 \atop j \neq i}^{n} \frac{1}{h} K^{2}\left(\frac{t_{i}-t_{j}}{h}\right) \varepsilon_{ik_1} \varepsilon_{ik_2}\varepsilon_{jk_1} \varepsilon_{jk_2}\right\}+\mathrm{o}_{P}(1) \\ & \equiv 2 S_{nk_1k_2}^3+\mathrm{o}_{P}(1) .
	\end{aligned}
	\end{equation*}
	Here $S_{nk_1k_2}^3$ is a standard U-statistic with the kernel
	\begin{equation*}
	H_{n}\left(z_{i}, z_{j}\right)=\frac{1}{h} K^{2}\left(\frac{t_{i}-t_{j}}{h}\right) \varepsilon_{ik_1} \varepsilon_{ik_2}\varepsilon_{jk_1} \varepsilon_{jk_2}.
	\end{equation*}
	
	We can easily find that $\mathrm{E}\left[\left\|H_{n}\left(z_{i}, z_{j}\right)\right\|^{2}\right]=\mathrm{o}(n)$ by direct computation. Applying Lemma~3.1 of \cite{Zheng.1996}, we have
	\begin{equation*}
	\begin{aligned}
	S_{nk_1k_2}^3&=\overline{r}_{nk_1k_2}+\mathrm{o}_{P}(1)\\
	&=\mathrm{E}\left[\frac{1}{h} K^{2}\left(\frac{t_{i}-t_{j}}{h}\right) \varepsilon_{ik_1} \varepsilon_{ik_2}\varepsilon_{jk_1} \varepsilon_{jk_2}\right] +\mathrm{o}_{P}(1)\\
	&=\int K^{2}(u) \mathrm{d} u \int\left[\sigma_{k_1k_2}(t)+(X_{k_1}(t)-F_{k_1}(t,\theta_0))(X_{k_2}(t)-F_{k_2}(t,\theta_0))\right]^{2} p^{2}(t) \mathrm{d} t+\mathrm{o}_{P}(1)\\
	&=\Sigma_{k_1k_2}^{\prime} / 2+\mathrm{o}_{P}(1).
	\end{aligned}
	\end{equation*}
	Thus
	\begin{equation*}
	\hat{\Sigma}=2 S_{nk_1k_2}^3+\mathrm{o}_{P}(1)=\Sigma^{\prime}+\mathrm{o}_{P}(1).
	\end{equation*}
\end{proof}

\subsection{Proof of Theorem~2.3 (under global alternatives)}

\begin{proof}
	The result is an easily derived consequence of Lemma~3 and Lemma~4.
\end{proof}

\subsection{Proof of Theorem~2.4 (under local alternatives)}

\begin{proof}
	
	Here we just focus on giving the limiting distribution of $V_n^F$ in the case that $p=1$. The arguments in the proof can be easily applied to handle multidimensional case and obtain the convergence result of $TM_{n}$.
	
	In the case that $p=1$, $V_{n}^F$ is similar to Zheng's statistic $V_n$ and it seems that it follows the results of Theorem 3 of \cite{Zheng.1996}. However, in that proof, the author gives (A.37) to show that the limit distribution of $V_n$ only depends on the limit distribution of $S_{7n}$, which is not enough since we need to show that $V_n-S_{7n}$ is $\mathrm{o}_{P}(n^{-1}h^{-1/2})$ instead of $\mathrm{o}_{P}(1)$. Actually, when the rate of $(\hat{\theta}-\theta)$ is slower than $1/\sqrt{n}$, the result of Theorem 3 is incorrect. We give the result of $V_n^F$ as follows.
	
	Similarly as the proof of Lemma~1, $V_n^F$ can be decomposed as
	
	\begin{equation*}
	\begin{aligned}
	V_{n}^F=& \frac{1}{n(n-1)} \sum_{i=1}^{n} \sum_{j=1 \atop j \neq i}^{n} \frac{1} K\left(\frac{t_{i}-t_{j}}{h}\right) e_{i} e_{j} \\
	=&\left[ \frac{1}{n(n-1)} \sum_{i=1}^{n} \sum_{j=1 \atop j \neq i}^{n} \frac{1} K\left(\frac{t_{i}-t_{j}}{h}\right) \varepsilon_{i} \varepsilon_{j}\right] \\
	&-2\left\{\frac{1}{n(n-1)} \sum_{i=1}^{n} \sum_{j=1 \atop j \neq i}^{n} \frac{1} K\left(\frac{t_{i}-t_{j}}{h}\right) \varepsilon_{i}\left[F\left(t_{j}, \hat{\theta}\right)-F\left(t_{j}, \theta_{0}\right)\right]\right\} \\
	&+\left\{\frac{1}{n(n-1)} \sum_{i=1}^{n} \sum_{j=1 \atop j \neq i}^{n} \frac{1} K\left(\frac{t_{i}-t_{j}}{h}\right) \right. \\
	&\left.\times\left[F\left(t_{i},\hat{\theta}\right)-F\left(t_{i}, \theta_{0}\right)\right]\left[F\left(t_{j}, \hat{\theta}\right)-F\left(t_{j}, \theta_{0}\right)\right]\right\} \\
	\equiv & S_{1 n}-2 S_{2 n}+S_{3 n}.
	\end{aligned}
	\end{equation*}
	
	For $S_{1n}$,
	\begin{equation*}
	\begin{aligned}
	S_{1n} =& \frac{1}{n(n-1)} \sum_{i=1}^{n} \sum_{j=1 \atop j \neq i}^{n} \frac{1}{h} K\left(\frac{t_{i}-t_{j}}{h}\right)\left[u_{i}+\delta_{n}  L\left(t_{i}\right)\right]\left[u_{j}+\delta_{n}  L\left(t_{j}\right)\right]\\
	=&\left[ \frac{1}{n(n-1)} \sum_{i=1}^{n} \sum_{j=1 \atop j \neq i}^{n} \frac{1}{h} K\left(\frac{t_{i}-t_{j}}{h}\right) u_{i} u_{j}\right]\\
	&+\delta_{n} \left[ \frac{1}{n(n-1)} \sum_{i=1}^{n} \sum_{j=1 \atop j \neq i}^{n} \frac{1}{h} K\left(\frac{t_{i}-t_{j}}{h}\right) u_{i} L\left(t_{j}\right)\right]
	\\&+\delta_{n}^{2}\left[\frac{1}{n(n-1)} \sum_{i=1}^{n} \sum_{j=1 \atop j \neq i}^{n} \frac{1}{h} K\left(\frac{t_{i}-t_{j}}{h}\right) L\left(t_{i}\right) L\left(t_{j}\right)\right]\\
	\equiv & Q_{1 n}+\delta_{n}  Q_{2 n}+\delta_{n}^{2}Q_{3 n}
	\end{aligned}
	\end{equation*}
	where $u_i=\varepsilon_i-\delta_{n} L\left(t_{i}\right)$.
	
	By a similar proof used for Theorem 3 of \cite{Zheng.1996}, we can easily show that
	$n h^{1 / 2} Q_{1n} \stackrel{\mathbf{d}}{\longrightarrow} \mathrm{N}(0, \Sigma)$, $\sqrt{n} Q_{2 n} \stackrel{\mathbf{d}}{\longrightarrow} \mathrm{N}\left(0, \mathrm{E}\left[\sigma^{2}\left(t_{i}\right) L^{2}\left(t_{i}\right) p^{2}\left(t_{i}\right)\right]\right)$, $Q_{3 n} \stackrel{\mathbf{P}}{\longrightarrow} \mathrm{E}\left[L^{2}\left(t_{i}\right) p\left(t_{i}\right)\right]$.
	
	For $S_{2n}$, we have
	\begin{equation*}
	\begin{aligned}
	S_{2 n}=& \frac{1}{n(n-1)} \sum_{i=1}^{n} \sum_{j=1 \atop j \neq i}^{n} \frac{1}{h} K\left(\frac{t_{i}-t_{j}}{h}\right) \varepsilon_{i}\left[F\left(t_{j}, \hat{\theta}\right)-F\left(t_{j}, \theta_{0}\right)\right] \\
	=&\frac{1}{n(n-1)} \sum_{i=1}^{n} \sum_{j=1 \atop j \neq i}^{n} \frac{1}{h} K\left(\frac{t_{i}-t_{j}}{h}\right) [u_{i}+\delta_{n} L\left(t_{i}\right)]\left[F\left(t_{j}, \hat{\theta}\right)-F\left(t_{j}, \theta_{0}\right)\right] \\
	=&\left\{\frac{1}{n(n-1)} \sum_{i=1}^{n} \sum_{j=1 \atop j \neq i}^{n} \frac{1}{h} K\left(\frac{t_{i}-t_{j}}{h}\right) u_{i}\left[F\left(t_{j}, \hat{\theta}\right)-F\left(t_{j}, \theta_{0}\right)\right]\right\}\\
	&+\left\{\frac{\delta_n}{n(n-1)} \sum_{i=1}^{n} \sum_{j=1 \atop j \neq i}^{n} \frac{1}{h} K\left(\frac{t_{i}-t_{j}}{h}\right) L\left(t_{i}\right)\left[F\left(t_{j}, \hat{\theta}\right)-F\left(t_{j}, \theta_{0}\right)\right]\right\}\\
	=&\left\{\frac{1}{n(n-1)} \sum_{i=1}^{n} \sum_{j=1 \atop j \neq i}^{n} \frac{1}{h} K\left(\frac{t_{i}-t_{j}}{h}\right) u_{i}\frac{\partial F\left(t_{j}, \theta_{0}\right)}{\partial \theta^{\mathrm{T}}}\left(\hat{\theta}-\theta_{0}\right)\right\}\\
	&+\left\{\left(\hat{\theta}-\theta_{0}\right)^{T}\frac{1}{n(n-1)} \sum_{i=1}^{n} \sum_{j=1 \atop j \neq i}^{n} \frac{1}{h} K\left(\frac{t_{i}-t_{j}}{h}\right) u_{i}\frac{\partial^{2} F\left(t_{j}, \theta\right)}{\partial \theta \partial \theta^{\mathrm{T}}}\left(\hat{\theta}-\theta_{0}\right)\right\}\\
	&+\left\{\frac{\delta_n}{n(n-1)} \sum_{i=1}^{n} \sum_{j=1 \atop j \neq i}^{n} \frac{1}{h} K\left(\frac{t_{i}-t_{j}}{h}\right) L\left(t_{i}\right)\frac{\partial F\left(t_{j}, \theta_{0}\right)}{\partial \theta^{\mathrm{T}}}\left(\hat{\theta}-\theta_{0}\right)\right\}\\
	&+\left\{\left(\hat{\theta}-\theta_{0}\right)^{T}\frac{\delta_n}{n(n-1)} \sum_{i=1}^{n} \sum_{j=1 \atop j \neq i}^{n} \frac{1}{h} K\left(\frac{t_{i}-t_{j}}{h}\right) L\left(t_{i}\right)\frac{\partial^{2} F\left(t_{j}, \theta\right)}{\partial \theta \partial \theta^{T}}\left(\hat{\theta}-\theta_{0}\right) \right\}\\
	&+\mathrm{O}_{P}(\|\theta_n-\theta_0\|^3)+\mathrm{O}_{P}(\delta_n \|\theta_n-\theta_0\|^3)\\
	\equiv & Q_{4 n}+Q_{5 n}+Q_{6 n}+Q_{7n}+\mathrm{O}_{P}(\|\theta_n-\theta_0\|^3)+\mathrm{O}_{P}(\delta_n \|\theta_n-\theta_0\|^3).
	\end{aligned}
	\end{equation*}
	
	According to Theorem~2.1, under the local alternatives,
	\begin{equation}
	\theta_n-\theta_0=H_{\dot{F}}^{-1} \frac{1}{n} \sum_{i=1}^{n} \varepsilon_{i}\frac{\partial F\left(t_{i}, \theta_{0}\right)}{\partial \theta}+ \delta_{n} H_{\dot{F}}^{-1} \mathrm{E}\left[L\left(t_{i}\right)\frac{\partial F\left(t_{i}, \theta_{0}\right)}{\partial \theta} \right]=\mathrm{O}_{P}(\delta_n)\label{A17}
	\end{equation}
	where
	$$H_{\dot{F}}=\mathrm{E}\left[\frac{\partial F\left(t_{i}, \theta_{0}\right)}{\partial \theta} \frac{\partial F\left(t_{i}, \theta_{0}\right)}{\partial \theta^{\prime}} \right].$$
	
	Thus, by a similar proof of Lemma~3.3d of \cite{Zheng.1996}, we can show $Q_{4n}=\mathrm{O}_{P}(n^{-1/2}\delta_n)$, $Q_{5n}=\mathrm{O}_{P}(n^{-1/2}\delta_n^2)$,  $Q_{6n}=\mathrm{O}_{P}(\delta_n^2)$, and $Q_{7n}=\mathrm{O}_{P}(\delta_n^3)$. It is easy to see that the leading term in the4 above decomposition is $Q_{6n}$. Here,
	\begin{equation*}
	\begin{aligned}
	Q_{6 n}=&\left\{\frac{\delta_n}{n(n-1)} \sum_{i=1}^{n} \sum_{j=1 \atop j \neq i}^{n} \frac{1}{h} K\left(\frac{t_{i}-t_{j}}{h}\right) L\left(t_{i}\right)\frac{\partial F\left(t_{j}, \theta_{0}\right)}{\partial \theta^{\prime}}\left(\hat{\theta}-\theta_{0}\right)\right\}\\
	=&\delta_n^2\mathrm{E}\left[L\left(t_{i}\right)\frac{\partial F\left(t_{i}, \theta_{0}\right)}{\partial \theta^{\prime}} p\left(t_{i}\right)\right]H^{-1}_{\dot{f}}\mathrm{E}\left[L\left(t_{i}\right)\frac{\partial F\left(t_{i}, \theta_{0}\right)}{\partial \theta} \right].
	\end{aligned}
	\end{equation*}
	On the other hand,
	\begin{equation*}
	\begin{aligned}
	S_{3 n}=&\frac{1}{n(n-1)} \sum_{i=1}^{n} \sum_{j=1 \atop j \neq i}^{n} \frac{1}{h} K\left(\frac{t_{i}-t_{j}}{h}\right)  \\
	&\times\left[F\left(t_{i},\hat{\theta}\right)-F\left(t_{i}, \theta_{0}\right)\right]\left[F\left(t_{j}, \hat{\theta}\right)-F\left(t_{j}, \theta_{0}\right)\right]\\
	=& \left(\hat{\theta}-\theta_{0}\right)^{T}\frac{1}{n(n-1)} \sum_{i=1}^{n} \sum_{j=1 \atop j \neq i}^{n} \frac{1}{h} K\left(\frac{t_{i}-t_{j}}{h}\right) \\
	&\times\frac{\partial F\left(t_{i}, \theta_{0}\right)}{\partial \theta}\frac{\partial F\left(t_{j}, \theta_{0}\right)}{\partial \theta^{\prime}}\left(\hat{\theta}-\theta_{0}\right)[1+\mathrm{o}_P(1)]\\
	=&\delta_n^2 \mathrm{E}\left[L\left(t_{i}\right)\frac{\partial F\left(t_{i}, \theta_{0}\right)}{\partial \theta^{\prime}}\right]  H^{-1}_{\dot{F}}  \\
	&\times \mathrm{E} \left[ \frac{\partial F\left(t_{i}, \theta_{0}\right)}{\partial \theta}\frac{\partial F\left(t_{i}, \theta_{0}\right)}{\partial \theta^{T}}p(t_i) \right] H^{-1}_{\dot{F}}\mathrm{E}\left[L\left(t_{i}\right)\frac{\partial F\left(t_{i}, \theta_{0}\right)}{\partial \theta}\right][1+\mathrm{o}_P(1)].
	\end{aligned}
	\end{equation*}
	Thus, when $\delta_n=n^{-1 / 2} h^{-m / 4}$, we have
	$n h^{m / 2} V_{n} \stackrel{\mathbf{P}}{\longrightarrow} \mathrm{N}\left(V, \Sigma\right)$
	where
	\begin{equation*}
	\begin{aligned}
	V=&\mathrm{E}\left[L^{2}\left(t_{i}\right) p\left(t_{i}\right)\right]-2\mathrm{E}\left[L\left(t_{i}\right)\frac{\partial F\left(t_{i}, \theta_{0}\right)}{\partial \theta^{\prime}} p\left(t_{i}\right)\right]H^{-1}_{\dot{F}}\mathrm{E}\left[L\left(t_{i}\right)\frac{\partial F\left(t_{i}, \theta_{0}\right)}{\partial \theta} \right]\\
	&+\mathrm{E}\left[L\left(t_{i}\right)\frac{\partial F\left(t_{i}, \theta_{0}\right)}{\partial \theta^{\prime}}\right]  H^{-1}_{\dot{F}} \mathrm{E} \left[ \frac{\partial F\left(t_{i}, \theta_{0}\right)}{\partial \theta}\frac{\partial F\left(t_{i}, \theta_{0}\right)}{\partial \theta^{\prime}}p(t_i) \right] H^{-1}_{\dot{f}}\mathrm{E}\left[L\left(t_{i}\right)\frac{\partial F\left(t_{i}, \theta_{0}\right)}{\partial \theta}\right]\\
	=&\mathrm{E}\left\{\left\{ L\left(t_{i}\right)-\frac{\partial F\left(t_{i}, \theta_{0}\right)}{\partial \theta^{\prime}} H^{-1}_{\dot{f}}\mathrm{E}\left[L\left(t_{i}\right)\frac{\partial F\left(t_{i}, \theta_{0}\right)}{\partial \theta}\right]  \right\}^2 p(t) \right\}.
	\end{aligned}
	\end{equation*}
	
	The extension to the multivariate case ($p >1$) is straightforward. Then the convergence result of $TM_{n}$ is easy to derive by using the convergence result of $V_n$.
\end{proof}
\section{The results in Section~3}
\subsection{Proof of Theorem~3.1}

\begin{proof}
	Under $H_{1n}^F$, since the local alternative model is $X(t) = F\left(t, \theta_{0}\right)+\delta_{n} L_1\left(t\right)$, we can deal with the derivatives on both side:
	\begin{equation*}	
	\begin{aligned}
	X^{\prime}(t) =& F^{\prime}\left(t, \theta_{0}\right)+\delta_{n} L^{\prime}\left(t\right)\\
	=&f\left(t, X(t)-\delta_{n} L\left(t\right);\theta_{0}\right)+\delta_{n} L^{\prime}\left(t\right)\\
	=&f\left(t, X(t);\theta_{0}\right)+\delta_{n} L^{\prime}\left(t\right)-\delta_{n} \frac{\partial f\left(t, X(t);\theta_{0}\right)}{\partial X^{\mathrm{T}}}L\left(t\right) + \mathrm{o}(\delta_{n})v_2(t).\\
	\end{aligned}
	\end{equation*}
	
	Thus we have $v_1(t)=L^{\prime}\left(t\right)-\frac{\partial f\left(t, X(t);\theta_{0}\right)}{\partial X^{\mathrm{T}}}L\left(t\right)$ and the former part of this theorem is proven. The latter part can be proven by contradiction.
\end{proof}

\subsection{Proof of Theorem~3.2}

\begin{proof}
	We can show $\hat{\theta}$ is a consistent estimator of $\theta^*$ by mimicking the proof of Theorem 1 of \cite{Liang.2008}. Since the proof is similar, here we omit the details. Next we focus on giving the asymptotically linear representation of $\hat{\theta}-\theta^*$ and its root-$n$ consistency.
	
	The two-stage collocation estimator $\hat{\theta}$ is defined as
	\begin{equation*}
	S_n(\theta)=\frac{1}{m} \sum_{j=1}^{m} \left[\hat{X}_k^{\prime}(t_j^{*})-f_k\left(t_j^*, \hat{X}(t_j^{*}),\theta \right)\right]^2\omega(t_j^*),
	\end{equation*}
	\begin{equation*}
	\hat{\theta}_n=\arg \min_{\theta} S_n(\theta).
	\end{equation*}
	Using Taylor expansion, we obtain
	\begin{equation*}
	\dot{S}_n(\hat{\theta}_n)-\dot{S}_n(\theta^*)=\ddot{S}_n(\tilde{\theta}_n)(\hat{\theta}_n-\theta^*)
	\end{equation*}
	where $\tilde{\theta}_n$ is a mid-value between $\hat{\theta}$ and $\theta^*$.
	Thus,
	\begin{equation}
	\begin{aligned} \hat{\theta}_n-\theta^* &=\ddot{S}_n(\tilde{\theta}_n)^{-1}[\dot{S}_n(\hat{\theta}_n)-\dot{S}_n(\theta^*)] \\
	&=-\ddot{S}_n(\tilde{\theta}_n)^{-1}\dot{S}_n(\theta^*). \label{main}
	\end{aligned}
	\end{equation}
	
	For $\dot{S}_n(\theta^*)$, we have
	\begin{equation}
	\begin{aligned}
	\dot{S}_n(\theta^*)=&-\frac{2}{m} \sum_{j=1}^{m} \left[ \hat{X}_k^{\prime}(t_i^*)-f_k\left(t_j^*,\hat{X}(t_j^*),\theta^* \right)\right]\omega(t_j^*) \frac{\partial f_k\left(t_j^*,\hat{X}(t_j^*),\theta^* \right)}{\partial \theta}\\
	=&-\frac{2}{m} \sum_{j=1}^{m} \left[ X^{\prime}_j(t_i^*)-f_k\left(t_j^*,\hat{X}(t_j^*),\theta^* \right)\right] \omega(t_j^*)\frac{\partial f_k\left(t_j^*,\hat{X}(t_j^*),\theta^* \right)}{\partial \theta}\\
	&-\frac{2}{m} \sum_{j=1}^{m} \left[ \hat{X}^{\prime}_k(t_j^*)-X_k^{\prime}(t_j^*)\right]\omega(t_j^*) \frac{\partial f_k\left(t_j^*,\hat{X}(t_j^*),\theta^* \right)}{\partial \theta}.\label{A4}\\
	\end{aligned}
	\end{equation}
	
	We consider the  asymptotic approximations of the estimator under the null, global and local alternative hypothesis.\\
	{\bf 1.  Under $H_{0k}$. } We note that under the null hypothesis,
	\begin{equation}
	\begin{aligned}
	&E_{p^*}\left\{ \left[X_k^{\prime}(t)-f_k\left(t, X(t), \theta\right)\right]^{2}w(t)\right\}\\
	=&E_{p^*}\left\{ \left[f_k\left(t, X(t), \theta_0\right)-f_k\left(t, X(t), \theta)\right)\right]^{2}w(t)\right\}\\
	\geq& E_{p^*}\left\{ \left[f_k\left(t, X(t), \theta_0\right)-f_k\left(t, X(t), \theta_0\right)\right]^{2}w(t)\right\}=0.\\
	\end{aligned}
	\end{equation}
	Therefore, $\theta^*=\theta_{0}$. Based on (\ref{A4}), we have
	\begin{equation}
	\begin{aligned}
	\dot{S}_n(\theta^*)=&-\frac{2}{m} \sum_{j=1}^{m} \left[ f_k\left(t_j^*,X(t_j^*),\theta_0 \right)-f_k\left(t_j^*,\hat{X}(t_j^*),\theta_0 \right)\right] \omega(t_j^*)\frac{\partial f_k\left(t_j^*,\hat{X}\left(t_{j}^*\right), \theta_0 \right)}{\partial \theta}\\
	&-\frac{2}{m} \sum_{j=1}^{m} \left[\hat{X}_k^{\prime}(t_j^*)-X_k^{\prime}(t_j^*)\right]\omega(t_j^*) \frac{\partial f_k\left(t_j^*,\hat{X}\left(t_{j}^*\right), \theta_0 \right)}{\partial \theta}\\
	\equiv&-2V_{01}-2V_{02}. \label{S1st}
	\end{aligned}
	\end{equation}
	Consider $V_{01}$. It can be decomposed as
	\begin{equation*}
	\begin{aligned}
	V_{01}=&\frac{1}{m} \sum_{j=1}^{m} \left[f_k\left(t_j^*, X(t_j^*),\theta_0 \right)-f_k\left(t_j^*,\hat{X}(t_j^*),\theta_0 \right)\right]\omega(t_j^*) \frac{\partial f_k\left(t_j^*,\hat{X}\left(t_{j}^*\right), \theta_0 \right)}{\partial \theta}\\
	=&\frac{1}{m} \sum_{j=1}^{m} \left[ f_k\left(t_j^*,X(t_j^*),\theta_0 \right)-f_k\left(t_j^*, \hat{X}(t_j^*),\theta_0 \right)\right] \omega(t_j^*)\\
	&\times \left[ \frac{\partial f_k\left(t_j^*,\hat{X}\left(t_{j}^*\right), \theta_0 \right)}{\partial \theta}-\frac{\partial f_k\left(t_j^*,X\left(t_{j}^*\right), \theta_0 \right)}{\partial \theta}\right]\\
	&+\frac{1}{m} \sum_{j=1}^{m} \left[f_k\left(t_j^*,X(t_j^*),\theta_0 \right)-f_k\left(t_j^*,\hat{X}(t_j^*),\theta_0 \right)\right] \omega(t_j^*)\frac{\partial f_k\left(t_j^*,X\left(t_{j}^*\right), \theta_0 \right)}{\partial \theta}\\
	\equiv &V_{011}+V_{012}=V_{012} [1+\mathrm{o}_P(1)],
	\end{aligned}
	\end{equation*}
	as $V_{011}=o_P(V_{012})$ by noting that,  $\hat{X}(t_j^*)$ is a consistent estimator of $X(t_j^*)$,
	\begin{equation}
	\frac{\partial f_k\left(t_j^*,\hat{X}\left(t_{j}^*\right), \theta_0 \right)}{\partial \theta}-\frac{\partial f_k\left(t_j^*,X\left(t_{j}^*\right), \theta_0 \right)}{\partial \theta}=\mathrm{o}_P(1).\label{A5}
	\end{equation}
	Then we consider $V_{012}$ decomposed as
	\begin{equation}
	\begin{aligned}
	V_{012}=&\frac{1}{m} \sum_{j=1}^{m} \left[ f_k\left( t_j^*,X(t_j^*),\theta_0 \right)-f_
	k \left( t_j^*,\hat{X}(t_j^*),\theta_0 \right) \right] \omega(t_j^*)\frac{\partial f_k\left(t_j^*, X\left(t_{j}^*\right), \theta_0 \right)}{\partial \theta}\\
	=&\frac{-1}{m} \sum_{j=1}^{m} \omega(t_j^*)\frac{\partial f_k\left(t_j^*,X\left(t_{j}^*\right), \theta_0 \right)}{\partial \theta} \frac{\partial f_k\left( t_j^*,X\left(t_{j}^*\right), \theta_0 \right)}{\partial X^\mathrm{T}}\\
	&\times \left[\hat{X}(t_j^*)-X(t_j^*)\right][1+\mathrm{o}_P(1)]. \label{v012}\\
	\end{aligned}
	\end{equation}
	For $V_{02}$, the decomposition is as
	\begin{equation}
	\begin{aligned}
	V_{02}=&\frac{1}{m} \sum_{j=1}^{m} \left[\hat{X}^{\prime}_k(t_j^*)-X_k^{\prime}(t_j^*)\right]\omega(t_j^*) \frac{\partial f_k\left( t_j^*,\hat{X}\left(t_{j}^*\right), \theta_0 \right)}{\partial \theta}\\
	=&\frac{1}{m} \sum_{i=1}^{m} \left[\hat{X}_k^{\prime}(t_j^*)-X_k^{\prime}(t_j^*)\right]\omega(t_j^*) \\
	&\times \left[ \frac{\partial f_k\left(t_j^*,\hat{X}\left(t_{j}^*\right), \theta_0 \right)}{\partial \theta}-\frac{\partial f_k\left(t_j^*,X\left(t_{j}^*\right), \theta_0 \right)}{\partial \theta}\right]\\
	&+\frac{1}{m} \sum_{j=1}^{m} \left[\hat{X}_k^{\prime}(t_j^*)-X_k^{\prime}(t_j^*)\right] \omega(t_j^*)\frac{\partial f_k\left( t_j^*, X\left(t_{j}^*\right), \theta_0 \right)}{\partial \theta}\\
	\equiv &V_{021}+V_{022}=V_{022} [1+\mathrm{o}_P(1)], \label{v022}
	\end{aligned}
	\end{equation}
	as $V_{021}$ can be proven to be $o_P(V_{022})$ by using (\ref{A5}).
	Together with (\ref{v012}) and (\ref{v022}), $\dot{S}_n(\theta^*)$ in (\ref{S1st}) has the linear approximation $(V_{012}+V_{022})(1+o_p(1)).$

	Next we consider the second order derivative $\ddot{S}_n(\tilde{\theta}_n)$ of $S_n$ with respect to $\theta$. We have
	\begin{equation*}
	\begin{aligned}
	\ddot{S}_n(\tilde{\theta}_n)=&-\frac{2}{m} \sum_{j=1}^{m}\left[ \hat{X}^{\prime}_k\left(t_{j}^*\right)-f_k\left(t_{j}^*,\hat{X}\left(t_{j}^*\right), \tilde{\theta}_n\right)\right]\omega(t_j^*) \frac{\partial^{2} f_k\left(t_{j}^*,\hat{X}\left(t_{j}^*\right), \tilde{\theta}_n\right)}{\partial \theta \partial \theta^{\mathrm{T}}}\\
	&+\frac{2}{m} \sum_{j=1}^{m} \omega(t_j^*)\frac{\partial f_k\left(t_{j}^*,\hat{X}\left(t_{j}^*\right), \tilde{\theta}_n\right)}{\partial \theta}\frac{\partial f_k\left(t_{j}^*,\hat{X}\left(t_{j}^*\right), \tilde{\theta}_n\right)}{\partial \theta^{\mathrm{T}}}\\
	=&2E_{p^*} \left[ \omega(t)\frac{\partial f_k\left( t,X\left(t\right), \theta_0 \right)}{\partial \theta}\frac{\partial f_k\left(t,X\left(t\right), \theta_0 \right)}{\partial \theta^{\mathrm{T}}}\right]+\mathrm{o}_P(1)\\
	\equiv &2H_{\dot{f}}+\mathrm{o}_P(1).
	\end{aligned}
	\end{equation*}
	Altogether, the linear approximation of $\hat{\theta}_{n}-\theta_{0}$ is as
	\begin{equation}
	\begin{aligned}
	&\hat{\theta}_{n}-\theta_{0}\\
	=&H_{\dot{f}}^{-1}\frac{1}{m} \sum_{j=1}^{m}\left\{ \left[\hat{X}_k^{\prime}(t_j^*)-X_k^{\prime}(t_j^*)\right] \omega(t_j^*)\frac{\partial f_k\left(t_j^*,X\left(t_{j}^*\right), \theta_0 \right)}{\partial \theta} \right.\\
	&\left.- \omega(t_j^*)\frac{\partial f_k\left(t_j^*,X\left(t_{j}^*\right), \theta_0 \right)}{\partial \theta} \frac{\partial f_k\left(t_j^*,X\left(t_{j}^*\right), \theta_0 \right)}{\partial X^\mathrm{T}} \left[\hat{X}(t_j^*)-X(t_j^*)\right]\right\}\left[ 1+o_P(1) \right]\\
	=&H_{\dot{f}}^{-1}\frac{1}{m} \sum_{j=1}^{m}\left[ \Delta_k(t_j^*) \omega(t_j^*)\frac{\partial f_k\left(t_j^*,X\left(t_{j}^*\right), \theta_0 \right)}{\partial \theta}\right. \\
	&\left.- \omega(t_j^*)\frac{\partial f_k\left(t_j^*,X\left(t_{j}^*\right), \theta_0 \right)}{\partial \theta} \frac{\partial f_k\left(t_j^*,X\left(t_{j}^*\right), \theta_0 \right)}{\partial X^\mathrm{T}} \Lambda\left(t_{j}^{*}\right)\right]\left[ 1+o_P(1) \right]\label{A8}
	\end{aligned}
	\end{equation}
	where
	\begin{equation}
	\begin{aligned}
	\Lambda\left(t_{j}^{*}\right)=&\hat{X}\left(t_{j}^{*}\right)-X\left(t_{j}^{*}\right)\\
	=&\left\{\frac{1}{n p\left(t_{j}^{*}\right)} \sum_{s=1}^{n} \frac{1}{h} K\left(\frac{t_{s}-t_{j}^{*}}{h}\right)\left[X\left(t_{s}\right)-X\left(t_{j}^{*}\right)-X^{\prime}\left(t_{j}^{*}\right)\left(t_{s}-t_{j}^{*}\right)\right]\right.\\
	&\left.+\frac{1}{n p\left(t_{j}^{*}\right)} \sum_{s=1}^{n} \frac{1}{h} K\left(\frac{t_{s}-t_{j}^{*}}{h}\right) \epsilon\left(t_{s}\right)\right\}\left[1+o_{P}(1)\right], \label{A6}
	\end{aligned}
	\end{equation}
	\begin{equation}
	\begin{aligned}
	\Delta_k(t_j^*)=&\left[\hat{X}_k^{\prime}(t_j^*)-X_k^{\prime}(t_j^*)\right]\\
	=&\left\{\frac{1}{n h^{2} \mu_{2}(K) p\left(t_{j}^{*}\right)} \sum_{s=1}^{n} \frac{1}{h} K\left(\frac{t_{s}-t_{j}^{*}}{h}\right)\left(t_{s}-t_{j}^{*}\right)\right.\\
	&\times\left[X_k\left(t_{s}\right)-X_k\left(t_{j}^{*}\right)-X_k^{\prime}\left(t_{j}^{*}\right)\left(t_{s}-t_{i}^{*}\right)-X_k^{(2)}\left(t_{j}^{*}\right) \frac{\left(t_{s}-t_{j}^{*}\right)^{2}}{2}\right]\\
	&\left.+\frac{1}{n h^{2} \mu_{2}(K) p\left(t_{j}^{*}\right)} \sum_{s=1}^{n} \frac{1}{h} K\left(\frac{t_{s}-t_{j}^{*}}{h}\right)\left(t_{s}-t_{j}^{*}\right) \epsilon_k\left(t_{s}\right)\right\}\left[1+o_{P}(1)\right]. \label{A7}
	\end{aligned}
	\end{equation}
	We now prove the root-$n$ consistency of $\hat \theta_n-\theta_0$. Denote
	\begin{eqnarray*}	
		\frac{\partial f_k\left(t_{j}^*\right)}{\partial \theta}&=&\frac{\partial f_k\left(t_j^*,X\left(t_{j}^*\right),\theta_0 \right)}{\partial \theta},\\
		M_1&=&\frac{1}{m} \sum_{j=1}^{m}\Delta_k(t_j^*) \omega(t_j^*)\frac{\partial f_k\left(t_{j}^*\right)}{\partial \theta},\\
		M_2&=&\frac{1}{m} \sum_{j=1}^{m}\omega(t_j^*)\frac{\partial f_k\left(t_{j}^*\right)}{\partial \theta} \frac{\partial f_k\left(t_{j}^*\right)}{\partial X^\mathrm{T}} \Lambda\left(t_{j}^{*}\right).
	\end{eqnarray*}

	Based on the asymptotic form of $\Delta_k(t_j^*)$ in (\ref{A7}), we can compute the conditional variance of $M_1$ as follows:
	\begin{equation}
	\begin{aligned}
	&Var\left(M_{1} | \mathfrak{D}\right) \\
	=& \frac{1}{m^{2}} \sum_{j=1}^{m} \frac{\partial f_k\left(t_{j}^{*}\right)}{\partial \theta} \frac{\omega^2\left(t_{j}^{*}\right) \sigma^{2}_k(t_{j}^{*})\nu_{2}(K) }{n h^{3} \mu_{2}^{2}(K) p\left(t_{j}^{*}\right)} \frac{\partial f_k\left(t_{j}^{*}\right)}{\partial \theta^{\mathrm{T}}} \\
	&+\frac{1}{m^{2} n^{2} h^{4} \mu^{2}\left(K^{2}\right)} \sum_{l \neq i}^{m} \sum_{s=1}^{n} \frac{\sigma^{2}_k(t_s)}{h^{2} p\left(t_{i}^{*}\right) p\left(t_{l}^{*}\right)} \omega\left(t_{i}^{*}\right) \omega\left(t_{l}^{*}\right) \\
	&\times K\left(\frac{t_{s}-t_{i}^{*}}{h}\right) K\left(\frac{t_{s}-t_{l}^{*}}{h}\right) \left(t_{s}-t_{i}^{*}\right)\left(t_{s}-t_{l}^{*}\right) \frac{\partial f_k\left(t_{i}^{*}\right)}{\partial \theta} \frac{\partial f_k\left(t_{l}^{*}\right)}{\partial \theta^{\mathrm{T}}}\\
	=&\frac{1}{n h^{4} \mu_{2}^{2}(K)} \int_{t} \int_{z_{1}} \int_{z_{2}} \frac{h^{2}\sigma^{2}_k(t) K\left(z_{1}\right) K\left(z_{2}\right) z_{1} z_{2}}{p\left(t+z_{1} h\right) p\left(t+z_{2} h\right)} \frac{\partial f_k\left(t+z_{1} h\right)}{\partial \theta} \frac{\partial f_k\left(t+z_{2} h\right)}{\partial \theta^\mathrm{T}}\\
	&\times p^*\left(t+z_{1} h\right) p^*\left(t+z_{2} h\right) w\left(t+z_{1} h\right) w\left(t+z_{2} h\right) d z_{1} d z_{2} p(t) d t\\
	&+\frac{ \nu_{2}(K)}{n m h^{3} \mu_{2}^{2}(K)} E_{p^*}\left[\frac{1}{p(t)}\left(\omega(t)\sigma_k(t)\frac{\partial f_k(t)}{\partial \theta} \right)^{\otimes 2}\right]+o_{P}\left[\left(n m h^{3}\right)^{-1}\right]\\
	=&\frac{1}{n} E_{p^*}\left\{\frac{\sigma^{2}_k(t)}{p(t) p^*(t)}\left[\frac{\partial}{\partial t}\left(\omega(t) p^*(t)\frac{\partial f_k(t)}{\partial \theta} \right)\right]^{\otimes 2}\right\}\\
	&+\frac{ \nu_{2}(K)}{n m h^{3} \mu_{2}^{2}(K)} E_{p^*}\left[\frac{1}{p(t)}\left(\omega(t)\sigma_k(t)\frac{\partial f_k(t)}{\partial \theta} \right)^{\otimes 2}\right]+o_{P}\left[n^{-1}+\left(n m h^{3}\right)^{-1}\right]\\
	=&\frac{1}{n} V_{22}+\frac{1}{n m h^{3}} V_{22}^{*}+o_{P}\left[n^{-1}+\left(n m h^{3}\right)^{-1}\right].\label{A9}
	\end{aligned}
	\end{equation}
	Following the similar steps, we derive the conditional variance of $M_2$ and the conditional covariance of $(M_1, M_2)$, based on (\ref{A6}) and (\ref{A7}),
	\begin{equation}
	\begin{aligned}
	Var\left(M_{2} | \mathfrak{D}\right)=&\frac{1}{n} E_{p^*}\left[\frac{p^*(t)}{p(t)}\left(\omega(t)\frac{\partial f_k(t)}{\partial \theta} \frac{\partial f_k(t)}{\partial X^{\mathrm{T}}} \Sigma^{\frac{1}{2}}(t)\right)^{\otimes 2}\right]\\
	&+O_{P}\left[(n m h)^{-1}\right]+o_{P}\left[n^{-1}+(n m h)^{-1}\right]\\
	=&\frac{1}{n}V_{11}+O_{P}\left[(n m h)^{-1}\right]+o_{P}\left[n^{-1}+(n m h)^{-1}\right],\label{A10}
	\end{aligned}
	\end{equation}
	\begin{equation}
	\begin{aligned}
	Cov\left(M_{1}, M_{2} | \mathfrak{D}\right)=&\frac{1}{\mu_{2}(K)} E_{p^*}\left[\frac{\omega(t)}{p(t)} \frac{\partial f_k(t)}{\partial \theta} \frac{\partial f_k(t)}{\partial X^{\mathrm{T}}} \Sigma_k \frac{\partial}{\partial t}\left(\omega(t) p^*(t)\frac{\partial f_k(t)}{\partial \theta} \right)\right]\\
	&+O_{P}\left[(n m h)^{-1}\right]+o_{P}\left[n^{-1}+(n m h)^{-1}\right]\\
	=&\frac{1}{n}V_{12}+O_{P}\left[(n m h)^{-1}\right]+o_{P}\left[n^{-1}+(n m h)^{-1}\right].\label{A11}
	\end{aligned}
	\end{equation}
	Combining (\ref{A8}), (\ref{A9}), (\ref{A10}) and (\ref{A11}), we conclude
	\begin{equation*}
	Var \left[\sqrt{n}\left(\hat{\theta}_n-\theta_{0}\right)\right] \longrightarrow  H_{\dot{f}}^{-1}\left(V_{11}+V_{22}-2 V_{12}\right) H_{\dot{f}}^{-1}.
	\end{equation*}
	This implies the root-$n$ consistency.
	
	{\bf 2. Under global alternative hypothesis.}  Since $X_k^{\prime}(t)=f_k\left(t, X\left(t\right), \theta_0 \right)+q(t)$, the minimizer $\theta^*=\theta_{1}$ is a value which is  possibly to be different from $\theta_0$. Here, base on (\ref{A4}), we have
	\begin{equation*}
	\begin{aligned}
	\dot{S}_n(\theta^*)=&-\frac{2}{m} \sum_{j=1}^{m} \left[ f_k\left(t_j^*, X(t_j^*),\theta_0 \right)-f_k\left(t_j^*,\hat{X}\left(t_{j}^*\right), \theta_1 \right) \right] \omega(t_j^*)\frac{\partial f_k\left(t_j^*,\hat{X}\left(t_{j}^*\right), \theta_1 \right)}{\partial \theta}\\
	&-\frac{2}{m} \sum_{j=1}^{m} \left[\hat{X}_k^{\prime}(t_j^*)-X_k^{\prime}(t_j^*)\right]\omega(t_j^*) \frac{\partial f_k\left(t_j^*,\hat{X}\left(t_{j}^*\right), \theta_1 \right)}{\partial \theta}\\
	&-\frac{2}{m} \sum_{j=1}^{m} q(t_j^*) \omega(t_j^*)\frac{\partial f_k\left(t_j^*,\hat{X}\left(t_{j}^*\right), \theta_1 \right)}{\partial \theta}\\
	\equiv &-2V_{21}-2V_{22}-2V_{23}.
	\end{aligned}
	\end{equation*}
	As the proof is very similar to the above, we will give the detail briefly somehow.
	For $V_{21}$, we have
	\begin{equation*}
	\begin{aligned}
	V_{21}=&\frac{1}{m} \sum_{j=1}^{m} \left[ f_k\left(t_j^*, X(t_j^*),\theta_0 \right)-f_k\left(t_j^*,\hat{X}\left(t_{j}^*\right), \theta_1 \right) \right] \omega(t_j^*)\\
	&\times \frac{\partial f\left( t_{j}^*,X\left(t_{j}^*\right), \theta_1 \right)}{\partial \theta}\left[ 1+o_P(1) \right]\\
	=&\frac{1}{m} \left\{ \sum_{j=1}^{m} \left[ f_k\left(t_j^*, X(t_j^*),\theta_1 \right)-f_k\left(t_j^*,\hat{X}\left(t_{j}^*\right), \theta_1 \right) \right] \omega(t_j^*)\frac{\partial f_k\left( t_{j}^*,X\left(t_{j}^*\right), \theta_1 \right)}{\partial \theta}\right.\\
	&+\left[ f_k\left( t_j^*, X(t_j^*),\theta_0 \right)-f_k\left(t_j^*, X(t_j^*),\theta_1 \right)\right] \omega(t_j^*)\\
	&\left. \times \frac{\partial f_k\left( t_j^*, X\left(t_{j}^*\right), \theta_1 \right)}{\partial \theta} \right\}\left[ 1+o_P(1) \right]\\
	=&\frac{1}{m} \sum_{j=1}^{m} \left\{ -\omega(t_j^*)\frac{\partial f_k\left( t_j^*,X\left(t_{j}^*\right), \theta_1 \right)}{\partial \theta} \frac{\partial f_k\left( t_j^*, X\left(t_j^*\right), \theta_1 \right)}{\partial X^\mathrm{T}}\left[\hat{X}(t_j^*)-X(t_j^*)\right] \right.\\
	&\left.+\left[f_k\left(t_j^*,X(t_j^*),\theta_0 \right)-f_k\left( t_j^*,X(t_j^*),\theta_1 \right)\right] \omega(t_j^*)\frac{\partial f_k\left(t_j^*, X\left(t_{j}^*\right), \theta_1 \right)}{\partial \theta} \right\}\left[ 1+o_P(1) \right].\\
	\end{aligned}
	\end{equation*}
	
	$V_{22}$ and $V_{23}$ can be computed as under the null hypothesis. Similarly, we compute $\ddot{S}_n(\theta)$ as follows
	\begin{equation*}
	\begin{aligned}
	\ddot{S}_n(\tilde{\theta}_n)=&-\frac{2}{m} \sum_{j=1}^{m}\left\{\hat{X}^{\prime}_k\left(t_{j}^*\right)-X_k^{\prime}(t_j^*)\right\}\omega(t_j^*) \frac{\partial^{2} f_k\left(t_j^*,\hat{X}\left(t_{j}^*\right), \tilde{\theta}_n\right)}{\partial \theta \partial \theta^{\mathrm{T}}}\\
	&-\frac{2}{m} \sum_{j=1}^{m}\left\{f_k\left( t_j^*, X\left(t_{j}^*\right), \theta_0\right)+q(t_j^*)-f_k\left(t_j^*, \hat{X}\left(t_{j}^*\right), \tilde{\theta}_n\right)\right\}\\
	&\times \omega(t_j^*) \frac{\partial^{2} f_k\left(t_j^*,\hat{X}\left(t_{j}^*\right), \tilde{\theta}_n\right)}{\partial \theta \partial \theta^{\mathrm{T}}}\\
	&+\frac{2}{m} \sum_{j=1}^{m} \omega(t_j^*)\frac{\partial f_k\left( t_j^*, \hat{X}\left(t_{j}^*\right), \tilde{\theta}_n\right)}{\partial \theta}\frac{\partial f_k\left( t_j^*, \hat{X}\left(t_{j}^*\right), \tilde{\theta}_n\right)}{\partial \theta^{\mathrm{T}}}\\
	=&2E_{p^*} \left\{\omega(t)\frac{\partial f_k\left(t, \hat{X}\left(t\right), \theta_1\right)}{\partial \theta}\frac{\partial f_k\left(t,\hat{X}\left(t\right), \theta_1\right)}{\partial \theta^{\mathrm{T}}}\right\}\\
	&-2E_{p^*} \left\{ \left[ f_k\left(t, X\left(t\right), \theta_0\right)+q(t)-f_k\left(t,\hat{X}\left(t\right), \theta_1\right)\right]\omega(t) \frac{\partial^{2} f_k\left( t,\hat{X}\left(t\right), \theta_1\right)}{\partial \theta \partial \theta^{\mathrm{T}}} \right\}\\
	& +o_P(1)\\
	\equiv & 2G+o_P(1).
	\end{aligned}
	\end{equation*}
	
	Combining the above results, the linear approximation of $\sqrt{n}(\hat{\theta}_{n}-\theta_{1})$ is as
	\begin{equation}
	\begin{aligned}
	&\sqrt{n}(\hat{\theta}_{n}-\theta_{1})\\
	=&G^{-1}\frac{\sqrt{n}}{m} \sum_{j=1}^{m}\left[ \Delta_k(t_j^*) \omega(t_j^*)\frac{\partial f_k\left( t_j^*, X\left(t_{j}^*\right), \theta_1 \right)}{\partial \theta}\right.\\
	&\left.- \omega(t_j^*)\frac{\partial f_k\left( t_j^*, X\left(t_{j}^*\right), \theta_1 \right)}{\partial \theta} \frac{\partial f_k\left( t_j^*, X\left(t_{j}^*\right), \theta_1 \right)}{\partial X^\mathrm{T}}
	\Lambda\left(t_{i}^{*}\right)\right]\\
	&+ \frac{\sqrt{n}}{m} \sum_{j=1}^{m}  \left[ f_k\left( t_j^*, X\left(t_{j}^*\right), \theta_0 \right)-f_k\left( t_j^*, X\left(t_{j}^*\right), \theta_1 \right)+q(t_j^*)\right]\\
	&\times \omega(t_j^*)G^{-1}\frac{\partial f_k\left( t_j^*, X\left(t_{j}^*\right), \theta_1 \right)}{\partial \theta}+o_P(1)\\
	=&G^{-1}\frac{\sqrt{n}}{m} \sum_{j=1}^{m}\left[ \Delta_k(t_j^*) \omega(t_j^*)\frac{\partial f_k\left( t_j^*, X\left(t_{j}^*\right), \theta_1 \right)}{\partial \theta}\right.\\
	&\left.- \omega(t_i^*)\frac{\partial f_k\left( t_j^*, X\left(t_{j}^*\right), \theta_1 \right)}{\partial \theta} \frac{\partial f_k\left( t_j^*, X\left(t_{j}^*\right), \theta_1 \right)}{\partial X^\mathrm{T}}
	\Lambda\left(t_{i}^{*}\right)\right]\\
	&+\sqrt{n} G^{-1} E_{p^*} \left\{ \left[ X_k^{\prime}(t)-f_k\left( t, X(t),\theta_1 \right)\right] \omega(t)\frac{\partial f_k\left( t, X\left(t\right), \theta_1 \right)}{\partial \theta} \right\}
	+o_P(1)\\
	=&G^{-1}\frac{\sqrt{n}}{m} \sum_{j=1}^{m}\left[ \Delta_k(t_j^*) \omega(t_j^*)\frac{\partial f_k\left( t_j^*, X\left(t_{j}^*\right), \theta_1 \right)}{\partial \theta} \right.\\
	&\left.- \omega(t_j^*)\frac{\partial f_k\left( t_j^*, X\left(t_{j}^*\right), \theta_1 \right)}{\partial \theta} \frac{\partial f_k\left( t_j^*, X\left(t_{j}^*\right), \theta_1 \right)}{\partial X^\mathrm{T}}
	\Lambda\left(t_{j}^{*}\right)\right]
	+o_P(1).\label{A28} \\
	\end{aligned}
	\end{equation}
	The last equation uses the formula that
	\begin{equation*}
	\begin{aligned}
	0&=\left.\frac{\partial E_{p^*}\left\{ w(t) [X_k^{\prime}(t)-f_k\left( t,X(t), \theta \right)]^{2}\right\}}{\partial \theta}\right|_{\theta=\theta_1}\\
	&=-2 E_{p^*} \left\{ \left[ X_k^{\prime}(t)-f_k\left( t, X(t),\theta_1 \right)\right] \omega(t)\frac{\partial f_k\left( t, X\left(t\right), \theta_1 \right)}{\partial \theta} \right\}.
	\end{aligned}
	\end{equation*}
	
	Again, by computing the conditional variance of $\sqrt{n}(\hat{\theta}_{n}-\theta_{1})$ using the similar methods as under the null hypothesis, we can derive $\sqrt{n}(\hat{\theta}_{n}-\theta_{1})=O_P(n^{-1/2})$.
	
	{\bf 3. Under local alternative $H_{1kn}^f$}. Since
	\begin{equation*}
	\begin{aligned}
	&\lim_{\delta_n \rightarrow 0} E_{p^*}\left\{ \left[X_k^{\prime}(t)-f_k\left(t, X(t), \theta\right)\right]^{2}w(t)\right\}\\
	\geq& \lim_{\delta_n \rightarrow 0} E_{p^*}\left\{ \left[f_k\left(t, X(t), \theta_0\right)+\delta_n l_k(t)-f_k\left(t, X(t), \theta_0\right)\right]^{2}w(t)\right\}=0,\\
	\end{aligned}
	\end{equation*}
	we have $\theta^*=\theta_{0}$. Based on (\ref{A4}), we can derive
	\begin{equation*}
	\begin{aligned}
	\dot{S}_n(\theta^*)=&-\frac{2}{m} \sum_{j=1}^{m} \left\{f_k\left( t_j^*, X\left(t_{j}^*\right), \theta_0 \right)+\delta_n l_k(t_i^*)-f_k\left( t_j^*, \hat{X}\left(t_{j}^*\right), \theta_0 \right)\right\}\\
	&\times \omega(t_j^*)\frac{\partial f_k\left( t_j^*, \hat{X}\left(t_{j}^*\right), \theta_0 \right)}{\partial \theta}\\
	&-\frac{2}{m} \sum_{j=1}^{m} \left[\hat{X}_k^{\prime}(t_j^*)-X_k^{\prime}(t_j^*)\right]\omega(t_j^*) \frac{\partial f_k\left( t_j^*, \hat{X}\left(t_{j}^*\right), \theta_0 \right)}{\partial \theta}\\
	=&-\frac{2}{m} \sum_{j=1}^{m} \left\{f_k\left( t_j^*, X\left(t_{j}^*\right), \theta_0 \right)-f_k\left( t_j^*, \hat{X}\left(t_{j}^*\right), \theta_0 \right)\right\} \\
	&\times \omega(t_j^*)\frac{\partial f_k\left( t_j^*, \hat{X}\left(t_{j}^*\right), \theta_0 \right)}{\partial \theta}\\
	&-\frac{2}{m} \sum_{j=1}^{m} \left[\hat{X}_k^{\prime}(t_j^*)-X_k^{\prime}(t_j^*)\right]\omega(t_j^*) \frac{\partial f_k\left[\hat{X}\left(t_{i}^*\right), \theta_0 \right]}{\partial \theta}\\
	&-\frac{2}{m} \sum_{j=1}^{m} \delta_n l_k(t_j^*) \omega(t_j^*)\frac{\partial f_k\left( t_j^*, \hat{X}\left(t_{j}^*\right), \theta_0 \right)}{\partial \theta}\\
	=&-2V_{11}-2V_{12}-2V_{13}.
	\end{aligned}
	\end{equation*}
	
	The limiting properties of $V_{11}$ and $V_{12}$ are same as $V_{01}$ and $V_{02}$ under the null hypothesis. As for $V_{13}$,
	\begin{equation*}
	\begin{aligned}
	V_{13}=&\frac{1}{m} \sum_{j=1}^{m} \delta_n l_k(t_j^*) \omega(t_j^*)\frac{\partial f_k\left( t_j^*,X\left(t_{j}^*\right), \theta_0 \right)}{\partial \theta}\left[ 1+o_P(1) \right].
	\end{aligned}
	\end{equation*}
	The second derivative $\ddot{S}_n(\tilde{\theta}_n)$ can similarly be decomposed as
	\begin{equation*}
	\begin{aligned}
	\ddot{S}_n(\tilde{\theta}_n)=&-\frac{2}{m} \sum_{j=1}^{m}\left[ \hat{X}^{\prime}_k\left(t_{j}^*\right)-f_k\left( t_j^*,\hat{X}\left(t_{j}^*\right), \tilde{\theta}_n\right)\right]\omega(t_j^*) \frac{\partial^{2} f_k\left( t_j^*,\hat{X}\left(t_{j}^*\right), \tilde{\theta}_n\right)}{\partial \theta \partial \theta^{\mathrm{T}}}\\
	&+\frac{2}{m} \sum_{j=1}^{m} \omega(t_j^*)\frac{\partial f_k\left( t_j^*,\hat{X}\left(t_{j}^*\right), \tilde{\theta}_n\right)}{\partial \theta}\frac{\partial f_k\left( t_j^*,\hat{X}\left(t_{j}^*\right), \tilde{\theta}_n\right)}{\partial \theta^{\mathrm{T}}}\\
	=&-\frac{2}{m} \sum_{j=1}^{m}\left[\hat{X}^{\prime}_k\left(t_{j}^*\right)-X_k^{\prime}(t_j^*)\right]\omega(t_j^*) \frac{\partial^{2} f_k\left( t_j^*,\hat{X}\left(t_{j}^*\right), \tilde{\theta}_n\right)}{\partial \theta \partial \theta^{\mathrm{T}}}\\
	&-\frac{2}{m} \sum_{j=1}^{m}\left\{f_k\left( t_j^*,\hat{X}\left(t_{j}^*\right), \tilde{\theta}_n\right)+\delta_n l(t_j^*)-f_k\left( t_j^*,\hat{X}\left(t_{j}^*\right), \tilde{\theta}_n\right)\right\}\\
	&\times \omega(t_j^*) \frac{\partial^{2} f_k\left( t_j^*,\hat{X}\left(t_{j}^*\right), \tilde{\theta}_n\right)}{\partial \theta \partial \theta^{\mathrm{T}}}\\
	&+\frac{2}{m} \sum_{j=1}^{m} \omega(t_j^*)\frac{\partial f_k\left( t_j^*,\hat{X}\left(t_{j}^*\right), \tilde{\theta}_n\right)}{\partial \theta}\frac{\partial f_k\left( t_j^*,\hat{X}\left(t_{j}^*\right), \tilde{\theta}_n\right)}{\partial \theta^{\mathrm{T}}}\\
	=&2E_{p^*} \left[\omega(t)\frac{\partial f_k\left(t, \hat{X}\left(t\right), \tilde{\theta}_n\right)}{\partial \theta}\frac{\partial f_k\left(t, \hat{X}\left(t\right), \tilde{\theta}_n\right)}{\partial \theta^{\mathrm{T}}}\right]+o_P(1)\\
	=&2H_{\dot{f}}+o_P(1).
	\end{aligned}
	\end{equation*}
	Therefore, the linear approximation of $\sqrt{n}(\hat{\theta}_{n}-\theta_{0})$ is
	\begin{equation}
	\begin{aligned}
	&\sqrt{n}(\hat{\theta}_{n}-\theta_{0})\\
	=&H_{\dot{f}}^{-1}\frac{\sqrt{n}}{m} \sum_{j=1}^{m}\left[ \Delta_k(t_j^*) \omega(t_j^*)\frac{\partial f_k\left( t_j^*, X\left(t_{j}^*\right), \theta_0 \right)}{\partial \theta}\right.\\
	&\left.- \omega(t_j^*)\frac{\partial f_k\left( t_j^*, X\left(t_{j}^*\right), \theta_0 \right)}{\partial \theta} \frac{\partial f_k\left( t_j^*, X\left(t_{j}^*\right), \theta_0 \right)}{\partial X^\mathrm{T}}
	\Lambda\left(t_{j}^{*}\right)\right]\\
	&+\delta_n \frac{\sqrt{n}}{m} \sum_{j=1}^{m} l_k(t_j^*) \omega(t_j^*)H_{\dot{f}}^{-1}\frac{\partial f_k\left( t_j^*, X\left(t_{j}^*\right), \theta_0 \right)}{\partial \theta}
	+o_P(1)\\
	=&H_{\dot{f}}^{-1}\frac{\sqrt{n}}{m} \sum_{j=1}^{m}\left[ \Delta_k(t_j^*) \omega(t_j^*)\frac{\partial f_k\left( t_j^*, X\left(t_{j}^*\right), \theta_0 \right)}{\partial \theta} \right.\\
	&\left.- \omega(t_j^*)\frac{\partial f_k\left( t_j^*, X\left(t_{j}^*\right), \theta_0 \right)}{\partial \theta} \frac{\partial f_k\left( t_j^*, X\left(t_{j}^*\right), \theta_0 \right)}{\partial X^\mathrm{T}}
	\Lambda\left(t_{j}^{*}\right)\right]\\
	&+\sqrt{n} \delta_n H_{\dot{f}}^{-1} E_g \left[ l_k(t_j^*) \omega(t_j^*)\frac{\partial f_k\left( t_j^*, X\left(t_{j}^*\right), \theta_0 \right)}{\partial \theta} \right]
	+o_P(1).\label{A16}\\
	\end{aligned}
	\end{equation}
	The root-$n$ consistency can also be derived by computing its variance as before.
	
	Further, the proof under $H_{1kn}^F$ is very similar to that under $H_{1kn}^f$,  we omit the detail.
\end{proof}	

\section{The results in Section~4}
\subsection{Proof of Lemma~5}

\begin{proof}
	This proof is an extension of the proof of Lemma~3.1 of \cite{Powell.1989}. To prove $\sqrt{n}(U_n-\hat{U}_n)=\mathrm{o}_P(1)$, it is sufficient to show $n E\left[\left\|U_{n}-\hat{U}_{n}\right\|^{2}\right]=\mathrm{o}(1)$. Write $r_{n}\left(z_{i}\right)=E\left(h_{n}| z_{i}\right)$ and $\theta_{n}=E\left[r_{n}\left(z_{i}\right)\right]=E\left[h_{n}\left(z_{1},\cdots,z_{m^*}\right)\right]$. Define
	\begin{equation*}
	q_{n}\left(z_{1},\cdots,z_{m^*}\right)=\left[h_{n}\left(z_{1},\cdots,z_{m^*}\right)-\sum_{i=1}^{m^*}r_{n}\left(z_{i}\right)+(m^*-1)\theta_{n}\right],
	\end{equation*}
	so that
	\begin{equation*}
	U_{n}-\hat{U}_{n}=\left(\begin{array}{c}{n} \\ {m^*}\end{array}\right)^{-1} \sum_{c} q_{n}\left(z_{i_1},\cdots,z_{i_{m^*}}\right) .
	\end{equation*}
	
	The expectation of the squared length of the vector $U_{n}-\hat{U}_{n}$ is \begin{equation*}
	E\left[\left\|U_{n}-\hat{U}_{n}\right\|^{2}\right]=\left(\begin{array}{c}{n} \\ {m^*}\end{array}\right)^{-2} \sum_{c_1} \sum_{c_2} E\left[q_{n}\left(z_{i_1},\cdots,z_{i_{m^*}}\right)^{\prime} q_{n}\left(z_{j_1},\cdots,z_{j_{m^*}}\right)\right].
	\end{equation*}
	
	It is easy to show that if $E\left[q_{n}\left(z_{i_1},\cdots,z_{i_{m^*}}\right)^{\prime} q_{n}\left(z_{j_1},\cdots,z_{j_{m^*}}\right)\right] \neq 0$, there are at least two same terms in $q_{n}\left(z_{i_1},\cdots,z_{i_{m^*}}\right)^{\prime}$ and $q_{n}\left(z_{j_1},\cdots,z_{j_{m^*}}\right)$. For example $i_1=j_1$ and $i_2=j_2$. Thus the number of nonzero terms in the sum is only of  order $\mathrm{O}(n^{2m^*-2})$ instead of $\mathrm{O}(n^{2m^*})$. Each nonzero term can be shown to be $\mathrm{o}(n)$ according to the condition. Consequently we have
	\begin{equation*}
	\begin{aligned} n E\left[\left\|U_{n}-\hat{U}_{n}\right\|^{2}\right] &=N\left(\begin{array}{c}{n} \\ {m^*}\end{array}\right)^{-2} \mathrm{O}\left(n^{2m^*-2}\right) \mathrm{o}(n) \\ &=\mathrm{o}(1), \end{aligned}
	\end{equation*}
	which is what we need.
	
\end{proof}

\subsection{Proof of Lemma~6}

\begin{proof}
	Denote
	\begin{equation*}
	\begin{aligned}
	\hat{e}_{ik} &=Y_{ik}-\hat{F}_k(t_i;\hat{\theta})\\
	&= Y_{ik}-X_k(t_0)-\int_{t_0}^{t_i} f_k\left(t,\hat{X}\left(t\right) ; \hat{\theta}\right)\mathrm{d}t .
	\end{aligned}
	\end{equation*}
	In the remaining part of this proof, we omit the subscript $k$ for notational simplicity.
	
	Decompose $V_{n}$ as
	\begin{equation*}
	\begin{aligned}
	V_{n}
	=& \frac{1}{n(n-1)} \sum_{i=1}^{n} \sum_{j=1 \atop j \neq i}^{n} \frac{1}{h} K\left(\frac{t_{i}-t_{j}}{h}\right) \hat{e}_{i} \hat{e}_{j}\\
	=& \frac{1}{n(n-1)} \sum_{i=1}^{n} \sum_{j=1 \atop j \neq i}^{n} \frac{1}{h} K\left(\frac{t_{i}-t_{j}}{h}\right)[\varepsilon_i+F(t_i;\theta_0)-\hat{F}(t_i;\hat{\theta})][\varepsilon_j+F(t_j;\theta_0)-\hat{F}(t_j;\hat{\theta})]\\
	=& \frac{1}{n(n-1)} \sum_{i=1}^{n} \sum_{j=1 \atop j \neq i}^{n} \frac{1}{h} K\left(\frac{t_{i}-t_{j}}{h}\right)\varepsilon_i\varepsilon_j\\
	&-\frac{2}{n(n-1)} \sum_{i=1}^{n} \sum_{j=1 \atop j \neq i}^{n} \frac{1}{h} K\left(\frac{t_{i}-t_{j}}{h}\right)\varepsilon_i[\hat{F}(t_j;\hat{\theta})-F(t_j;\theta_0)]\\
	&+\frac{1}{n(n-1)} \sum_{i=1}^{n} \sum_{j=1 \atop j \neq i}^{n} \frac{1}{h} K\left(\frac{t_{i}-t_{j}}{h}\right)[\hat{F}(t_i;\hat{\theta})-F(t_i;\theta_0)][\hat{F}(t_j;\hat{\theta})-F(t_j;\theta_0)]\\
	\equiv &V_{1}-2V_{2}+V_{3}.
	\end{aligned}
	\end{equation*}
	Now to prove that $nh^{1/2}V_{2}$ and $nh^{1/2}V_3$ are $o_P(1)$ and then $nh^{1/2} V_n=nh^{1/2}V_{1}+o_p(1)$. \\
	
	{\bf Step 1.} Consider $V_{2}$ first which has the following:
	\begin{equation}
	\begin{aligned}
	V_{2}=&\frac{1}{n(n-1)} \sum_{i=1}^{n} \sum_{j=1 \atop j \neq i}^{n} \frac{1}{h} K\left(\frac{t_{i}-t_{j}}{h}\right)\varepsilon_i[\hat{F}(t_j;\hat{\theta})-F(t_j;\theta_0)]\\
	=&\frac{1}{n(n-1)} \sum_{i=1}^{n} \sum_{j=1 \atop j \neq i}^{n} \frac{1}{h} K\left(\frac{t_{i}-t_{j}}{h}\right)\varepsilon_i[\hat{F}(t_j;\hat{\theta})+F(t_j;\hat{\theta})-F(t_j;\hat{\theta})-F(t_j;\theta_0)]\\
	=&\frac{1}{n(n-1)} \sum_{i=1}^{n} \sum_{j=1 \atop j \neq i}^{n} \frac{1}{h} K\left(\frac{t_{i}-t_{j}}{h}\right)\varepsilon_i[F(t_j;\hat{\theta})-F(t_j;\theta_0)]\\
	&+\frac{1}{n(n-1)} \sum_{i=1}^{n} \sum_{j=1 \atop j \neq i}^{n} \frac{1}{h} K\left(\frac{t_{i}-t_{j}}{h}\right)\varepsilon_i[\hat{F}(t_j;\hat{\theta})-F(t_j;\hat{\theta})]\\
	\equiv &V_{21}+V_{22}.\label{A34}
	\end{aligned}
	\end{equation}
	
	For $V_{21}$, we have
	\begin{equation}
	\begin{aligned}
	V_{21}=&\frac{1}{n(n-1)} \sum_{i=1}^{n} \sum_{j=1 \atop j \neq i}^{n} \frac{1}{h} K\left(\frac{t_{i}-t_{j}}{h}\right)\varepsilon_i[F(t_j;\hat{\theta})-F(t_j;\theta_0)]\\
	=&\frac{1}{n(n-1)} \sum_{i=1}^{n} \sum_{j=1 \atop j \neq i}^{n} \frac{1}{h} K\left(\frac{t_{i}-t_{j}}{h}\right)\varepsilon_i\frac{\partial F\left(t_{j}; \theta_{0}\right)}{\partial \theta^{\mathrm{T}}}\left(\hat{\theta}-\theta_{0}\right)[1+o_P(1)],\label{A33}\\
	\end{aligned}
	\end{equation}
	which can be proven to be $O_P(n^{-1})$ using Lemma~3.3b of \cite{Zheng.1996}.
	
	For $V_{22}$, we have
	\begin{equation}
	\begin{aligned}
	&V_{22}\\
	=&\frac{1}{n(n-1)} \sum_{i=1}^{n} \sum_{j=1 \atop j \neq i}^{n} \frac{1}{h} K\left(\frac{t_{i}-t_{j}}{h}\right)\varepsilon_i[\hat{F}(t_j;\hat{\theta})-F(t_j;\hat{\theta})]\\
	=&\frac{1}{n(n-1)} \sum_{i=1}^{n} \sum_{j=1 \atop j \neq i}^{n} \frac{1}{h} K\left(\frac{t_{i}-t_{j}}{h}\right)\varepsilon_i \int_{t_0}^{t_j} \left\{f[t,\hat{X}\left(t\right) ; \hat{\theta}] - f[t,X\left(t\right) ; \hat{\theta}]\right\}\mathrm{d}t \\
	=&\frac{1}{n(n-1)} \sum_{i=1}^{n} \sum_{j=1 \atop j \neq i}^{n} \frac{1}{h} K\left(\frac{t_{i}-t_{j}}{h}\right)\varepsilon_i \int_{t_0}^{t_j} \left\{ \frac{\partial f[t,X\left(t\right) ; \hat{\theta}]}{\partial X}[\hat{X}\left(t\right)-X\left(t\right)] \right\}\mathrm{d}t \\
	=&\frac{1}{n(n-1)} \sum_{i=1}^{n} \sum_{j=1 \atop j \neq i}^{n} \frac{1}{h} K\left(\frac{t_{i}-t_{j}}{h}\right)\varepsilon_i\int_{t_0}^{t_j} \left\{ \frac{\partial f[t,X\left(t\right) ; \hat{\theta}]}{\partial X}\right.\\
	&\left.\left[\frac{\frac{1}{n^2h^2}\sum_k^n \sum_l^n \left[(\frac{t-t_k}{h})^2K(\frac{t-t_k}{h})K(\frac{t-t_l}{h})Y_l - \frac{t-t_k}{h}K(\frac{t-t_k}{h}) \frac{t-t_l}{h}K(\frac{t-t_l}{h})Y_l \right]}{\frac{1}{n^2h^2}\sum_k^n \sum_l^n \left[(\frac{t-t_k}{h})^2K(\frac{t-t_k}{h})K(\frac{t-t_l}{h}) - \frac{t-t_k}{h}K(\frac{t-t_k}{h}) \frac{t-t_l}{h}K(\frac{t-t_l}{h})\right]}-X\left(t\right)\right] \right\}\mathrm{d}t \\
	=&\frac{1}{n(n-1)} \sum_{i=1}^{n} \sum_{j=1 \atop j \neq i}^{n} \frac{1}{h} K\left(\frac{t_{i}-t_{j}}{h}\right)\varepsilon_i\int_{t_0}^{t_j} \left\{ \frac{\partial f[t,X\left(t\right) ; \hat{\theta}]}{\partial X}\frac{N_n(t)-X\left(t\right)M_n(t)}{M_n(t)} \right\}\mathrm{d}t \\
	=&\frac{1}{n(n-1)} \sum_{i=1}^{n} \sum_{j=1 \atop j \neq i}^{n} \frac{1}{h} K\left(\frac{t_{i}-t_{j}}{h}\right)\varepsilon_i\int_{t_0}^{t_j} \left\{ \frac{\partial f[t,X\left(t\right) ; \hat{\theta}]}{\partial X}\right.\\
	&\left.\left[\frac{N_n(t)-X\left(t\right)M_n(t)}{M_n(t)}+\frac{N_n(t)-X\left(t\right)M_n(t)}{M(t)}-\frac{N_n(t)-X\left(t\right)M_n(t)}{M(t)}\right] \right\}\mathrm{d}t \\
	=&\frac{1}{n(n-1)} \sum_{i=1}^{n} \sum_{j=1 \atop j \neq i}^{n} \frac{1}{h} K\left(\frac{t_{i}-t_{j}}{h}\right)\varepsilon_i\int_{t_0}^{t_j} \left\{ \frac{\partial f[t,X\left(t\right) ; \hat{\theta}]}{\partial X}\frac{N_n(t)-X\left(t\right)M_n(t)}{M(t)} \right\}\mathrm{d}t \\
	&+\frac{1}{n(n-1)} \sum_{i=1}^{n} \sum_{j=1 \atop j \neq i}^{n} \frac{1}{h} K\left(\frac{t_{i}-t_{j}}{h}\right)\varepsilon_i\int_{t_0}^{t_j} \left\{ \frac{\partial f[t,X\left(t\right) ; \hat{\theta}]}{\partial X}\right.\\
	&\left.\left[\frac{N_n(t)-X\left(t\right)M_n(t)}{M_n(t)}-\frac{N_n(t)-X\left(t\right)M_n(t)}{M(t)}\right] \right\}\mathrm{d}t \\
	\equiv &V_{221}+V_{222}, \label{A32}
	\end{aligned}
	\end{equation}
	where $M_n(t)$, $N_n(t)$ and $M(t)$ are defined in Lemma \ref{LP.2}.
	
	Then, $V_{221}$ can be written in a form of the difference between two V-statistics:
	\begin{equation*}
	V_{221}=\frac{n}{n-1}V_{221}^1-\frac{1}{n-1}V_{222}^2
	\end{equation*}
	where
	\begin{equation*}
	V_{221}^1=\frac{1}{n^2} \sum_{i=1}^{n} \sum_{j=1}^{n} \frac{1}{h} K\left(\frac{t_{i}-t_{j}}{h}\right)\varepsilon_i\int_{t_0}^{t_j} \left\{ \frac{\partial f[t,X\left(t\right) ; \hat{\theta}]}{\partial X}\frac{N_n(t)-X\left(t\right)M_n(t)}{M(t)} \right\}\mathrm{d}t
	\end{equation*}
	and
	\begin{equation*}
	V_{221}^2=\frac{1}{n} \sum_{j=1}^{n} \frac{1}{h} K\left(0\right)\varepsilon_j\int_{t_0}^{t_j} \left\{ \frac{\partial f[t,X\left(t\right) ; \hat{\theta}]}{\partial X}\frac{N_n(t)-X\left(t\right)M_n(t)}{M(t)} \right\}\mathrm{d}t.
	\end{equation*}
	For $V_{221}^1$, the corresponding kernel function is
	\begin{equation*}
	\begin{aligned}
	&H_n(z_i,z_j,z_k,z_l)\\
	=& \frac{1}{24} \sum_{P} \frac{1}{h h_0^2} K\left(\frac{t_{i}-t_{j}}{h}\right)\varepsilon_i \\
	&\int_{t_0}^{t_j} \frac{1}{\mu_2(K) p^2(t)} \frac{\partial f[t,X\left(t\right) ; \theta_0]}{\partial X}[ (\frac{t-t_k}{h})^2K(\frac{t-t_k}{h})K(\frac{t-t_l}{h})(Y_l-X\left(t\right)) \\
	&- \frac{t-t_k}{h}K(\frac{t-t_k}{h}) \frac{t-t_l}{h}K(\frac{t-t_l}{h})(Y_l-X\left(t\right))] \mathrm{d}t.\\
	\end{aligned}
	\end{equation*}
	Its second order moment  is
	\begin{equation*}
	\begin{aligned}
	&E\left[H_n^2(z_i,z_j,z_k,z_l)\right]\\
	\leq & E\{ \frac{1}{h^2 h_0^4} K^2\left(\frac{t_{i}-t_{j}}{h}\right)\varepsilon_i^2 \\
	&\{\int_{t_0}^{t_j} \frac{1}{\mu_2(K) p^2(t)} \frac{\partial f[t,X\left(t\right) ; \theta_0]}{\partial X}[ (\frac{t-t_k}{h})^2K(\frac{t-t_k}{h})K(\frac{t-t_l}{h})(Y_l-X\left(t\right)) \\
	&- \frac{t-t_k}{h}K(\frac{t-t_k}{h}) \frac{t-t_l}{h}K(\frac{t-t_l}{h})(Y_l-X\left(t\right))] \mathrm{d}t \}^2\\
	=&O(\frac{1}{hh_0^2})=o(n).
	\end{aligned}
	\end{equation*}
	Thus the condition of Theorem~1 of \cite{MartinsFilho.2006} is satisfied and the limiting distribution of $V_{221}^1$  is equivalent to the relevant U-statistic. The application of Lemma~5 can yield the limiting distribution of this V-statistic by computing the projection of the relevant U-statistic.
	
	The conditional expectation of $H_n$ given $z_i$ has the following result:
	\begin{equation}
	\begin{aligned}
	&r(z_i)\\
	=&\mathrm{E}(V_{221}|z_i)\\
	=&\frac{1}{4 h h_0^2} \int \int \int K\left(\frac{t_{i}-t_{j}}{h}\right)\varepsilon_i \\
	&\int_{t_0}^{t_j} \frac{1}{\mu_2(K) p^2(t)} \frac{\partial f[t,X\left(t\right) ; \theta_0]}{\partial X}[ (\frac{t-t_k}{h})^2K(\frac{t-t_k}{h})K(\frac{t-t_l}{h})(X_l-X(t)) \\
	&- \frac{t-t_k}{h}K(\frac{t-t_k}{h}) \frac{t-t_l}{h}K(\frac{t-t_l}{h})(X_l-X(t))] \mathrm{d}t p(t_k)p(t_j)p(t_l) \mathrm{d}t_j\mathrm{d}t_l\mathrm{d}t_k\\
	=&\frac{1}{4 h h_0^2} \int \int \int \int_{t_0}^{t_j} K\left(\frac{t_{i}-t_{j}}{h}\right)\varepsilon_i \\
	& \frac{1}{\mu_2(K) p^2(t)} \frac{\partial f[t,X\left(t\right) ; \theta_0]}{\partial X}[ (\frac{t-t_k}{h})^2K(\frac{t-t_k}{h})K(\frac{t-t_l}{h})(X_l-X(t)) \\
	&- \frac{t-t_k}{h}K(\frac{t-t_k}{h}) \frac{t-t_l}{h}K(\frac{t-t_l}{h})(X_l-X(t))] p(t_k)p(t_j)p(t_l) \mathrm{d}t\mathrm{d}t_j\mathrm{d}t_l\mathrm{d}t_k\\
	=&\frac{1}{4}  \int \int_{t_0}^{t_i+u_jh}\int \int K\left(u_j\right)\varepsilon_i \frac{1}{\mu_2(K) p^2(t)} \frac{\partial f[t,X\left(t\right) ; \theta_0]}{\partial X}[ u_k^2K(u_k)K(u_l)\\
	&\times(X(t-u_lh_0)-X(t))
	p(t-u_kh_0)p(t-u_lh_0)  - u_k u_l K(u_k) K(u_l)\\
	&\times (X(t-u_lh_0)-X(t))p(t-u_kh_0)p(t-u_lh_0)] p(t_i+u_jh) \mathrm{d}u_l\mathrm{d}u_k\mathrm{d}t\mathrm{d}u_j\\
	=&\frac{h_0^2}{8}  \int \int_{t_0}^{t_i+u_jh} K\left(u_j\right)\varepsilon_i \frac{1}{\mu_2(K) p^2(t)} \frac{\partial f[t,X\left(t\right) ; \theta_0]}{\partial X}[2X^{\prime}(t)p^{\prime}(t)[\mu_2(K)]^2p(t)\\
	&+X^{(2)}(t)[\mu_2(K)]^2p^2(t)-2X^{\prime}(t)p^{\prime}(t)[\mu_2(K)]^2p(t)]\\ &p(t_i+u_jh)\mathrm{d}t\mathrm{d}u_j+\mathrm{o}_P(h_0^2)\\
	=&\frac{h_0^2}{8}  \int  K\left(u_j\right)\varepsilon_i p(t_i+u_jh) \int_{t_0}^{t_i+u_jh} \frac{\partial f[t,X\left(t\right) ; \theta_0]}{\partial X}\\
	&X^{(2)}(t)\mu_2(K) \mathrm{d}t\mathrm{d}u_j+\mathrm{o}_P(h_0^2)\\
	=&\frac{h_0^2}{8}  \int  K\left(u_j\right)\varepsilon_i p(t_i+u_jh) \left[R(t_i+u_jh)-R(t_0)\right]\mathrm{d}u_j+\mathrm{o}_P(h_0^2)\\
	=&\frac{h_0^2}{8} \varepsilon_i p(t_i) [R(t_i)-R(t_0)] + \mathrm{o}_P(h_0^2) + \mathrm{O}_P(h^2h_0^2).\label{A12} \\
	\end{aligned}
	\end{equation}
	
	Using this conditional expectation, we can obtain the limiting distribution of the projection of the relevant U-statistic. Thus, according to Lemma~5, we know that $V_{221}^1=\mathrm{O}_P(\frac{h_0^2}{\sqrt{n}})$ and then $nh^{1/2} [n/(n-1)]V_{221}^1=\mathrm{O}_P(h_0^2\sqrt{nh})=o_P(1)$ under the condition that $h_0=\mathrm{o}(n^{-1/4}h^{-1/4})$. Using similar method for the V-statistic $V_{221}^2$, it can be proven that $nh^{1/2} [1/(n-1)] V_{221}^2$ is also $\mathrm{o}_P(1)$. Therefore,
	\begin{equation}
	nh^{\frac{1}{2}}V_{221}=nh^{\frac{1}{2}}\left(\frac{n}{n-1}V_{221}^1-\frac{1}{n-1}V_{222}^2\right)=\mathrm{o}_P(1). \label{A30}
	\end{equation}
	Turn to $V_{222}$. We have
	\begin{equation}
	\begin{aligned}
	&V_{222}\\
	=&\frac{1}{n(n-1)} \sum_{i=1}^{n} \sum_{j=1 \atop j \neq i}^{n} \frac{1}{h} K\left(\frac{t_{i}-t_{j}}{h}\right)\varepsilon_i\int_{t_0}^{t_j} \left\{ \frac{\partial f[t,X\left(t\right) ; \hat{\theta}]}{\partial X} \right.\\
	&\left.[N_n(t)-X\left(t\right)M_n(t)][\frac{1}{M_n(t)}-\frac{1}{M(t)}] \right\}\mathrm{d}t \\
	=&\frac{1}{n(n-1)} \sum_{i=1}^{n} \sum_{j=1 \atop j \neq i}^{n} \frac{1}{h} K\left(\frac{t_{i}-t_{j}}{h}\right)\varepsilon_i\int_{t_0}^{t_j} \left\{ \frac{\partial f[t,X\left(t\right) ; \hat{\theta}]}{\partial X} \right.\\
	&\left.[N_n(t)-X\left(t\right)M_n(t)]\frac{M(t)-M_n(t)}{M(t)M_n(t)} \right\}\mathrm{d}t \\
	=&\frac{1}{n(n-1)} \sum_{i=1}^{n} \sum_{j=1 \atop j \neq i}^{n} \frac{1}{h} K\left(\frac{t_{i}-t_{j}}{h}\right)\varepsilon_i\int_{t_0}^{t_j} \left\{ \frac{\partial f[t,X\left(t\right) ; \hat{\theta}]}{\partial X} \right.\\
	&\left. [\frac{N_n(t)}{M_n(t)}-X\left(t\right)]\frac{M(t)-M_n(t)}{M(t)} \right\}\mathrm{d}t \\
	\leq& \frac{1}{n(n-1)} \sum_{i=1}^{n} \sum_{j=1 \atop j \neq i}^{n} \frac{1}{h} K\left(\frac{t_{i}-t_{j}}{h}\right)\varepsilon_i\int_{t_0}^{t_j} \left\{ \left\| \frac{\partial f[t,X\left(t\right) ; \hat{\theta}]}{\partial X} \right\|\right. \\
	&\left. \left[\sup_t \left\| X_n(t)-X\left(t\right) \right\| \right] \frac{\left[\sup_t \left| M_n(t)-M\left(t\right) \right| \right]}{\left|M(t)\right|}\right\}\mathrm{d}t\\
	=&\mathrm{O}_P[a_n^2(h_0)].\label{A13}
	\end{aligned}
	\end{equation}
	
	Under the condition that $a_n^2(h_0)=o(n^{-1}h^{-1/2})$, we have
	\begin{equation}
	nh^{\frac{1}{2}}V_{222}=nh^{\frac{1}{2}}\mathrm{O}_P[a_n^2(h_0)]=\mathrm{o}_P(1). \label{A31}
	\end{equation}
	Combining  (\ref{A34}), (\ref{A33}), (\ref{A32}), (\ref{A30}) and (\ref{A31}), we conclude
	\begin{equation}
	nh^{1/2}V_{2}=o_P(1). \label{A37}
	\end{equation}
	
	{\bf Step 2.} Deal with $V_3$, which  can be decomposed as
	\begin{equation}
	\begin{aligned}
	V_{3} =& \frac{1}{n(n-1)} \sum_{i=1}^{n} \sum_{j=1 \atop j \neq i}^{n} \frac{1}{h} K\left(\frac{t_{i}-t_{j}}{h}\right)[\hat{F}(t_i;\hat{\theta})-F(t_i;\theta_0)][\hat{F}(t_j;\hat{\theta})-F(t_j;\theta_0)]\\
	=& \frac{1}{n(n-1)} \sum_{i=1}^{n} \sum_{j=1 \atop j \neq i}^{n} \frac{1}{h} K\left(\frac{t_{i}-t_{j}}{h}\right)[\hat{F}(t_i;\hat{\theta})-F(t_i;\hat{\theta})+F(t_i;\hat{\theta})-F(t_i;\theta_0)]\\
	&\times[\hat{F}(t_j;\hat{\theta})-F(t_j;\hat{\theta})+F(t_j;\hat{\theta})-F(t_j;\theta_0)]\\
	=&\frac{1}{n(n-1)} \sum_{i=1}^{n} \sum_{j=1 \atop j \neq i}^{n} \frac{1}{h} K\left(\frac{t_{i}-t_{j}}{h}\right)[F(t_i;\hat{\theta})-F(t_i;\theta_0)][F(t_j;\hat{\theta})-F(t_j;\theta_0)]\\
	&+\frac{2}{n(n-1)} \sum_{i=1}^{n} \sum_{j=1 \atop j \neq i}^{n} \frac{1}{h} K\left(\frac{t_{i}-t_{j}}{h}\right)[\hat{F}(t_i;\hat{\theta})-F(t_i;\hat{\theta})][F(t_j;\hat{\theta})-F(t_j;\theta_0)]\\
	&+\frac{1}{n(n-1)} \sum_{i=1}^{n} \sum_{j=1 \atop j \neq i}^{n} \frac{1}{h} K\left(\frac{t_{i}-t_{j}}{h}\right)[\hat{F}(t_i;\hat{\theta})-F(t_i;\hat{\theta})][\hat{F}(t_j;\hat{\theta})-F(t_j;\hat{\theta})]\\
	\equiv &V_{31}+2V_{32}+V_{33}. \label{A35}
	\end{aligned}
	\end{equation}
	Discuss them one by one. For $V_{31}$, we have
	\begin{equation}
	\begin{aligned}
	V_{31} =& \frac{1}{n(n-1)} \sum_{i=1}^{n} \sum_{j=1 \atop j \neq i}^{n} \frac{1}{h} K\left(\frac{t_{i}-t_{j}}{h}\right)[F(t_i;\hat{\theta})-F(t_i;\theta_0)][F(t_j;\hat{\theta})-F(t_j;\theta_0)]\\
	=&\left(\hat{\theta}-\theta_{0}\right)^{\mathrm{T}} \frac{1}{n(n-1)} \sum_{i=1}^{n} \sum_{j=1 \atop j \neq i}^{n} \frac{1}{h} K\left(\frac{t_{i}-t_{j}}{h}\right)\\
	&\times \frac{\partial F\left(t_{i}; \theta_{0}\right)}{\partial \theta} \frac{\partial F\left(t_{j}; \theta_{0}\right)}{\partial \theta^{\mathrm{T}}}\left(\hat{\theta}-\theta_{0}\right)[1+\mathrm{o}_P(1)]\\
	=&\mathrm{O}_P(n^{-1}).\label{A14}
	\end{aligned}
	\end{equation}
	
	For $V_{32}$, we have
	\begin{equation}
	\begin{aligned}
	V_{32} =& \frac{1}{n(n-1)} \sum_{i=1}^{n} \sum_{j=1 \atop j \neq i}^{n} \frac{1}{h} K\left(\frac{t_{i}-t_{j}}{h}\right)[\hat{F}(t_i;\hat{\theta})-F(t_i;\hat{\theta})][F(t_j;\hat{\theta})-F(t_j;\theta_0)]\\
	=&\left(\hat{\theta}-\theta_{0}\right)^{\mathrm{T}} \frac{1}{n(n-1)} \sum_{i=1}^{n} \sum_{j=1 \atop j \neq i}^{n} \frac{1}{h} K\left(\frac{t_{i}-t_{j}}{h}\right)\\
	&\times \frac{\partial F\left(t_{j}; \theta_{0}\right)}{\partial \theta} \left\{ \int_{t_0}^{t_j}  \frac{\partial f[t,X\left(t\right) ; \hat{\theta}]}{\partial X}[\hat{X}\left(t\right)-X\left(t\right)] \mathrm{d}t \right\} [1+\mathrm{o}_P(1)]\\
	\leq &\left\|\hat{\theta}-\theta_{0}\right\| \frac{1}{n(n-1)} \sum_{i=1}^{n} \sum_{j=1 \atop j \neq i}^{n} \frac{1}{h} K\left(\frac{t_{i}-t_{j}}{h}\right)(t_i-t_0) \left\| \frac{\partial F\left(t_{j}; \theta_{0}\right)}{\partial \theta}\right\| \\
	&\times \left[ \sup_t \left\| \frac{ \partial f[t,X\left(t\right) ; \hat{\theta}]}{\partial X}\right\| \right]\left[\sup_t \left\|\hat{X}\left(t\right)-X\left(t\right)\right\|\right] [1+\mathrm{o}_P(1)]\\
	=&\mathrm{O}_P[a_n(h_0)n^{-1/2}]. \label{A38}
	\end{aligned}
	\end{equation}
	Similarly, for $V_{33}$,
	\begin{equation}
	\begin{aligned}
	V_{33} =& \frac{1}{n(n-1)} \sum_{i=1}^{n} \sum_{j=1 \atop j \neq i}^{n} \frac{1}{h} K\left(\frac{t_{i}-t_{j}}{h}\right)[\hat{F}(t_i;\hat{\theta})-F(t_i;\hat{\theta})][\hat{F}(t_j;\hat{\theta})-F(t_j;\hat{\theta})]\\
	=& \frac{1}{n(n-1)} \sum_{i=1}^{n} \sum_{j=1 \atop j \neq i}^{n} \frac{1}{h} K\left(\frac{t_{i}-t_{j}}{h}\right) \left\{ \int_{t_0}^{t_i}  \frac{\partial f[t,X\left(t\right) ; \hat{\theta}]}{\partial X}[\hat{X}\left(t\right)-X\left(t\right)] \mathrm{d}t \right\}\\
	&\times \left\{ \int_{t_0}^{t_j}  \frac{\partial f[t,X\left(t\right) ; \hat{\theta}]}{\partial X}[\hat{X}\left(t\right)-X\left(t\right)] \mathrm{d}t \right\}(1+o_P(1))\\
	\leq& \frac{1}{n(n-1)} \sum_{i=1}^{n}  \sum_{j=1 \atop j \neq i}^{n} \frac{1}{h} K\left(\frac{t_{i}-t_{j}}{h}\right) (t_i-t_0)(t_j-t_0) \\
	&\times \left[ \sup_t \left\| \frac{ \partial f[t,X\left(t\right) ; \hat{\theta}]}{\partial X}\right\| \right]^2\left[\sup_t \left\|\hat{X}\left(t\right)-X\left(t\right)\right\|\right]^2[1+\mathrm{o}_P(1)]\\
	=&\mathrm{O}_P[a_n^2(h_0)].\label{A15}
	\end{aligned}
	\end{equation}
	Summarizing the results in (\ref{A35})-(\ref{A15}) and reminding of the conditions on bandwidths, we can show that
	\begin{equation}
	nh^{1/2}V_{3}=o_P(1). \label{A36}
	\end{equation}
	
	Together (\ref{A37}) with (\ref{A36}), the lemma is proved.
	
\end{proof}

\subsection{Proof of Theorem~4.1 (under the null hypothesis)}

\begin{proof}
	The result is an easy consequence of Lemma~6 and Theorem~1 of \cite{Zheng.1996}.
\end{proof}

\subsection{Proof of Theorem~4.2 (under global alternatives)}

\begin{proof}
	The result is an easy consequence of Lemma~6 and Theorem~2 of \cite{Zheng.1996}.
\end{proof}

\subsection{Proof of Theorem~4.3 (under local alternatives)}

\begin{proof}
	According to Lemma~6, we only need to study the convergence result of $V_{1n}$. The derivation of it is similar as in the proof of Theorem~2.4, thus we omit the detailed proof. Notice that since here we use the two-step collocation estimator instead of nonlinear least squares estimator, the linear approximation (\ref{A16}) should be used to replace (\ref{A17}).
\end{proof}
\section{The results in Section~5}
\subsection{Proof of Lemma~7}

\begin{proof}
	
	In this proof we use $f_k(t)$ to write $f_k(t, X(t); \theta_0)$ and $\hat{f}_k(t)$ to  $f_k(t, \hat{X}(t); \hat{\theta})$ for notational simplicity. The preudo-residual is
	\begin{equation*}
	\hat{e}_f(t_d)=\hat{X}_k^{\prime}(t_d)-\hat{f}_k(t_d)= \frac{\hat{h}_k^{\prime}(x) \hat{p}(t_d)-\hat{h}_k(t_d) \hat{p}^{\prime}(t_d)-\hat{p}^2\left(t_{d}\right)\hat{f}_k(t_d)}{\hat{p}^2\left(t_{d}\right)}.
	\end{equation*}
	
	We then have
	\begin{equation}
	\begin{aligned}
	V_{n}^f=&\frac{1}{nh^2}\sum^{n}_{d=1} \left\{ \left[\hat{X}_k^{\prime}(t_d)-\hat{f}_k(t_d) \right]\hat{p}^2\left(t_{d}\right)\right\}^2 \\
	=&\frac{1}{nh^2}\sum^{n}_{d=1} \left\{\left[ \hat{X}^{\prime}_k(t_d)-\hat{f}_k(t_d)+f_k(t_d)-f_k(t_d) \right] \hat{p}^2\left(t_{d}\right)\right\}^2 \\
	=&\frac{1}{nh^2}\sum^{n}_{d=1} \left\{\left[\hat{X}^{\prime}_k(t_d)-f_k(t_d)\right]\hat{p}^2\left(t_{d}\right)\right\}^2 + \frac{1}{nh^2}\sum^{n}_{d=1} \left\{ \left[\hat{f}_k(t_d)-f_k(t_d)\right]\hat{p}^2\left(t_{d}\right)\right\}^2\\
	&-\frac{2}{nh^2}\sum^{n}_{d=1} \left[\hat{X}^{\prime}_k(t_d)-f_k(t_d)\right]\left[\hat{f}_k(t_d)-f_k(t_d)\right]\hat{p}^2\left(t_{d}\right)\\
	\equiv &V_{1n}^f+V_{3n}^f-2V_{2n}^f .\label{A21}
	\end{aligned}
	\end{equation}
	
	We will prove that $\sqrt nV_{1n}^f$ is the leading term with a limiting distribution and $V_{2n}^f$ and $V_{3n}^f$ are $o_p(n^{-1/2}).$
	Decompose $V_{1n}^f$ as follows,
	\begin{equation*}
	\begin{aligned}
	V_{1n}^f=&\frac{1}{nh^2}\sum^{n}_{s=1} \left\{\left[\hat{X}^{\prime}_k(t_s)-f_k(t_s)\right]\hat{p}^2\left(t_{s}\right)\right\}^2 \\
	=& \frac{1}{nh^2}\sum^{n}_{s=1}\left[\hat{h}_k^{\prime}(t_d) \hat{p}(t_s)-\hat{h}_k(t_s) \hat{p}^{\prime}(t_s)-\hat{p}^2\left(t_{s}\right)f_k(t_s)\right]^2 \\
	=&\frac{1}{nh^2}\sum^{n}_{s=1} \left\{ \frac{1}{n^2}\sum^{n}_{i=1} \sum^{n}_{j=1} \left[\frac{1}{h^3} K^{\prime}\left(\frac{t_{s}-t_{i}}{h}\right)K\left(\frac{t_{s}-t_{j}}{h}\right) (Y_{ik}-Y_{jk}) \right. \right.\\
	&\left. \left. - \frac{1}{h^2}K\left(\frac{t_{s}-t_{i}}{h}\right)K\left(\frac{t_{s}-t_{j}}{h}\right)f_k(t_s)\right] \right\}^2\\
	=&\frac{1}{nh^2}\sum^{n}_{k=1} \left\{ \frac{1}{n^2}\sum^{n}_{i=1} \sum^{n}_{j=1}K\left(\frac{t_{s}-t_{j}}{h}\right) \left[\frac{1}{h^3} K^{\prime}\left(\frac{t_{s}-t_{i}}{h}\right) (Y_{ik}-Y_{jk}) \right. \right.\\
	&\left.\left. - \frac{1}{h^2}K\left(\frac{t_{s}-t_{i}}{h}\right)f_k(t_s)\right] \right\}^2\\
	=&\frac{1}{n^5h^2}\sum^{n}_{s=1}\sum^{n}_{a=1} \sum^{n}_{b=1}\sum^{n}_{c=1} \sum^{n}_{s=1}K\left(\frac{t_{s}-t_{a}}{h}\right)K\left(\frac{t_{s}-t_{b}}{h}\right) \left[\frac{1}{h^3} K^{\prime}\left(\frac{t_{s}-t_{c}}{h}\right) (Y_{ck}-Y_{ak}) \right.\\
	&\left. - \frac{1}{h^2}K\left(\frac{t_{s}-t_{c}}{h}\right)f_k(t_s)\right]\left[\frac{1}{h^3} K^{\prime}\left(\frac{t_{s}-t_{d}}{h}\right) (Y_{dk}-Y_{bk})  - \frac{1}{h^2}K\left(\frac{t_{s}-t_{d}}{h}\right)f_k(t_s)\right]\\
	\end{aligned}
	\end{equation*}
	
	Define
	\begin{equation*}
	\begin{aligned}
	&H_n^{\prime}(z_a,z_b,z_c,z_d,z_s)\\
	=&\frac{1}{h^2} K\left(\frac{t_{s}-t_{a}}{h}\right)K\left(\frac{t_{s}-t_{b}}{h}\right) \left[\frac{1}{h^3} K^{\prime}\left(\frac{t_{s}-t_{c}}{h}\right) (Y_{ck}-Y_{ak}) - \frac{1}{h^2}K\left(\frac{t_{s}-t_{c}}{h}\right)f_k(t_k)\right]\\
	&\left[\frac{1}{h^3} K^{\prime}\left(\frac{t_{s}-t_{d}}{h}\right) (Y_{dk}-Y_{bk})  - \frac{1}{h^2}K\left(\frac{t_{s}-t_{d}}{h}\right)f_k(t_s)\right]
	\end{aligned}
	\end{equation*}
	and define $H_n$ as the symmetry form of $H_n^{\prime}$.
	\begin{equation}
	H_n(z_a,z_b,z_c,z_d,z_s)=\frac{1}{5!}\sum_{P}H_n^{\prime}(z_a,z_b,z_c,z_d,z_s).\label{A27}
	\end{equation}
	Thus $V_{1n}^f$ is actually a V-statistic with the kernel $H_n$ of order $5$. Since
	\begin{equation*}
	\begin{aligned}
	&E\left[H_n^2(z_a,z_b,z_c,z_d,z_s)\right]\\
	\leq& E\left[H_n^{\prime 2}(z_a,z_b,z_c,z_d,z_s)\right]\\
	=&E\left\{ E \left\{\frac{1}{h^4} K^2\left(\frac{t_{s}-t_{a}}{h}\right)K^2\left(\frac{t_{s}-t_{b}}{h}\right) \right. \right.\\
	&\left[\frac{1}{h^3} K^{\prime}\left(\frac{t_{s}-t_{c}}{h}\right) (Y_{ck}-Y_{ak}) - \frac{1}{h^2}K\left(\frac{t_{s}-t_{c}}{h}\right)f_k(t_s)\right]^2\\
	&\left. \left. \left[\frac{1}{h^3} K^{\prime}\left(\frac{t_{s}-t_{d}}{h}\right) (Y_{dk}-Y_{bk})  - \frac{1}{h^2}K\left(\frac{t_{s}-t_{d}}{h}\right)f_k(t_s)\right]^2 | t_a,t_b,t_c,t_d,t_s \right\} \right\}\\
	\end{aligned}
	\end{equation*}
	\begin{equation*}
	\begin{aligned}
	=&E \left\{\frac{1}{h^4} K^2\left(\frac{t_{s}-t_{a}}{h}\right)K^2\left(\frac{t_{s}-t_{b}}{h}\right) \right.\\
	&\left[\frac{1}{h^{12}} K^{\prime 2}\left(\frac{t_{s}-t_{c}}{h}\right)K^{\prime 2}\left(\frac{t_{s}-t_{d}}{h}\right) (V_{ck}+V_{ak}-2X_{ak}X_{ck})(V_{dk}+V_{bk}-2X_{dk}X_{bk}) \right. \\
	&+ \frac{1}{h^8}K^2\left(\frac{t_{s}-t_{c}}{h}\right)K^2\left(\frac{t_{s}-t_{d}}{h}\right)f_k^4(t_s)\\
	&+ \frac{4}{h^{10}}K^{\prime}\left(\frac{t_{s}-t_{c}}{h}\right)K\left(\frac{t_{s}-t_{c}}{h}\right)K^{\prime}\left(\frac{t_{s}-t_{d}}{h}\right)K\left(\frac{t_{s}-t_{d}}{h}\right)f_k^2(t_s)  (X_{ck}-X_{ak})(X_{dk}-X_{bk})\\
	&+\frac{1}{h^{10}} K^{\prime 2}\left(\frac{t_{s}-t_{c}}{h}\right)K^{2}\left(\frac{t_{s}-t_{d}}{h}\right) (V_{ck}+V_{ak}-2X_{ak}X_{ck})f_k^2(t_s)\\
	&+\frac{1}{h^{10}} K^{\prime 2}\left(\frac{t_{s}-t_{d}}{h}\right)K^{2}\left(\frac{t_{s}-t_{c}}{h}\right) (V_{dk}+V_{bk}-2X_{dk}X_{bk})f_k^2(t_s)\\
	&-\frac{2}{h^{11}} K^{\prime 2}\left(\frac{t_{s}-t_{c}}{h}\right)K^{\prime}\left(\frac{t_{s}-t_{d}}{h}\right)K\left(\frac{t_{s}-t_{d}}{h}\right) (V_{ck}+V_{ak}-2X_{ak}X_{ck})(X_{dk}-X_{bk})f_k(t_s)\\
	&-\frac{2}{h^{11}} K^{\prime 2}\left(\frac{t_{s}-t_{d}}{h}\right)K^{\prime}\left(\frac{t_{s}-t_{c}}{h}\right)K\left(\frac{t_{s}-t_{c}}{h}\right) (X_{ck}-X_{ak})(V_{dk}+V_{bk}-2X_{dk}X_{bk})f_k(t_s)\\
	&-\frac{2}{h^{9}} K^{2}\left(\frac{t_{s}-t_{c}}{h}\right)K^{\prime}\left(\frac{t_{s}-t_{d}}{h}\right)K\left(\frac{t_{s}-t_{d}}{h}\right) (X_{dk}-X_{bk})f^3_k(t_s)\\
	&\left. \left.-\frac{2}{h^{9}} K^{2}\left(\frac{t_{s}-t_{d}}{h}\right)K^{\prime}\left(\frac{t_{s}-t_{c}}{h}\right)K\left(\frac{t_{s}-t_{c}}{h}\right) (X_{ck}-X_{ak})f^3_k(t_s) \right] \right\}\\
	=&\int K^2\left(u_a\right)K^2\left(u_b\right) \left[\frac{1}{h^{12}} K^{\prime 2}\left(u_c\right)K^{\prime 2}\left(u_d\right)M_1M_2 + \frac{1}{h^8}K^2\left(u_c\right)K^2\left(u_d\right)f^4_k(t_s)\right.\\
	&+  \frac{4}{h^{10}}K^{\prime}\left(u_c\right)K\left(u_c\right)K^{\prime}\left(u_d\right)K\left(u_d\right)f^2 _k(t_s) M_3M_4\\
	&+\frac{1}{h^{10}} K^{\prime 2}\left(u_c\right)K^{2}\left(u_d\right) M_1f^2_k(t_s)+\frac{1}{h^{10}} K^{\prime 2}\left(u_d\right)K^{2}\left(u_c\right) M_2f^2_k(t_s)\\
	&-\frac{2}{h^{11}} K^{\prime 2}\left(u_c\right)K^{\prime}\left(u_d\right)K\left(u_d\right) M_1M_4f(t_s)-\frac{2}{h^{11}} K^{\prime 2}\left(u_d\right)K^{\prime}\left(u_c\right)K\left(u_c\right) M_2M_3f(t_s)\\
	&\left. -\frac{2}{h^{9}} K^{2}\left(u_c\right)K^{\prime}\left(u_d\right)K\left(u_d\right) M_4f^3_k(t_s)-\frac{2}{h^{9}} K^{2}\left(u_d\right)K^{\prime}\left(u_c\right)K\left(u_c\right) M_3f^3_k(t_s) \right] \\
	&p(t_s-u_ah)p(t_s-u_bh)p(t_s-u_ch)p(t_s-u_dh)p(t_s)\mathrm{d}u_a\mathrm{d}u_b\mathrm{d}u_c\mathrm{d}u_d\mathrm{d}t_s\\
	=&O(\frac{1}{h^{12}})=o(n)\\
	\end{aligned}
	\end{equation*}
	where $M_1=V_k(t_s-u_ch)+V_k(t_s-u_ah)-2X_k(t_s-u_ah)X_k(t_s-u_ch)$ , $M_2=V_k(t_s-u_dh)+V_k(t_s-u_bh)-2X_k(t_s-u_dh)X_k(t_s-u_bh)$, $M_3=X_k(t_s-u_ch)-X_k(t_s-u_ah)$, $M_4=X_k(t_s-u_dh)-X_k(t_s-u_bh)$ and $V_i=V(t_i)=E(Y_i^2|t_i)$. The condition of Theorem 1 of \cite{MartinsFilho.2006} is satisfied. Thus we have
	\begin{equation}
	U_n-V_{1n}^f=o_P(n^{-1/2}) \label{A22}
	\end{equation}
	where $U_n$ is the corresponding U-statistic with the kernel $H_n$. Next we consider the limiting properties of $U_n$.  Application of Lemma~5 to $U_n(z_a,z_b,z_c,z_d,z_s)$ with some tedious computation, we can let the projection of $H_n(z_a,z_b,z_c,z_d,z_s)$, which can be computed as:
	\begin{equation*}
	\begin{aligned}
	&r_n(z_s)=E\left[H_n(z_a,z_b,z_c,z_d,z_s)|z_s\right]\\
	=&E\left[\frac{1}{5!}\sum_{P}H_n^{\prime}(z_a,z_b,z_c,z_d,z_s)|z_s\right]\\
	=&E\left\{\frac{1}{5h^2} K\left(\frac{t_{s}-t_{a}}{h}\right)K\left(\frac{t_{s}-t_{b}}{h}\right) \left[\frac{1}{h^3} K^{\prime}\left(\frac{t_{s}-t_{c}}{h}\right) (Y_{ck}-Y_{ak}) - \frac{1}{h^2}K\left(\frac{t_{s}-t_{c}}{h}\right)f_k(t_s)\right] \right.\\
	&\left[\frac{1}{h^3} K^{\prime}\left(\frac{t_{s}-t_{d}}{h}\right) (Y_{dk}-Y_{bk})  - \frac{1}{h^2}K\left(\frac{t_{s}-t_{d}}{h}\right)f_k(t_s)\right]\\
	&+\frac{2}{5h^2} K\left(\frac{t_{a}-t_{s}}{h}\right)K\left(\frac{t_{a}-t_{b}}{h}\right) \left[\frac{1}{h^3} K^{\prime}\left(\frac{t_{a}-t_{c}}{h}\right) (Y_{ck}-Y_{sk}) - \frac{1}{h^2}K\left(\frac{t_{a}-t_{c}}{h}\right)f_k(t_a)\right]\\
	&\left[\frac{1}{h^3} K^{\prime}\left(\frac{t_{a}-t_{d}}{h}\right) (Y_{dk}-Y_{bk})  - \frac{1}{h^2}K\left(\frac{t_{a}-t_{d}}{h}\right)f_k(t_a)\right]\\
	&+\frac{2}{5h^2} K\left(\frac{t_{c}-t_{a}}{h}\right)K\left(\frac{t_{c}-t_{b}}{h}\right) \left[\frac{1}{h^3} K^{\prime}\left(\frac{t_{c}-t_s}{h}\right) (Y_{sk}-Y_{ak}) - \frac{1}{h^2}K\left(\frac{t_{c}-t_s}{h}\right)f(t_c)\right]\\
	&\left[\frac{1}{h^3} K^{\prime}\left(\frac{t_{c}-t_{d}}{h}\right) (Y_{dk}-Y_{bk})  - \frac{1}{h^2}K\left(\frac{t_{c}-t_{d}}{h}\right)f(t_c)\right]|z_k\}\\
	=&\int \{ \frac{1}{5h^2} K\left(\frac{t_s-t_{a}}{h}\right)K\left(\frac{t_s-t_{b}}{h}\right) \left[\frac{1}{h^3} K^{\prime}\left(\frac{t_s-t_{c}}{h}\right) (X_{ck}-X_{ak}) - \frac{1}{h^2}K\left(\frac{t_s-t_{c}}{h}\right)f_k(t_s)\right]\\
	&\left[\frac{1}{h^3} K^{\prime}\left(\frac{t_s-t_{d}}{h}\right) (X_{dk}-X_{bk})  - \frac{1}{h^2}K\left(\frac{t_s-t_{d}}{h}\right)f_k(t_s)\right]\\
	&+\frac{2}{5h^2} K\left(\frac{t_{a}-t_s}{h}\right)K\left(\frac{t_{a}-t_{b}}{h}\right) \left[\frac{1}{h^3} K^{\prime}\left(\frac{t_{a}-t_{c}}{h}\right) (X_{ck}-Y_{sk}) - \frac{1}{h^2}K\left(\frac{t_{a}-t_{c}}{h}\right)f_k(t_a)\right]\\
	&\left[\frac{1}{h^3} K^{\prime}\left(\frac{t_{a}-t_{d}}{h}\right) (X_{dk}-X_{bk})  - \frac{1}{h^2}K\left(\frac{t_{a}-t_{d}}{h}\right)f_k(t_a)\right]\\
	&+\frac{2}{5h^2} K\left(\frac{t_{c}-t_{a}}{h}\right)K\left(\frac{t_{c}-t_{b}}{h}\right) \left[\frac{1}{h^3} K^{\prime}\left(\frac{t_{c}-t_s}{h}\right) (Y_{sk}-X_{ak}) - \frac{1}{h^2}K\left(\frac{t_{c}-t_s}{h}\right)f_k(t_c)\right]\\
	&\left.\left[\frac{1}{h^3} K^{\prime}\left(\frac{t_{c}-t_{d}}{h}\right) (X_{dk}-X_{bk})  - \frac{1}{h^2}K\left(\frac{t_{c}-t_{d}}{h}\right)f_k(t_c)\right] \right\} p(t_a)p(t_b)p(t_c)p(t_d)\mathrm{d}t_a\mathrm{d}t_b\mathrm{d}t_c\mathrm{d}t_d\\
	\end{aligned}
	\end{equation*}
	
	\begin{equation*}
	\begin{aligned}
	=&\int  \frac{1}{5h^2} K\left(u_a\right)K\left(u_b\right) \left[\frac{1}{h^3} K^{\prime}\left(u_c\right) (X_k(t_s-u_ch)-X_k(t_s-u_ah)) - \frac{1}{h^2}K\left(u_c\right)f_k(t_s)\right]\\
	&\left[\frac{1}{h^3} K^{\prime}\left(u_d\right) (X_k(t_s-u_dh)-X_k(t_s-u_bh))  - \frac{1}{h^2}K\left(u_d\right)f_k(t_s)\right]\\
	&+\frac{2}{5h^2} K\left(u_a\right)K\left(u_a-u_b\right) \left[\frac{1}{h^3} K^{\prime}\left(u_c-u_a\right) (X_k(t_s-u_ch)-Y_{sk}) - \frac{1}{h^2}K\left(u_c-u_a\right)f_k(t_s-u_ah)\right]\\
	&\left[\frac{1}{h^3} K^{\prime}\left(u_d-u_a\right) (X_k(t_s-u_dh)-X_k(t_s-u_bh))  - \frac{1}{h^2}K\left(u_d-u_a\right)f_k(t_s-u_ah)\right]\\
	&+\frac{2}{5h^2} K\left(u_a-u_c\right)K\left(u_b-u_c\right) \left[\frac{1}{h^3} K^{\prime}\left(-u_c\right) (Y_{sk}-X_k(t_s-u_ah)) - \frac{1}{h^2}K\left(u_c\right)f_k(t_s-u_ch)\right]\\
	&\left[\frac{1}{h^3} K^{\prime}\left(u_d-u_c\right) (X_k(t_s-u_dh)-X_k(t_s-u_bh))  - \frac{1}{h^2}K\left(u_d-u_c\right)f_k(t_s-u_ch)\right] \\ &p(t_s-u_ah)p(t_s-u_bh)p(t_s-u_ch)p(t_s-u_dh)h^4\mathrm{d}u_a\mathrm{d}u_b\mathrm{d}u_c\mathrm{d}u_d\\
	=&\frac{1}{15}\int u^3K^{\prime}(u)\mathrm{d}u X_k^{(4)}(t_s)\left[X_k(t_s)-Y_k(t_s)\right]p^4(t_s)\left[ 1+o_P(1)\right].\\
	\end{aligned}
	\end{equation*}
	The last equation uses  Taylor expansion and the properties on the kernel function: $\int K(u)\mathrm{d}u=1$, $\int u K(u)\mathrm{d}u=0$, $\int K^{\prime}(u)\mathrm{d}u=0$, $\int u K^{\prime}(u)\mathrm{d}u=-1$.
	
	Let $r(z_s)=\frac{1}{15}\int u^3K^{\prime}(u)\mathrm{d}u X_k^{(4)}(t_s)\left[X_k(t_s)-Y_k(t_s)\right]p^4(t_s)$. Further,
	\begin{equation*}
	E\left[r(z_k)\right]=0
	\end{equation*}
	and
	\begin{equation*}
	\begin{aligned}
	&Var\left[r(z_k)\right]\\
	=&E\left\{\frac{1}{15}\int u^3K^{\prime}(u)\mathrm{d}u X_k^{(4)}(t_s)\left[X_k(t_s)-Y_k(t_s)\right]p^4(t_s)\right\}^2\\
	=&\frac{1}{225}\left[\int u^3K^{\prime}(u)\mathrm{d}u\right]^2\int [X_k^{(4)}(t_s)]^2[V_k(t_s)-X_k^2(t_s)]p^8(t_s)\mathrm{d}t_s\\
	=&\frac{1}{225}\left[\int u^3K^{\prime}(u)\mathrm{d}u\right]^2\int [X_k^{(4)}(t_s)]^2\sigma_k^2(t_s) p^8(t_s)\mathrm{d}t_s .\\
	\end{aligned}
	\end{equation*}
	The limiting null distribution of $\sqrt{n}\left[U_n-E(U_n)\right]$ is the same as that of $\frac{5}{\sqrt{n}}\sum_{s=1}^{n}r(z_s)$. By the Lindeberg-Levy central limit theorem, we have
	\begin{equation}
	\sqrt{n}\left[U_n-E(U_n)\right] \stackrel{\mathbf{d}}{\longrightarrow} N\left(0,\frac{1}{9}\left[\int u^3K^{\prime}(u)\mathrm{d}u\right]^2\int [X_k^{(4)}(t_s)]^2\sigma_k^2(t_s) p^8(t_s)\mathrm{d}t_s\right).\label{A23}
	\end{equation}
	
	The limiting null distribution of $V_{1n}^f$ can then be derived by combining (\ref{A22}) and (\ref{A23}).
	
	We now prove that $V_{2n}^f=o_P(n^{-1/2})$. Decompose it as
	\begin{equation*}
	\begin{aligned}
	V_{2 n}^{f}=& \frac{1}{n h^{2}} \sum_{s=1}^{n}\left[\hat{X}_k^{\prime}\left(t_{s}\right)-f_k\left(t_{s}\right)\right]\left[\hat{f}_k\left(t_{s}\right)-f_k\left(t_{s}\right)\right]\hat{p}^{4}\left(t_{s}\right) \\
	=& \frac{1}{n h^{2}} \sum_{s=1}^{n}\left[\hat{X}_k^{\prime}\left(t_{s}\right)-f_k\left(t_{s}\right)\right]\left[f_k\left(t_{s},\hat{X}(t_s),\hat{\theta}\right)-f_k\left(t_{s},X(t_s),\hat{\theta}\right)\right]\hat{p}^{4}\left(t_{s}\right)\\
	&+\frac{1}{n h^{2}} \sum_{s=1}^{n}\left[\hat{X}_k^{\prime}\left(t_{s}\right)-f_k\left(t_{s}\right)\right]\left[f_k\left(t_{s},X(t_s),\hat{\theta}\right)-f_k\left(t_{s}\right)\right]\hat{p}^{4}\left(t_{s}\right) \\
	\equiv &S_{1}+S_{2}.
	\end{aligned}
	\end{equation*}
	For $S_1$, we have
	\begin{equation*}
	\begin{aligned}
	S_{1}=&\frac{1}{n h^{2}} \sum_{s=1}^{n}\left[\hat{X}_k^{\prime}\left(t_{s}\right)-f_k\left(t_{s}\right)\right]\left[f_k\left(t_{s},\hat{X}(t_s),\hat{\theta}\right)-f_k\left(t_{s},X(t_s),\hat{\theta}\right)\right]\hat{p}^{4}\left(t_{s}\right)\\
	=&\frac{1}{n h^{2}} \sum_{s=1}^{n}\left[\hat{h}_k^{\prime}\left(t_{s}\right) \hat{p}\left(t_{s}\right)-\hat{h}_k\left(t_{s}\right) \hat{p}^{\prime}\left(t_{s}\right)-\hat{p}^{2}\left(t_{s}\right) f_k\left(t_{s}\right)\right]\\
	&\times \left[f_k\left(t_{s},\hat{X}(t_s),\hat{\theta}\right)-f_k\left(t_{s},X(t_s),\hat{\theta}\right)\right]\hat{p}^{2}\left(t_{s}\right)\\
	=&\left\{ \frac{1}{n h^{2}} \sum_{s=1}^{n}\left[ \hat{h}_k^{\prime}\left(t_{s}\right) \hat{p}\left(t_{s}\right)-\hat{h}_k\left(t_{s}\right) \hat{p}^{\prime}\left(t_{s}\right)-\hat{p}^{2}\left(t_{s}\right) f_k\left(t_{s}\right)\right] \right.\\
	&\left.\times \frac{\partial f_k(t_s, X(t_s) ; \theta_0)}{\partial X^{\mathrm{T}}}\frac{N_{n}(t_s)-M_{n}(t_s)X(t_s)}{M_{n}(t_s)}\hat{p}^2(t_s)\right\}[1+o_P(1)]\\
	=&\left\{ \frac{1}{n h^{2}} \sum_{s=1}^{n}\left[ \hat{h}_k^{\prime}\left(t_{s}\right) \hat{p}\left(t_{s}\right)-\hat{h}_k\left(t_{s}\right) \hat{p}^{\prime}\left(t_{s}\right)-\hat{p}^{2}\left(t_{s}\right) f_k\left(t_{s}\right)\right] \right.\\
	&\left.\times \frac{\partial f_k(t_s, X(t_s) ; \theta_0)}{\partial X^{\mathrm{T}}}\left[\frac{N_{n}(t_s)-M_{n}(t_s)X(t_s)}{M_{n}(t_s)} -\frac{N_{n}(t_s)-M_{n}(t_s)X(t_s)}{M(t_s)}\right]\hat{p}^2(t_s)\right.\\
	&+\frac{1}{n h^{2}} \sum_{s=1}^{n}\left[ \hat{h}_k^{\prime}\left(t_{s}\right) \hat{p}\left(t_{s}\right)-\hat{h}_k\left(t_{s}\right) \hat{p}^{\prime}\left(t_{s}\right)-\hat{p}^{2}\left(t_{s}\right) f_k\left(t_{s}\right)\right] \\
	&\left.\times \frac{\partial f_k(t_s, X(t_s) ; \theta_0)}{\partial X^{\mathrm{T}}}\frac{N_{n}(t_s)-M_{n}(t_s)X(t_s)}{M(t_s)}\hat{p}^2(t_s) \right\}[1+\mathrm{o}_P(1)]\\
	\equiv& (S_{11}+S_{12})[1+\mathrm{o}_P(1)].
	\end{aligned}
	\end{equation*}
	Here the subscript $0$ represents that the corresponding estimator uses the bandwidth $h_0$ instead of $h$.
	
	By using Lemmas~\ref{LP.1}-\ref{LP.3}, we have
	\begin{equation*}
	\begin{aligned}
	S_{11}=& \frac{1}{n h^{2}} \sum_{s=1}^{n}\left[ \hat{h}_k^{\prime}\left(t_{s}\right) \hat{p}\left(t_{s}\right)-\hat{h}_k\left(t_{s}\right) \hat{p}^{\prime}\left(t_{s}\right)-\hat{p}^{2}\left(t_{s}\right) f_k\left(t_{s}\right)\right]\\
	&\times \frac{\partial f_k(t_s, X(t_s) ; \theta_0)}{\partial X^{\mathrm{T}}}\left[\frac{N_{n}(t_s)-M_{n}(t_s)X(t_s)}{M_{n}(t_s)} -\frac{N_{n}(t_s)-M_{n}(t_s)X(t_s)}{M(t_s)}\right]\hat{p}^2(t_s)\\
	=&\frac{1}{n h^{2}} \sum_{s=1}^{n}\left[ \hat{h}_k^{\prime}\left(t_{s}\right) \hat{p}\left(t_{s}\right)-\hat{h}_k\left(t_{s}\right) \hat{p}^{\prime}\left(t_{s}\right)-\hat{p}^{2}\left(t_{s}\right) f_k\left(t_{s}\right)\right]\\
	&\times \frac{\partial f_k(t_s, X(t_s) ; \theta_0)}{\partial X^{\mathrm{T}}}\hat{p}^2(t_s)\frac{M(t_s)-M_{n}(t_s)}{M_{n}(t_s)M(t_s)}[N_{n}(t_s)-M_{n}(t_s)X(t_s)]\\
	=&\frac{1}{n h^{2}} \sum_{s=1}^{n}\left[ \hat{X}_k^{\prime}(t_s)-X^{\prime}_k(t_s) \right]\frac{\partial f_k(t_s, X(t_s) ; \theta_0)}{\partial X^{\mathrm{T}}}\hat{p}^4(t_s)\frac{M(t_s)-M_{n}(t_s)}{M(t_s)}[\hat{X}(t_s)-X(t_s)]\\
	\leq & \frac{\Gamma_{\theta}}{h^{2}} \left[ \sup_t\left|\hat{X}_k^{\prime}(t)-X^{\prime}_k(t)\right|\right] \left[ \sup_t |M(t)-M_{n}(t)| \right] \left[\sup_t\left\|\hat{X}(t)-X(t)\right\|\right]\\
	=&O_P[a_n^2(h_0)]=o_P(n^{-1/2}).
	\end{aligned}
	\end{equation*}
	Since $S_{12}$ is also a V-statistic, a similar argument for proving $V_{1n}^f$ can yield $S_{12}=o_P(n^{-1/2})$ when $h_0 = o(h)$  Therefore we conclude that $S_1=o_P(n^{-1/2})$.
	
	Then we decompose $S_2$ as
	\begin{equation*}
	\begin{aligned}
	S_{2}=& \frac{1}{n h^{2}} \sum_{s=1}^{n}\left[ \hat{X}_{k}^{\prime}\left(t_{s}\right)-f_k\left(t_{s}\right)\right]\left[f_k\left(t_{s},X(t_s),\hat{\theta}\right)-f_k\left(t_{s}\right)\right]\hat{p}^{4}\left(t_{s}\right) \\
	=& \left\{ \frac{1}{n h^{2}} \sum_{s=1}^{n}\left[ \hat{h}_k^{\prime}\left(t_{s}\right) \hat{p}\left(t_{s}\right)-\hat{h}_k\left(t_{s}\right) \hat{p}^{\prime}\left(t_{s}\right)-\hat{p}^{2}\left(t_{s}\right) f_k\left(t_{s}\right)\right]\right.\\
	&\left.\times \frac{\partial f_k(t_s, X(t_s) ; \theta_0)}{\partial \theta}\hat{p}^2_{n}(t_s) \right\}(\hat{\theta}-\theta_0)[1+o_P(1)]\\
	\equiv &S_{21}(\hat{\theta}-\theta_0)[1+o_P(1)].\\
	\end{aligned}
	\end{equation*}
	Again, $S_{21}$ is a V-statistic at the rate of  $O_P(h^{-2}n^{-1/2})$. Thus we have
	\begin{equation*}
	S_{2}=S_{21}(\hat{\theta}-\theta_0)[1+o_P(1)]=O_P(h^{-2}n^{-1/2})O_P(n^{-1/2})=o_P(n^{-1/2}).
	\end{equation*}
	Altogether, $V_{2 n}^{f}=o_P(n^{-1/2}).$
	
	The remaining part is to show $V_{3 n}^{f}=o_p(n^{-1/2})$. Note that
	\begin{equation*}
	\begin{aligned}
	V_{3 n}^{f}=& \frac{1}{n h^{2}} \sum_{s=1}^{n}\left[\hat{f}_k(t_s)-f_k(t_s, X(t_s) ; \hat{\theta})+f_k(t_s, X(t_s) ; \hat{\theta})-f_k(t_s)\right]^{2}\hat{p}^{4}\left(t_s\right) \\
	=&\frac{1}{n h^{2}} \sum_{s=1}^{n}\left[f_k(t_s, \hat{X}(t_s) ; \hat{\theta})-f_k(t_s, X(t_s) ; \hat{\theta})\right]^{2}\hat{p}^{4}\left(t_s\right)\\
	&+\frac{1}{n h^{2}} \sum_{s=1}^{n}\left[f_k(t_s, X(t_s) ; \hat{\theta})-f_k(t_s)\right]^{2}\hat{p}^{4}\left(t_s\right)\\
	&+\frac{2}{n h^{2}} \sum_{s=1}^{n}\left[f_k(t_s, \hat{X}(t_s) ; \hat{\theta})-f_k(t_s, X(t_s) ; \hat{\theta})\right]\left[f_k(t_s, X(t_s) ; \hat{\theta})-f_k(t_s)\right]\hat{p}^{4}\left(t_s\right)\\
	\equiv & S_3+S_4+2S_5.\\
	\end{aligned}
	\end{equation*}
	
	It is easy to have
	\begin{equation*}
	\begin{aligned}
	S_3 \leq & \frac{1}{n h^{2}} \sum_{s=1}^{n}\left\|\frac{\partial f_k(t_s, X(t_s) ; \hat{\theta})}{\partial X^{\mathrm{T}}}\right\|^2\left[ \sup_t\left\|\hat{X}(t)-X(t)\right\|^{2}\right]\left[\sup_{t} p^{4}\left(t\right)\right]\left[1+o_P(1)\right]\\
	=&O_P[a_n^2(h_0)h^{-2}],
	\end{aligned}
	\end{equation*}
	and
	\begin{equation*}
	\begin{aligned}
	S_4 \leq & \frac{1}{n h^{2}} \sum_{s=1}^{n}\left\|\frac{\partial f_k(t_s, X(t_s) ; \theta_0)}{\partial \theta^{\mathrm{T}}}\right\|^2\left\| \hat{\theta}-\theta_0\right\|^2\left[\sup_{t} p^{4}\left(t\right)\right]\left[1+o_P(1)\right]\\
	=& O_P(n^{-1}h^{-2}).
	\end{aligned}
	\end{equation*}
	
	It is clear that $2|S_5|$ is bounded by $S_3+S_4$. Therefore we have
	\begin{equation*}
	V_{3 n}^{f}=O_P[a_n^2(h_0)h^{-2}]+O_P(n^{-1}h^{-2})=o_P(n^{-1/2}).
	\end{equation*}
	
	Summarizing the above results, we conclude
	\begin{equation*}
	\sqrt{n}\left(V_n^f-V\right) \stackrel{\mathbf{d}}{\longrightarrow} N(0,\Sigma^f).
	\end{equation*}
	
\end{proof}

\subsection{Proof of Lemma~8}

\begin{proof}
	The statistic $\hat{S}$ can be decomposed into three terms:
	\begin{equation}
	\begin{aligned}
	\hat{S}=&\frac{1}{h^2} \int \left[ f_k(t, \hat{X}(t) ; \hat{\theta})-\hat{X}_k^{\prime}(t) \right]^{2}\hat{p}^5(t) \mathrm{d} t\\
	=&\frac{1}{h^2} \int \left[ f_k(t, \hat{X}(t) ; \hat{\theta})-X_k^{\prime}(t)+X_k^{\prime}(t)-\hat{X}_k^{\prime}(t) \right]^{2}\hat{p}^5(t) \mathrm{d} t\\
	=&\frac{1}{h^2} \int\left[f_k(t, \hat{X}(t) ; \hat{\theta})-X_k^{\prime}(t)\right]^2\hat{p}^5(t) \mathrm{d} t\\
	&+\frac{1}{h^2}\int \left[X_k^{\prime}(t)-\hat{X}_k^{\prime}(t)\right]^{2}\hat{p}^5(t) \mathrm{d} t\\
	&+\frac{2}{h^2}\int \left[f_k(t, \hat{X}(t) ; \hat{\theta})-X_k^{\prime}(t)\right]\left[X_k^{\prime}(t)-\hat{X}_k^{\prime}(t)\right]\hat{p}^5(t)\mathrm{d} t\\
	\equiv &V_{1}+V_{2}+2V_{3}.\label{A24}\\
	\end{aligned}
	\end{equation}
	
	To prove that  $V_{1}=o_p(n^{-1/2})$, we  note that
	\begin{equation*}
	\begin{aligned}
	V_{1}=&\frac{1}{h^2} \int\left[ f_k(t, \hat{X}(t) ; \hat{\theta})-X_k^{\prime}(t)\right]^2 \hat{p}^5(t) \mathrm{d} t\\
	= & \frac{1}{h^2} \int\left[f_k(t, \hat{X}(t) ; \hat{\theta})-f_k(t, X(t) ; \hat{\theta})+f_k(t, X(t) ; \hat{\theta})-X_k^{\prime}(t)\right]^2 \hat{p}^5(t) \mathrm{d} t\\
	= &\frac{1}{h^2} \int\left[f_k(t, \hat{X}(t) ; \hat{\theta})-f_k(t, X(t) ; \hat{\theta})\right]^2 \hat{p}^5(t) \mathrm{d} t\\
	&+\frac{1}{h^2} \int\left[f_k(t, X(t) ; \hat{\theta})-X_k^{\prime}(t)\right]^2 \hat{p}^5(t) \mathrm{d} t\\
	&+\frac{2}{h^2} \int\left[f_k(t, \hat{X}(t) ; \hat{\theta})-f_k(t, X(t) ; \hat{\theta})\right]\left[f_k(t, X(t) ; \hat{\theta})-X_k^{\prime}(t)\right]\hat{p}^5(t)\mathrm{d} t\\
	\equiv &V_{11}+V_{12}+2V_{13}.\\
	\end{aligned}
	\end{equation*}
	
	Note that  $V_{11}$ is
	\begin{equation*}
	\begin{aligned}
	V_{11}=&\frac{1}{h^2} \int\left[f_k(t, \hat{X}(t) ; \hat{\theta})-f_k(t, X(t) ; \hat{\theta})\right]^2 \hat{p}^5(t) \mathrm{d} t\\
	=&\frac{1}{h^2} \int \left\{ \frac{\partial f_k(t, X(t) ; \hat{\theta})}{\partial X^{\mathrm{T}}}[\hat{X}(t)-X(t)][1+o_P(1)] \right\}^2 \hat{p}^5(t) \mathrm{d} t\\
	\leq &\frac{1}{h^2} \int \left\|\frac{\partial f_k(t, X(t) ; \hat{\theta})}{\partial X^{\mathrm{T}}}\right\|^2\left[ \sup_t\left\|\hat{X}(t)-X(t)\right\|^{2}\right] \hat{p}^5(t) \mathrm{d} t [1+o_P(1)]\\
	=&O_P[a_n^2(h_0)h^{-2}].
	\end{aligned}
	\end{equation*}
	For $V_{12}$, we have
	\begin{equation*}
	\begin{aligned}
	V_{12}=&\frac{1}{h^2} \int\left[f_k(t, X(t) ; \hat{\theta})-X_k^{\prime}(t)\right]^2\hat{p}^5(t) \mathrm{d} t\\
	=&\frac{1}{h^2} \int\left\{ \frac{\partial f_k(t, X(t) ; \theta)}{\partial \theta^{\mathrm{T}}}(\hat{\theta}-\theta_0)[1+o_P(1)]\right\}^2\hat{p}^5(t) \mathrm{d} t\\
	=&O_P(n^{-1}h^{-2}).
	\end{aligned}
	\end{equation*}
	Further, $2|V_{13}|$ is bounded by $V_{11}+V_{12}$. Therefore we have
	\begin{equation}
	V_{1}=O_P[a_n^2(h_0)h^{-2}]+O_P(n^{-1}h^{-2})=o_P(n^{-1/2}).\label{A25}
	\end{equation}
	For $V_2$, we have	
	\begin{equation}
	\begin{aligned}
	V_{2}=&\frac{1}{h^2}\int \left[ X_k^{\prime}(t)-\hat{X}_k^{\prime}(t)\right]^{2} \hat{p}^5(t)\mathrm{d} t\\
	\leq & \frac{1}{h^2} \int \left[ \sup |X_k^{\prime}(t)-\hat{X}_k^{\prime}(t)|\right]^2 \hat{p}^5(t) \mathrm{d} t\\
	\leq & (T-t_0) \left[\sup |X^{\prime}_k(t)-\hat{X}_k^{\prime}(t)|\right]^2\left[\sup_t \left|p(t)\right| \right]^5[1+o_P(1)]\\
	=&O_P[b_n^2(h_1)h^{-2}]=o_P(n^{-1/2}).\label{A26}
	\end{aligned}
	\end{equation}
	Again, $V_{3}$ is bounded by $V_{1}+V_{2}$. Thus we combine (\ref{A24})-(\ref{A26}) and conclude
	\begin{equation*}
	\hat{S}=o_P(n^{-1/2}).
	\end{equation*}

\end{proof}

\subsection{Proof of Lemma~9}

\begin{proof}
	To estimate $Var\left[r(z_s)\right]=E\left[r^2(z_s)\right]$, let
	\begin{equation*}
	\begin{aligned}
	&\hat{r}_n(z_s)=\frac{1}{n^*}\sum_{i=1}^{\lfloor \frac{n-1}{4} \rfloor} H_n(z_{1i},z_{2i},z_{3i},z_{4i},z_{s}),\\
	\end{aligned}
	\end{equation*}
	where $H_n$ is defined as (\ref{A27}) and $n^*=\lfloor \frac{n-1}{4} \rfloor$.
	
	Some elementary calculations yield that $E\left[\left\|\hat{w}_{n}\left(z_{s}\right)-\hat{r}_{n}\left(z_{s}\right)\right\|^{2}\right]=\mathrm{o}_P(1)$. Next we give the second moment consistency of  $\hat{r}_{n}\left(z_{s}\right)$  to $r(z_s)$ as the follows.
	Note that
	\begin{equation*}
	\begin{aligned} &E\left[\left\|\hat{r}_{n}\left(z_{s}\right)-r_{n}\left(z_{s}\right)\right\|^{2}\right]\\
	=&E\left[\operatorname{Var}\left(\hat{r}_{n}\left(z_{s}\right) | z_{s}\right)\right] \\
	=&\frac{1}{n^*} E\left[\operatorname{Var}\left(H_{n}\left(z_{1i},z_{2i},z_{3i},z_{4i},z_{s}\right) | z_{s}\right)\right] \\
	\leqslant & \frac{1}{n^*} E\left[\left\|H_{n}\left(z_{1i},z_{2i},z_{3i},z_{4i},z_{s}\right)\right\|^{2}\right]\\
	=&O\left(1 / (n h^{12})\right)=o(1),
	\end{aligned}
	\end{equation*}
	and $E\left[\left\|r_{n}\left(z_{s}\right)-r\left(z_{s}\right)\right\|^{2}\right]=O\left(h^{2}\right)=o(1)$ according to (\ref{A28}).
	
	Altogether, $\hat{w}_n^2(z_i)-r^2(z_i)=o_P(1)$. SLLN thus gives
	\begin{equation*}
	\begin{aligned}
	&\frac{1}{n-1} \sum_{i=1}^{n} \hat{w}_{n}^2\left(z_{s}\right) - (\frac{1}{n}\sum_{i}^{n}\hat{w}_{n}\left(z_{s}\right))^2 \\
	=&\frac{1}{n-1} \sum_{i=1}^{n} r^2\left(z_{s}\right) +o_{P}(1) \\ =&E\left[r^2\left(z_{s}\right)\right]+o_{P}(1).
	\end{aligned}
	\end{equation*}
\end{proof}

\subsection{Proof of Theorem~5.1 (under the null hypothesis)}

\begin{proof}
	The result is the consequence of Lemma~7, Lemma~8 and Lemma~9.
\end{proof}

\subsection{Proof of Theorem~5.2 (under global alternatives)}

\begin{proof}
	As shown in the proof of Lemma~7, we can decompose $V_n^f$ in (\ref{A21}) to three terms, in which, the second and third terms are asymptotically negligible and the limiting null distribution only depends on the first
	term $V^f_{1n}$. This term is also a V-statistic satisfying the condition of Theorem~1 of \cite{MartinsFilho.2006} and Lemma~5. Therefore we apply Lemma~5 to investigate the asymptotic properties. We compute the projection of the corresponding U-statistic as
	\begin{equation*}
	\begin{aligned}
	&r_n(z_s)=E\left[H_n(z_a,z_b,z_c,z_d,z_s)|z_s\right]\\
	=&E\left[\frac{1}{5!}\sum_{P}H_n^{\prime}(z_a,z_b,z_c,z_d,z_s)|z_s\right]\\
	=&\frac{1}{h^2}\left[f_k(t_s)-X_k^{\prime}(t_s)\right]^2p^4(t_s)+\frac{2}{5h^2}\left[f_k^{\prime}(t_s)-X_k^{(2)}(t_s)\right](Y_{sk}-X_{sk})p^4(t_s)+o_P(1/h^2)\\
	\equiv&r(z_s)+o_P(1/h^2).\\
	\end{aligned}
	\end{equation*}
	It is easy to see that
	\begin{equation*}
	E\left[r(z_s)\right]=\frac{1}{h^2}\int \left[f_k(t_s)-X_k^{\prime}(t_s)\right]^2p^4(t_s)\mathrm{d}t_s,
	\end{equation*}
	and
	\begin{equation*}
	\begin{aligned}
	&Var\left[r(z_s)\right]\\=&E\left[r^2(z_s)\right]-E\left[r(z_s)\right]^2\\
	=&E\{\frac{1}{h^4}\left[f_k(t_s)-X_k^{\prime}(t_s)\right]^4p^8(t_s)+\frac{4}{25h^4}\left[f_k^{\prime}(t_s)-X_k^{(2)}(t_s)\right]^2(Y_{sk}-X_{sk})^2p^8(t_s)\\
	&+\frac{2}{5h^4}\left[f_k(t_s)-X_k^{\prime}(t_s)\right]^2\left[f_k^{\prime}(t_s)-X_k^{(2)}(t_s)\right](Y_{sk}-X_{sk})p^8(t_s)\}\\
	&-\frac{1}{h^4}\left\{\int \left[f_k(t_s)-X_k^{\prime}(t_s)\right]^2p^4(t_s)\mathrm{d}t_s \right\}^2\\
	=&\int\{ \frac{1}{h^4}\left[f_k(t_s)-X_k^{\prime}(t_s)\right]^4p^8(t_s)+\frac{4}{25h^4}\left[f_k^{\prime}(t_s)-X_k^{(2)}(t_s)\right]^2 \sigma_k^2(t_s) p^8(t_s)\}\mathrm{d}t_s\\
	&-\frac{1}{h^4}\left\{\int \left[f_k(t_s)-X_k^{\prime}(t_s)\right]^2p^4(t_s)\mathrm{d}t_s \right\}^2\\
	= & \frac{1}{h^4}\Sigma_f^{\prime}.
	\end{aligned}
	\end{equation*}
	Then, recalling that the subscripts $\tilde{n}1$ and $\tilde{n}2$ mean the first and second subsample,
	\begin{equation*}
	\frac{\sqrt{\tilde{n}}h^2(V^{f}_{\tilde{n}1}-V^{f}_{\tilde{n}2})}{\sqrt{2\Sigma^f}} \stackrel{\mathbf{d}}{\longrightarrow} N(0,1).
	\end{equation*}
	By a similar proof of Lemma~9, $h^4\widehat{\Sigma}^f$ is consistent to $\Sigma_f^{\prime}$. Resembling to the proof of Lemma~8, we can derive $\hat{S} - \frac{1}{h^2}\int (f_k(t, X(t); \theta_1)-X_k^{\prime}(t))^2p(t)^5\mathrm{d}t=o_p(1)$. Then under the global alternatives
	\begin{equation*}
	\begin{aligned}
	GM/\sqrt{\tilde{n}}=&\frac{(V^{f}_{\tilde{n}1}-V^{f}_{\tilde{n}2}+c\hat{S})}{\sqrt{2\widehat{\Sigma}^f}}\\
	=&\frac{(h^2V^{f}_{\tilde{n}1}-h^2V^{f}_{\tilde{n}2}+ch^2\hat{S})}{\sqrt{2h^4\widehat{\Sigma}^f}}\\
	\stackrel{\mathbf{P}}{\longrightarrow}& \frac{c\int (f(t, X(t); \theta_1)-X^{\prime}(t))^2p(t)^5\mathrm{d}t}{\sqrt{2\Sigma^{f^{\prime}}}}.
	\end{aligned}
	\end{equation*}
	That is, $GM_{n}$ diverges to infinity at the rate of $\sqrt {n}$ in probability.
\end{proof}

\subsection{Proof of Theorem~5.3 (under local alternatives)}

\begin{proof}
	Here we give the proof under $H_{1kn}^f$. Due to the similarity, the proof under $H_{1kn}^F$ is omitted here.
	
	Similarly as the proof of Theorem~5.2, we can show
	\begin{equation*}
	\frac{\sqrt{\tilde{n}}(V^{f}_{\tilde{n}1}-V^{f}_{\tilde{n}2})}{\sqrt{2\widehat{\Sigma}^f}} \stackrel{\mathbf{d}}{\longrightarrow} N(0,1).
	\end{equation*}
	Recall
	\begin{equation}
	GM_{n}/(\tilde{n}^{1/2}h^{-2}\delta_n^2)=\frac{\sqrt{\tilde{n}}(V^{f}_{\tilde{n}1}-V^{f}_{\tilde{n}2})}{\tilde{n}^{1/2}\delta_n^2 h^{-2}\sqrt{2\widehat{\Sigma}^f}}+\frac{ch^2\hat{S}}{\delta_n^2\sqrt{2\widehat{\Sigma}^f}}.\label{A18}
	\end{equation}
	When $\tilde{n}^{1/4}h^{-1}\delta_n \rightarrow \infty$,
	\begin{equation}
	\frac{\sqrt{\tilde{n}}(V^{f}_{\tilde{n}1}-V^{f}_{\tilde{n}2})}{\tilde{n}^{1/2}\delta_n^2 h^{-2}\sqrt{2\widehat{\Sigma}^f}}=o_P(1).
	\end{equation}
	Now we compute the bias correction term $\hat{S}$. Under the local alternatives,
	\begin{equation*}
	\begin{aligned}
	\hat{S}=&\frac{1}{h^2} \int(f_k(t, \hat{X}(t) ; \hat{\theta})+\delta_n l_k(t)- \delta_n l_k(t) -\hat{X}_k^{\prime}(t))^{2} \mathrm{d} t\\
	=&\frac{1}{h^2} \int(f_k(t, \hat{X}(t) ; \hat{\theta})+\delta_n l_k(t) -\hat{X}_k^{\prime}(t))^{2} \mathrm{d} t\\
	+&\frac{1}{h^2} \int\delta_n^2 l_k^2(t) \mathrm{d} t-\frac{2}{h^2}\int \delta_n l_k(t) (f_k(t, \hat{X}(t) ; \hat{\theta})+\delta_n l_k(t) -\hat{X}_k^{\prime}(t))\mathrm{d} t\\
	\equiv &V_{1}+V_{2}-2V_{3}.\\
	\end{aligned}
	\end{equation*}
	Note that
	\begin{equation*}
	\begin{aligned}
	V_{1}=&\frac{1}{h^2} \int(f_k(t, \hat{X}(t) ; \hat{\theta})+\delta_n l_k(t) +X_k^{\prime}(t)-X_k^{\prime}(t) -\hat{X}_k^{\prime}(t))^{2} \mathrm{d} t\\
	=&\frac{1}{h^2} \int(f_k(t, \hat{X}(t) ; \hat{\theta})+\delta_n l_k(t)-X^{\prime}_k(t))^{2} \mathrm{d} t+\frac{1}{h^2} \int(X^{\prime}_k(t)-\hat{X}_k^{\prime}(t))^{2} \mathrm{d} t\\
	&+\frac{2}{h^2}\int (X_k^{\prime}(t)-\hat{X}_k^{\prime}(t))(f_k(t, \hat{X}(t) ; \hat{\theta})+\delta_n l_k(t)-X_k^{\prime}(t)) \mathrm{d} t\\
	\equiv &V_{11}+V_{12}+2V_{13}.
	\end{aligned}
	\end{equation*}
	Rewrite $V_{11}$ as
	\begin{equation*}
	\begin{aligned}
	V_{11}=&\frac{1}{h^2} \int(f_k(t, \hat{X}(t) ; \hat{\theta})+\delta_n l_k(t)-X_k^{\prime}(t))^2 \mathrm{d} t\\
	= & \frac{1}{h^2} \int(f_k(t, \hat{X}(t) ; \hat{\theta})+\delta_n l_k(t)-f_k(t, X(t) ; \hat{\theta})+f_k(t, X(t) ; \hat{\theta})-X_k^{\prime}(t))^2 \mathrm{d} t\\
	= &\frac{1}{h^2} \int(f_k(t, \hat{X}(t) ; \hat{\theta})-f_k(t, X(t) ; \hat{\theta}))^2 \mathrm{d} t\\
	&+\frac{1}{h^2} \int(f_k(t, X(t) ; \hat{\theta})+\delta_n l_k(t)-X_k^{\prime}(t))^2 \mathrm{d} t\\
	&+\frac{2}{h^2} \int(f_k(t, \hat{X}(t) ; \hat{\theta})-f_k(t, X(t) ; \hat{\theta}))(f_k(t, X(t) ; \hat{\theta})+\delta_n l_k(t)-X_k^{\prime}(t))\mathrm{d} t\\
	\equiv &V_{111}+V_{112}+V_{113}.\\
	\end{aligned}
	\end{equation*}
	
	By Taylor expansion, we can show $V_{111}$ is negligible at the rate $O_P[a_n^2(h_0^2)h^{-2}]=o_P(n^{-1/2})$.  According to Theorem~3.2, we have that $V_{112}$ is not negligible:
	\begin{equation}
	\begin{aligned}
	V_{112}=&\frac{1}{h^2} \int\left[ f_k(t, X(t) ; \hat{\theta})+\delta_n l_k(t)-X_k^{\prime}(t)\right]^2 \mathrm{d} t\\
	=&\frac{1}{h^2} \int \left\{ \frac{\partial f_k[t, X(t) ; \theta]}{\partial \theta}(\hat{\theta}-\theta_0)[1+o_P(1)]\right\}^2 \mathrm{d} t\\
	=&\frac{\delta_n^2}{h^2}\left\{H_{\dot{f}}^{-1} E_{p^*}\left[l_k\left(t\right) \omega_k\left(t\right) \frac{\partial f_{k}\left(t,X\left(t\right), \theta_{0}\right)}{\partial \theta}\right]\right\}^{\mathrm{T}}\\
	&\times \left[ \int \frac{\partial f_k(t, X(t) ; \theta_0)}{\partial \theta} \frac{\partial f_k(t, X(t) ; \theta_0)}{\partial \theta^{\mathrm{T}}} \mathrm{d} t\right]\\
	&\times \left\{H_{\dot{f}}^{-1} E_{p^*}\left[l_k\left(t\right) \omega_k\left(t\right) \frac{\partial f_{k}\left(t,X\left(t\right), \theta_{0}\right)}{\partial \theta}\right]\right\}+o_P(1).\label{A19}
	\end{aligned}
	\end{equation}
	An analogous calculation yields that $V_{113}$ is negligible as
	\begin{equation*}
	V_{113}=O_P[a_n(h_0)\delta_n h^{-2}]= o_P(\delta_n^2 h^{-2}).
	\end{equation*}
	To prove that $V_{12}$ and $V_{13}$ are negligible. We have
	\begin{equation*}
	\begin{aligned}
	V_{12}=&\frac{1}{h^2}\int (X_k^{\prime}(t)-\hat{X}_k^{\prime}(t))^{2} \mathrm{d} t\\
	=&O_P[b_n^2(h_1)h^{-2}]=o_P(n^{-1/2}),
	\end{aligned}
	\end{equation*}
	\begin{equation*}
	\begin{aligned}
	V_{13}=&\frac{1}{h^2}\int (X_k^{\prime}(t)-\hat{X}_k^{\prime}(t))(f_k(t, \hat{X}(t) ; \hat{\theta})+\delta_n l_k(t)-X_k^{\prime}(t)) \mathrm{d} t\\
	=&\frac{1}{h^2}\int (X_k^{\prime}(t)-\hat{X}_k^{\prime}(t))(f_k(t, \hat{X}(t) ; \hat{\theta}) -f_k(t, X(t) ; \theta_0)) \mathrm{d} t\\
	=&\frac{1}{h^2}\int (X_k^{\prime}(t)-\hat{X}_k^{\prime}(t))(f_k(t, X(t) ; \hat{\theta})-f_k(t, X(t) ; \theta_0)) \mathrm{d} t\\
	&+\frac{1}{h^2}\int (X_k^{\prime}(t)-\hat{X}_k^{\prime}(t))(f_k(t, \hat{X}(t) ; \hat{\theta})-f_k(t, X(t) ; \hat{\theta})) \mathrm{d} t\\
	=&O_P(b_n(h_1) \delta_n h^{-2})+O_P(b_n(h_1) a_n(h_0) h^{-2}) =o_P(\delta_n^2 h^{-2}).
	\end{aligned}
	\end{equation*}
	Altogether, in $V_1$, only  $V_{112}$ is a non-negligible term.
	Turn to $V_{2}$. We have
	\begin{equation}
	V_2=\frac{\delta_n^2}{h^2} \int l_k^2(t) \mathrm{d} t.\label{A40}
	\end{equation}
	which is also a non-negligible term.
	
	Finally, we focus on $V_{3}$ to derive that
	\begin{equation}
	\begin{aligned}
	V_{3}=&\frac{1}{h^2}\int \delta_n l_k(t) (f_k(t, \hat{X}(t) ; \hat{\theta})+\delta_n l_k(t) -\hat{X}^{\prime}(t))\mathrm{d} t\\
	=&\frac{1}{h^2}\int \delta_n l_k(t) (f_k(t, \hat{X}(t) ; \hat{\theta})+\delta_n l_k(t)-X_k^{\prime}(t)+X_k^{\prime}(t) -\hat{X}_k^{\prime}(t))\mathrm{d} t\\
	=&\frac{1}{h^2}\int \delta_n l_k(t) (f_k(t, X(t) ; \hat{\theta})-f_k(t, X(t) ; \theta_0))\mathrm{d} t[1+o_P(1)]\\
	=&\frac{\delta_n^2}{h^2} \left[\int l_k(t) \frac{\partial f_k(t, X(t) ; \theta_0)}{\partial \theta^{\mathrm{T}}} \mathrm{d} t \right] H_{\dot{f}}^{-1} E_{p^*}\left[ l_k\left(t\right) \omega_k\left(t\right) \frac{\partial f_{k}\left(t,X\left(t\right), \theta_{0}\right)}{\partial \theta}\right]\\
	&+o_P(1).\label{A20}
	\end{aligned}
	\end{equation}
	
	Summarizing the above results, the leading term of $\hat{S}$ is  the sum $V_{112}+V_{2}-2V_{3}$. Combining (\ref{A18})-(\ref{A20}), we can show that
	\begin{equation*}
	\begin{aligned}
	GM_{n}/(\tilde{n}^{1/2}h^{-2}\delta_n^2)=&\frac{\sqrt{\tilde{n}}(V^{f}_{\tilde{n}1}-V^{f}_{\tilde{n}2})}{\tilde{n}^{1/2}\delta_n^2 h^{-2}\sqrt{2\widehat{\Sigma}^f}}+\frac{ch^2\hat{S}}{\delta_n^2\sqrt{2\widehat{\Sigma}^f}}\\
	=&\frac{ch^2(V_{112}+V_2-2V_3)}{\delta_n^2\sqrt{2\widehat{\Sigma}^f}}+o_P(1)\\
	=&c\mu_{4}/ \sqrt{2\Sigma^f}+o_P(1)
	\end{aligned}	
	\end{equation*}
	where
	\begin{equation}
	\begin{aligned}
	\mu_{4}=&\left\{H_{\dot{f}}^{-1} E_{p^*}\left[l_k\left(t\right) \omega_k\left(t\right) \frac{\partial f_{k}\left(t,X\left(t\right), \theta_{0}\right)}{\partial \theta}\right]\right\}^{\mathrm{T}}\\
	&\times \left[ \int \frac{\partial f_k(t, X(t) ; \theta_0)}{\partial \theta} \frac{\partial f_k(t, X(t) ; \theta_0)}{\partial \theta^{\mathrm{T}}} \mathrm{d} t\right]\\
	&\times \left\{H_{\dot{f}}^{-1} E_{p^*}\left[l_k\left(t\right) \omega_k\left(t\right) \frac{\partial f_{k}\left(t,X\left(t\right), \theta_{0}\right)}{\partial \theta}\right]\right\}+\int l_k^2(t) \mathrm{d} t\\
	&-2\left[\int l_k(t) \frac{\partial f_k(t, X(t) ; \theta_0)}{\partial \theta^{\mathrm{T}}} \mathrm{d} t \right] H_{\dot{f}}^{-1} E_{p^*}\left[ l_k\left(t\right) \omega_k\left(t\right) \frac{\partial f_{k}\left(t,X\left(t\right), \theta_{0}\right)}{\partial \theta}\right].
	\end{aligned}
	\end{equation}
	
	Similarly we can easily prove the result (5.15) of Theorem~5.3 under the condition that $\delta_{n}=\tilde{n}^{-1 / 4} h$.
\end{proof}

\bibliographystyle{apalike}
\bibliography{document}